\documentclass[a4paper,11pt,twoside]{article}
\usepackage{pst-fill,pst-grad} 
\usepackage{textcomp}
\usepackage[english]{babel}
\usepackage[utf8]{inputenc} 
\usepackage[titletoc]{appendix}
\usepackage{titlesec}
\usepackage{graphicx}
\usepackage{amsmath}
\usepackage{float} 
\usepackage{fancyhdr}  
\usepackage[matrix,arrow,curve]{xy}
\usepackage{pstricks} 
\usepackage{amsmath,amsfonts,verbatim,afterpage,theorem,euscript,mathrsfs,amssymb}
\usepackage{amsfonts}
\usepackage{amssymb}
\usepackage{array}
\usepackage{dsfont}
\usepackage[colorlinks=true,linkcolor=blue,citecolor=red]{hyperref}
\usepackage{authblk}
\usepackage{color} 
\newtheorem{Definition}{Definition}[section] 
\newtheorem{DefinitionP}{Definition}
\newtheorem{Proposition}{Proposition}[section]
\newtheorem{PropositionP}{Proposition}
\newtheorem{Lemme}{Lemma}[section]
\newtheorem{Theoreme}{Theorem}[section]
\newtheorem{TheoremeP}{Theorem}

\def \fe{\vec{f}} 
\def \vg{\vec{g}}
\def \vu{\vec{u}}
\def \vv{\vec{v}}
\def \vw{\vec{w}}
\def \P{\mathbb{P}}
\def \U{\vec{U}}
\def \W{\vec{W}}
\def \V{\vec{V}}
\def \R{\mathbb{R}}

\def \Rt{\mathbb{R}^3}

\def \finpv{\hfill $\blacksquare$  \\ \newline }
\def \pv{{\bf{Proof.}}~} 

\def \ds{\displaystyle}

\title{\bf On the long-time behavior for a damped Navier-Stokes-Bardina model}

\author[1]{ Manuel Fernando Cortez\footnote{ manuel.cortez@epn.edu.ec}}
\author[2]{ Oscar Jarr\'in\footnote{corresponding author: oscar.jarrin@udla.edu.ec}} 

\affil[1]{\scriptsize Departamento de Matem\'aticas, Escuela Politécnica Nacional,  Ladron de Guevera E11-253, Quito, Ecuador}
\affil[2]{\scriptsize Dirección General  de Vinculaci\'on e  Investigación  (DGVI),
	Universidad de las Américas,
	Calle José Queri s/n y Av. Granados. Bloque 7, Tercer Piso, Quito, Ecuador} 
\date{\today}
\begin{document} 
\maketitle
\tableofcontents 
%%%%%%%%%%%%%%%%%%%%%%%%%%%%%%%%%%%%%%%%%%%%%%
\begin{abstract}
 In this paper,  we consider a damped Navier-Stokes-Bardina model posed on the whole three-dimensional. These equations have an important  physical motivation and they arise from some oceanic model.  From the mathematical point of view, they write down as the  well-know Navier-Stokes equations with an additional  nonlocal operator in their nonlinear transport term, and moreover, with an additional damping term depending of a parameter $\beta>0$. We study first the existence and \emph{uniqueness} of global in time weak solutions in the \emph{energy space}. Thereafter, our main objective is to describe  the \emph{long time behavior} of these solutions. For this,  we use some tools in the theory of   dynamical systems to prove the existence of a \emph{global attractor}, which is \emph{compact subset} in the energy space attracting all the weak solutions when the time goes to infinity. Moreover, we derive an upper bound for the  \emph{fractal dimension} of the global attractor associated to these equations.  

Finally, we find a range of  values for the  damping parameter $\beta>0$, where we are  able to give an acutely description of the internal structure of the global attractor. More precisely,  we prove that the global attractor only contains the stationary (time-independing) solution of the  damped Navier-Stokes-Bardina equations.   \\[3mm]  
\end{abstract} 
%\tableofcontents
%%%%%%%%%%%%%%%%%%%%%%%%%%%%%%%%%%%%%%%%%%%%%%
%%%%%%%%%%%%%%%%%%%%%%%%%%%%%%%%%%%%%%%%%%%%%%
\section{Introduction} 
The theory of partial differential equations (PDEs) is a broad research field, rapidly growing in close connections with other mathematical disciplines and applied sciences. The connections between the theories of dynamical systems and PDEs will be explored from several points of view. Infinite-dimensional dynamical systems generated by evolutionary PDEs provide the most immediate examples of interplay between the two theories. Extensions of well established results and techniques from finite-dimensional dynamical systems (invariant manifolds and bifurcations) have proved very useful in qualitative studies of PDEs. On the other hand, specific questions for PDEs brought about stimulating problems in the theory of dynamical systems, such as the existence of finite-dimensional attractors and their behavior under (regular or singular) perturbations. Entire (or eternal) solutions, which emerged as key objects in these problems, have long served as organizing centers for qualitative investigations of dissipative evolutionary PDEs and they continue to play an important role in other modern approaches to PDEs. 

\medskip

It is well-known that the  dynamic of an incompressible fluid, which we will assume on the whole space $\Rt$, is successfully described by the   classical, homogeneous and  incompressible Navier-Stokes equations: 
\begin{equation}\label{NS}
\partial_t \vv + \text{div}(\vv \otimes \vv) - \nu \Delta \vv + \vec{\nabla} q =0, \quad  \text{div}(\vv)=0, 
\end{equation}
where, $\vv : [0, +\infty[\times \Rt \to \Rt$ and   $q: [0, +\infty[\times \Rt \to \R$ are the velocity of the fluid and the pressure term respectively, the parameter $\nu>0$   represents the kinetic viscosity parameter which we will  keep fix. Moreover, the equation $\text{div}(\vv)=0$ describes the  fluid's incompressibility. 

\medskip 

Although the Navier-Stokes equations are a relevant physical model used in many applications, see for instance the book \cite{Temam1}, the mathematical theory of these equations is not yet sufficient to prove the global well-posedness of the so-called Leray's solutions, and the uniqueness of these solutions remains a very challenging open question. In order to contour this problem, researchers who are investigating the use  of the Navier-Stokes equations in practical applications \cite{Bardina1} have applied some operators to these equations, to obtain regularized versions of the Navier-Stokes equations where the weak solutions are well-posed.  

\medskip 

The idea is to introduce a regularized velocity field $\vu(t,x)$ in terms of the original velocity field $\vv(t,x)$ given by the equation (\ref{NS}) leading a variety of  so-called $\alpha-$ models for the Navier-Stokes equations. See for instance the Chapter $17$ of the book \cite{PGLR1}.   In this context, for a parameter $\alpha>0$, J. Bardina, J. H. Ferziger, \& W. C. Reynolds introduced in \cite{Bardina}  the  operator $(\cdot)_\alpha$,  given by solving the  Helmholtz equation: 
 
\begin{equation*}
-\alpha^2 \Delta (\varphi)_\alpha + (\varphi)_\alpha =\varphi. 
\end{equation*}

The operator $(\cdot)_{\alpha}$ is  also called the \emph{filtering/averaging} operator, due to the fact that this operator allows us to obtain an accurate  model describing the large-scale motion of the fluid while filtering or averaging  the fluid motion at scales smaller than $\alpha$.

\medskip 

On the other hand, from the mathematical point of view, we may observe in the whole space $\Rt$, and denoting by $I_d$ the identity operator,  the expression $(\varphi)_\alpha$ is given by 
\begin{equation}\label{filtering} 
(\varphi)_\alpha= (-\alpha^2 \Delta + I_d)^{-1} \varphi,
\end{equation} where the action of the operator $\ds{(-\alpha^2 \Delta + I_d)^{-1} }$, also known as the \emph{Bessel Potential}, could be easily defined in the Fourier variable as $\ds{\mathcal{F}\left((-\alpha^2 \Delta + I_d)^{-1} \varphi\right)(\xi)=(\alpha^2 \vert \xi \vert^2+1)^{-1}\widehat{\varphi}(\xi)}$ for all $\varphi \in \mathcal{S}(\Rt)$. For a more exhaustive study of \emph{Bessel Potentials} see the Chapter $6$ of the book \cite{Grafakos}.  

\medskip 

Applying the  operator $(\cdot)_\alpha$ to the Navier-Stokes equations (\ref{NS}),  we get the following equations for the regularized velocity $(\vv)_{\alpha}$, and the regularized pressure $(q)_\alpha$:   
\begin{equation*}
\partial_t (\vv)_{\alpha} + \text{div}((\vv \otimes \vv))_{\alpha} - \nu \Delta (\vv)_\alpha + \vec{\nabla} (q)_{\alpha} = 0, \quad  \text{div}((\vv)_\alpha)=0.
\end{equation*}

%From the mathematical point of view, we observe that the filtering  operator $(\cdot)_\alpha$ given in (\ref{filtering}) regularizes the functions in the framework of the non homogeneous Sobolev spaces $H^s(\Rt)$ (with $s\in \mathbb{R}$)  and this fact  allow us to prove the \emph{uniqueness of weak solutions in the energy space}, as we will see later (in Theorem \ref{Th-WP} and the comments below) in detail. \\

However, there is still a  problem to overcome. These regularized Navier-Stokes equations  are not a closed system in the sense that, due to the regularized nonlinear term: $\ds{\text{div}((\vv \otimes \vv))_{\alpha}}$, the system  does not write down only in terms of the regularized velocity $(\vv)_\alpha$ and the regularized pressure $(q)_\alpha$. More precisely, we have the identity 
\begin{equation*}\label{non-lin-term}
\text{div}((\vv \otimes \vv))_{\alpha} = \text{div} \left( (\vv)_\alpha \otimes (\vv)_{a} \right)+ \text{div} (R(\vv,\vv)),     
 \end{equation*}
where the remainder   term 
\begin{equation*}
 R(\vv,\vv)= (\vv \otimes \vv)_{\alpha}- (\vv)_{\alpha}\otimes (\vv)_{\alpha},   
\end{equation*}
is known as the Reynolds stress tensor. See, \emph{e.g.}, \cite{Cao} for more details on this term. 

\medskip 

In order to obtain a closed system, W. Layton \& R. Lewandowski propose in  \cite{Layton1} to approximate the Reynolds stress tensor as follows:
\begin{equation*}
R(\vv,\vv) \approx \left( (\vv)_{\alpha} \otimes (\vv)_\alpha \right)_{\alpha}  -  (\vv)_{\alpha} \otimes (\vv)_\alpha.    
\end{equation*}
Hence, replacing this approximation of the term $R(\vv,\vv)$ in the identity for the nonlinear term above, we obtain  the following approximation for the nonlinear term:
\begin{equation*}
\text{div}((\vv \otimes \vv))_{\alpha} \approx  \text{div}\left(\left(  (\vv)_\alpha \otimes (\vv)_\alpha\right)_{\alpha}\right).    
\end{equation*}
%revisado

This approximation of the nonlinear term  has been successfully used in many practical applications   \cite{Cao,Chow};  and it finally  leads us to  the following  closed system: 
\begin{equation*}
\partial_t (\vu)_{\alpha} + \text{div}\left(\left(  (\vv)_\alpha \otimes (\vv)_\alpha\right)_{\alpha}\right) - \nu \Delta (\vu)_\alpha + \vec{\nabla} (p)_{\alpha} = 0.  \quad  \text{div}((\vu)_\alpha)=0. 
\end{equation*}
To simplify the notation, we shall write the regularized functions $\ds{(\vv)_\alpha}$ and $ \ds{(q)_\alpha}$ as $\vu$ and $p$ respectively; and we  thus obtain  the so-called  \emph{Navier-Stokes-Bardina} model: 

\begin{equation}\label{Bardina-classical}
\partial_t \vu + \text{div}((\vu \otimes \vu)_\alpha) - \nu \Delta \vu + \vec{\nabla} p =0, \quad  \text{div}(\vu)=0, \quad  \vu(0,\cdot)=\vu_0,   \quad \alpha>0, \,\,  \nu>0.
\end{equation}
where the initial $\vu_0$ denotes the (regularized) velocity of the fluid at the time $t=0$. 

\medskip 

The first mathematical studies for these equations  were given by  W. Layton \& R. Lewandowski in  \cite{Layton2} in the space-periodic case.  In this setting, for an initial datum $\vu_0 \in H^1(\mathbb{T}^2)$, where for $L>0$, $\mathbb{T}^3$ denotes the periodic box $[0, 2\pi L]^{3}$, and  for all time $T>0$, it is proven the  existence and the \emph{uniqueness of a weak solution} $(\vu, p)$ which verifies $ \vu \in L^{\infty}([0,T[, H^{1}(\mathbb{T}^3)) \cap L^{2}([0,T[, \dot{H}^{2}(\mathbb{T}^3))$. More recently, using  a variational formulation of the equation (\ref{Bardina-classical}),  and some \emph{a priori} energy estimates,  L. C. Berselli \& R. Lewandowski extended in \cite{Berselli} the previous well-posedness results ( given for the periodic case) to the whole space $\Rt$. See also the Chapter $17.5$ of the book \cite{PGLR1} for an alternative proof of this result. Moreover, in the same work \cite{Berselli},  the regularity of weak solutions is also studied and it is proven that, for all time $t>0$, the unique weak solution  $(\vu(t,\cdot), p(t,\cdot))$ of the Navier-Stokes-Bardina equations (\ref{Bardina-classical}) belongs to the Sobolev space $H^{m}(\Rt)$ for all entire $m \geq 0$.   

\medskip 

On the other hand, another relevant issue  for the Navier-Stokes Bardina's model is the  \emph{asymptotic properties} of weak solutions when times goes to infinity. From a physical perspective,   this problem is also interesting when the Navier-Stokes Bardina's model is used to perform numerical simulations related to the turbulence description \cite{Bardina1}. The main idea is to consider an external force term $\fe$ in the equations Navier-Stokes-Bardina's model, which is assumed to be  a stationary (time-independing) vector field.  This stationarity assumption is a simplification of the physical  model. Indeed, the idea behind of this physical  model is that we will assume that a time independing  external source  acts on the fluid and put its in a perpetual turbulent state. In this scenario, we are thus interested in understanding the behavior of the velocity $\vu(t,\cdot)$ when the time $t$ is large enough. It is also worth mention that if we consider a time-dependent force then we will need to find an appropriate  time interval in which the fluid is turbulent and this is a highly non-trivial issue. See  the Chapter $1$ of \cite{Jarrin} for a more detailed discussion. Thus, for a  stationary external force $\fe:\Rt \to \Rt$,  we obtain the following forced equations:

\begin{equation*}
\partial_t \vu + \text{div}((\vu \otimes \vu)_\alpha) - \nu \Delta \vu + \vec{\nabla} p =\fe, \quad  \text{div}(\vu)=0,  \quad  \vu(0,\cdot)=\vu_0, \quad \alpha>0, \,\,  \nu>0.
\end{equation*}

In the space-periodic setting, the long-time behavior of solutions was studied in \cite{Cao}. In that work, using the \emph{Poincaré} inequality the authors obtain good controls on the quantity $\Vert \vu(t,\cdot)\Vert_{H^1(\mathbb{T}^3)}$ when the time $t$ goes to infinity. More precisely, for a constant $\eta>0$, and for all time $t\geq 0$, they get the following estimate:  

\begin{equation*}
\Vert \vu(t,\cdot) \Vert^{2}_{L^2(\mathbb{T}^3)}  +  \alpha\, \Vert \vu (t,\cdot)\Vert^{2}_{\dot{H}^1(\mathbb{T}^3)} \lesssim  e^{-\eta\, t} \left( \Vert \vu_0 \Vert^{2}_{L^2(\mathbb{T}^3)}+\alpha \Vert \vu_0  \Vert^{2}_{\dot{H}^1(\mathbb{T}^3)} \right) + \Vert \fe \Vert^{2}_{L^2(\mathbb{T}^3)}+\alpha \Vert \fe  \Vert^{2}_{\dot{H}^1(\mathbb{T}^3)}. 
\end{equation*}
Hence, for $t> 0$ large enough, we may observe that the norm $\Vert \vu(t,\cdot)\Vert_{H^1(\mathbb{T}^3)}$, expressed by the \emph{equivalent} quantity  $\ds{\Vert \vu(t,\cdot) \Vert^{2}_{L^2(\mathbb{T}^3)}  +  \alpha\, \Vert \vu (t,\cdot)\Vert^{2}_{\dot{H}^1(\mathbb{T}^3)}}$, is controlled \emph{uniformly in time}   by the quantity involving the external force: $\ds{\Vert \fe \Vert^{2}_{L^2(\mathbb{T}^3)}+\alpha \Vert \fe  \Vert^{2}_{\dot{H}^1(\mathbb{T}^3)}}$; and  this control is  one of the key tools to apply the classical theory of \emph{dynamical systems} to study the long time behavior of solutions $\vu(t,\cdot)$. Specifically, the authors prove the existence of a \emph{global attractor} for the Navier-Stokes-Bardina equations with space-periodic conditions. In Section \ref{Sec:Results} below,  we recall the definition of a \emph{global attractor} and introduce  more detail all the tools in the theory  of  dynamical systems used to perform this study.

\medskip 

Now, in the non-periodic setting of the whole space $\Rt$, and due to the lack of the \emph{Poincaré} inequality, we can only obtain the following  not so useful estimate in time  (see the details in  Appendix ()) 

\begin{equation*}
 \Vert \vu(t,\cdot) \Vert^{2}_{L^2(\mathbb{R}^3)}  +  \alpha\, \Vert \vu (t,\cdot)\Vert^{2}_{\dot{H}^1(\mathbb{R}^3)} \lesssim   \Vert \vu_0 \Vert^{2}_{L^2(\mathbb{R}^3)}+\alpha \Vert \vu_0  \Vert^{2}_{\dot{H}^1(\mathbb{R}^3)}  + t \left( \Vert \fe \Vert^{2}_{L^2(\mathbb{R}^3)}+\alpha \Vert \fe  \Vert^{2}_{\dot{H}^1(\mathbb{R}^3)} \right),  
\end{equation*}
and here we  clearly loose any control in time when $t$ goes to infinity.  To contour this problematic,  some previous works related to the study of the long-time behavior for the Navier-Stokes equations and related models  \cite{Constantin2,Ilyn,Jarrin} suggest to compensate  the lack of the Poincaré inequality by adding in the forced equation  a supplementary \emph{damping term} of the form $-\beta \vu$, where $\beta >0$ is a damping parameter. It is worth mention that another (merely technical) damping terms can be considered to study the long-time behavior of the Navier-Stokes equations on the whole space \cite{Chamorro,Liua}. However, we will consider here the damping term $-\beta \vu$, since this term  has a physical meaning as a drag-friction term in some oceanic models \cite{Pedlosky}, and then, it is also interesting  from the physical point of view. 

\medskip 

Thus, we shall consider here the following Cauchy problem of the  \emph{damped  Navier-Stokes-Bardina's model} for  incompressible fluids in the whole space $\Rt$: 

\begin{equation}\label{Bardina}
\left\{ \begin{array}{ll}\vspace{2mm}
\partial_t \vu + \text{div}((\vu \otimes \vu)_\alpha) - \nu \Delta \vu + \vec{\nabla} p = \fe -\beta \vu, \quad \alpha>0, \,\, \beta>0,\,\, \nu>0, \\ \vspace{2mm}
\text{div}(\vu)=0,\\
\vu(0,\cdot)= \vu_0,  \quad \text{div}(\vu_0)=0.
\end{array} 
 \right.
\end{equation}

When $\beta=0$,  the system (\ref{Bardina}) writes down as the classical  Navier-Stokes-Bardina's model. However, all the results we obtained in this article deeply depend on the the parameter $\beta>0$ and, from now on, we will focus on this model, where, our main objective is to \emph{describe the long-time behavior of its solutions}.  

%%%%%%%%%%%%%%%%%%%%%%%%%%%%%%%%%%%%%%%%%%%%%%

\section{Statement of the results}\label{Sec:Results}
Our first main result is devoted to the well-posedness of equations (\ref{Bardina}) in the energy space. 
\begin{TheoremeP}\label{Th-WP} Let $\vu_0 \in H^1(\Rt)$ be the initial datum such that $div(\vu_0)=0$. Moreover, let $\fe \in H^1(\Rt)$ be a stationary external force such that $div(\fe)=0$. Then, for all $\alpha >0$ and $\beta>0$ there exists a couple of functions $\vu=\vu_{\alpha,\beta} \in L^{\infty}([0,+\infty[, H^1(\Rt)) \cap L^{2}_{loc}([0, +\infty[, \dot{H}^2(\Rt))$, and $p=p_{\alpha,\beta} \in L^{2}_{loc}([0,+\infty[, H^3(\Rt))$, such that $(\vu, p)$ is the \emph{unique} weak solution of   (\ref{Bardina}). Moreover, this solution verifies the following energy \emph{equality} for all $t \geq 0$:
\begin{equation}\label{energ}
\begin{split}
\Vert & \vu(t,\cdot)\Vert^{2}_{L^2}+\alpha^2 \Vert \vu(t,\cdot)\Vert^{2}_{\dot{H}^1}+ 2\nu \int_{0}^{t}\Vert \vu(s,\cdot)\Vert^{2}_{\dot{H}^{1}} ds + 2 \alpha^2 \int_{0}^{t}\Vert \vu(s,\cdot)\Vert^{2}_{\dot{H}^{2}} ds = \Vert \vu_0 \Vert^{2}_{L^2}+\alpha^2\Vert \vu_0 \Vert^{2}_{\dot{H}^1} \\
&  + 2 \int_{0}^{t} \left( \langle \fe, \vu(s,\cdot) \rangle_{L^2\times L^2} + \alpha^2 \langle \vec{\nabla} \otimes \fe, \vec{\nabla} \otimes \vu(s,\cdot)\rangle_{L^2\times L^2} \right) ds - 2\beta \int_{0}^{t}  \left(\Vert \vu(s,) \Vert^{2}_{L^2}+\alpha^2 \Vert \vu(s,\cdot) \Vert^{2}_{\dot{H}^{1}}\right) ds. 
\end{split}
\end{equation}	
\end{TheoremeP}	

Comparing this result with the classical result on the existence of Leray's weak solutions for the Navier-Stokes equations (\ref{NS}), we may observe two main differences. First, we have here the \emph{uniqueness} of the weak solution in the energy space, and moreover, this solution verifies an energy \emph{equality}.  These  facts are due to  filtering operator $\ds{(\cdot)_\alpha}$   (defined in expression (\ref{filtering})) in the nonlinear transport term. As mentioned,  the filtering operator regularizes the classical nonlinear term in the  framework of the non homogeneous Sobolev spaces, providing the weak solutions of the equations (\ref{Bardina}) these good properties. In particular, the uniqueness of weak solutions is the key idea to study their long-time behavior and we will get back to this point later. 

\medskip 

Another important contribution of the filtering operator $(\cdot)_\alpha$ is the fact that, following some of the ideas in \cite{Berselli}, for $t>0$ the regularity (in the spatial variable) of the  weak solution $\vu(t,\cdot)$ constructed in this theorem could be improved to Sobolev spaces of higher order (provided that the external force is regular enough). However, the natural regularity given by the energy space is enough to study the asymptotic time behavior of the solution $\vu(t,\cdot)$. 

\medskip 

Concerning  the damping term $-\beta \vu$,  is it  worth mention this term   allows us to derive the following controls in time. These estimates will be very useful when we studying the large-time behavior of weak solutions.

\begin{PropositionP}\label{Prop:time-control} Within the framework of Theorem \ref{Th-WP}, the solution $\vu(t,\cdot)$ verifies the following estimates. 
\begin{enumerate}
\item[1)] For  all $t\geq 0$:  
\begin{equation}\label{control-tiempo}
\Vert \vu(t,\cdot)\Vert^{2}_{L^2}+\alpha^2\Vert \vu(t,\cdot)\Vert^{2}_{\dot{H}^1} \leq \left( \Vert \vu_0 \Vert^{2}_{L^2}+\alpha^2\Vert \vu_0 \Vert^{2}_{\dot{H}^1}\right)e^{-\beta\, t} + \frac{4}{\beta^2} \left( \Vert \fe \Vert^{2}_{L^2}+\alpha^2\Vert \fe \Vert^{2}_{\dot{H}^1}\right).
\end{equation}
\item[2)] For all $t \geq 0$ and $T\geq 0$:
\begin{equation}\label{control-tiempo2}
\nu \int_{t}^{t+T} \Vert \vu(s,\cdot)\Vert^{2}_{\dot{H}^1} ds + \alpha^2 \int_{t}^{t+T} \Vert \vu(s,\cdot)\Vert^{2}_{\dot{H}^2}ds \leq \frac{2 T}{\beta}\left( \Vert \fe \Vert^{2}_{L^2}+\alpha^2 \Vert \fe \Vert^{2}_{\dot{H}^1}\right)+ \Vert \vu(t,\cdot)\Vert^{2}_{L^2}+\alpha^2\Vert \vu(t,\cdot)\Vert^{2}_{\dot{H}^1}.    
\end{equation}
\end{enumerate}
\end{PropositionP}

The uniqueness of  global weak solutions constructed in Theorem \ref{Th-WP}  will be one of the key tools to study their long-time asymptotic behavior. Inspired by \cite{Ilyn}, our approach will be given  through the language  of the dynamical systems.  Before to state our next result,  for the sake of a complete exposition, we will set first some notations and we will recall some classical definitions in the theory of dynamical systems.  We refer to \cite{Raugel,Temam} and the references therein for more details. 

\medskip 

From now on, we fix the filtering parameter $\alpha >0$. Then, we define  the space $\mathcal{H}^{1}_{\alpha}(\Rt)$ as the Banach space of \emph{divergence free} vector fields $\vg \in H^1(\Rt)$   with the norm  $\ds{\Vert \vg \Vert^{2}_{\mathcal{H}^{1}_{\alpha}}=\Vert \vg \Vert^{2}_{L^2}+\alpha^{2} \Vert \vg \Vert^{2}_{\dot{H}^1}}$.   This \emph{equivalent} norm on the space $H^1(\Rt)$ naturally appears in the estimate (\ref{control-tiempo}), which will be  useful in the proof of our next result below. 

\medskip 

The key link between the equation (\ref{Bardina}) and the framework of dynamical systems is the fact that  for a given (stationary and divergence free) external force $\fe \in  \mathcal{H}^{1}_{\alpha}(\Rt)$,  the equation   (\ref{Bardina}) defines  a \emph{semigroup} acting on the Banach space $\mathcal{H}^{1}_{\alpha}(\Rt)$. More precisely,  for a time $t\geq 0$ we define the semigroup  $\ds{S(t)}: \mathcal{H}^{1}_{\alpha}(\Rt) \to \mathcal{H}^{1}_{\alpha}(\Rt)$ as 
\begin{equation}\label{def-semigrupo}
S(t)\vu_0 = \vu(t,\cdot), \quad \text{for all} \quad \vu_0 \in \mathcal{H}^{1}_{\alpha}(\Rt),
\end{equation} where $\vu(t,\cdot)$ is the \emph{unique}  global weak solution of equation (\ref{Bardina}) arising from the initial datum $\vu_0$ and constructed in Theorem \ref{Th-WP}. Due to the uniqueness of weak solutions, it is easy to verify that the family $\ds{(S(t))_{t\geq 0}}$ given  in (\ref{def-semigrupo}) defines  a (strongly continuous) semigroup on the  Banach space $\mathcal{H}^{1}_{\alpha}(\Rt)$. 

\medskip 

The study of the long-time asymptotics of weak solutions for equation (\ref{Bardina}) can be treated through the study of the semigroup $S(t)$ when $t \to +\infty$. More precisely, our aim is to prove that this semigroup has a \emph{global attractor} whose definition we recall as follows. 

\begin{DefinitionP}\label{Def:atractor} Let $\ds{(S(t))_{t\geq 0}}$ be the semigroup given in (\ref{def-semigrupo}) acting on the Banach space $\mathcal{H}^{1}_{\alpha}(\Rt)$. A global attractor for the semigroup $\ds{(S(t))_{t\geq 0}}$ is a set $\mathcal{A} \subset \mathcal{H}^{1}_{\alpha}(\Rt)$ which verifies: 
\begin{enumerate}
\item[$1)$] The set $\mathcal{A}$ is compact in $\mathcal{H}^{1}_{\alpha}(\Rt)$. 
\item[$2)$] The set 	$\mathcal{A}$ is strictly invariant: for all time $t\geq 0$ we have $S(t)\mathcal{A}=\mathcal{A}$. 
\item[$3)$] For every bound set $B \subset \mathcal{H}^{1}_{\alpha}(\Rt)$ and for every neighborhood  $\mathcal{V} \subset \mathcal{H}^{1}_{\alpha}(\Rt)$ of the set  $\mathcal{A}$, there exists a time $T=T(B, \mathcal{V})>0$, depending on the set $B$ and the neighborhood $\mathcal{V}$, such that for all time $t >T$ we have $S(t)B \subset \mathcal{V}$.  \end{enumerate}		
\end{DefinitionP}	

%%%%%%%%%%%%

In this definition  we focus on point $3)$ to remark that, roughly speaking, a global attractor \emph{attires}  the image through the semigroup $S(t)$  of all bounded set $B \subset \mathcal{H}^{1}_{\alpha}(\Rt)$ when $t\to +\infty$. This property allows us to have a better comprehension of the long-time behavior of weak solutions for the equation (\ref{Bardina}). Indeed, we observe first  that for all initial data $\vu_0 \in \mathcal{H}^{1}_{\alpha}(\Rt)$, setting the bounded set  $B = \{ \vu_0\}$, and moreover, setting  $\mathcal{V}$ any neighborhood of the attractor $\mathcal{A}$, then we find  that the solution $\vu(t,\cdot)$ of the equation (\ref{Bardina}) (arising from the initial datum $\vu_0$) lies in the neighborhood $\mathcal{V}$ from a time $T=T(\vu_0, \mathcal{V})$. Consequently, from \emph{any initial datum} $\vu_0$  the solution $\vu(t,\cdot)$ is as close to the attractor $\mathcal{A}$ as we  want when $t\to +\infty$. 

\medskip 

 In our second main result, we prove the existence of a global attractor for the semigroup $(S(t))_{t\geq 0}$ associated to equation (\ref{Bardina}). 
 
 \begin{TheoremeP}\label{Th-Atractor}  Let $\fe \in \mathcal{H}^{1}_{\alpha}(\Rt)$ be a stationary and divergence free external force. Let $(S(t))_{t\geq 0}$ be the semigroup associated to (\ref{Bardina}) defined in (\ref{def-semigrupo}). Then, this semigroup has a global attractor $\mathcal{A}_{\fe}$, depending on the external force $\fe$,  given by Definition \ref{Def:atractor}. 
 \end{TheoremeP}
 
Once we have proven the existence of a global attractor for the equation (\ref{Bardina}), our general aim is to give more precisely descriptions of the set $\mathcal{A}_{\fe}$. Thus,  to present all our results in an orderly fashion, we divided them in four sections below.

\subsection*{Fractal dimension of the global attractor}

As the global attractor is a compact set of $\mathcal{H}^{1}_{\alpha}(\Rt)$, it is natural to measure its \emph{size} in some sense; and this is the aim of our next result.  Specifically,  we prove that the global attractor   $\mathcal{A}_{\fe}$ has a finite fractal dimension, and moreover, we derive an explicit upper bound for this dimension in terms of the parameter $\alpha, \beta, \nu$ and the norm of the external force $\ds{\Vert \fe \Vert_{\mathcal{H}^{1}_{\alpha}}}$. For this, we will quickly recall the definition of the fractal dimension through the so-called \emph{box-counting} method. For more details we refer the reader to \cite{Babin,Temam} and the references therein. 

\medskip  
 
Let $ \mathcal{A}_{\fe}  \subset \mathcal{H}^{1}_{\alpha}(\Rt)$ be the global attractor associated to equation (\ref{Bardina}). Then, by the Hausdorff criterion, for every $\epsilon > 0$
it can be covered by the finite number of $\epsilon -$balls in $\mathcal{H}^{1}_{\alpha}(\Rt)$. Let $N_{\epsilon}(\mathcal{A}_{\fe})$ be the minimal number of such balls. We thus have the following definition. 
\begin{DefinitionP}\label{def-fractal-dim-attractor}   
The fractal (box-counting) dimension of the attarctor $\mathcal{A}_{\fe}$ in $\mathcal{H}^{1}_{\alpha}(\Rt)$ is defined via the following expression:
\begin{equation*}
\mbox{dim} \left(\mathcal{A}_{\fe}\right)= \limsup_{\epsilon \to 0} \frac{\ln \ N_{\epsilon}(\mathcal{A}_{\fe})}{\ln\left( \frac{1}{\epsilon}\right)}.
\end{equation*}
\end{DefinitionP}    

In our third result we derive an explicit upper bound of  $\ds{\mbox{dim} \left(\mathcal{A}_{\fe}\right)}$.

\begin{TheoremeP}\label{Th-dim-attractor}
Let the assumptions of Theorem \eqref{Th-Atractor} hold. Then the fractal dimension of the global attractor $\mathcal{A}$ associates to equation \eqref{Bardina}  satisfies the following estimate:
\begin{equation}\label{dimfrac}
\mbox{dim} \left(\mathcal{A}_{\fe}\right) \leq  c(\alpha,\beta,\nu)\max\left(\Vert \fe \Vert^{14/5}_{\mathcal{H}^{1}_{\alpha}}, \Vert \fe \Vert^{2}_{\mathcal{H}^{1}_{\alpha}}  \right),  
\end{equation}
where the constant $c(\alpha,\beta,\nu) > 0$, depending only of the parameters $\alpha,\beta$ and $\nu$, is explicitly  given in the formula (\ref{defi-constant-Th}).
\end{TheoremeP} 

We observe in this estimate that the fractal dimension of the global attractor $\mathcal{A}_{\fe}$ is essentially  controlled for above by the size of the external force $\fe$ in the space $\mathcal{H}^{1}_{\alpha}(\Rt)$. This type of control was also pointed out in \cite{Ilyn}, for the case of the two-dimensional and damped Navier-Stokes equations, while,  for the Navier-Stokes-Bardina's model in the space-periodic case, similar upper bounds on $\mbox{dim} \left(\mathcal{A}_{\fe}\right) $ were established in \cite{Cao}. Finally, let us  mention that this is a first estimation for an upper bound of  $\ds{\mbox{dim} \left(\mathcal{A}_{\fe}\right)}$ and the optimally of this upper bound, or moreover the derivation of some lower bounds, are matter of further investigations. 

\subsection*{Internal structure of the global attractor} 

We are also interested in characterizing the global attractor $\mathcal{A}_{\fe}$. By \cite{Babin},\cite{Babin2} and \cite{Raugel}  it is well-known that   the global attractor  can be described  through a particular kind of solution for the equations  (\ref{Bardina}). Such solutions are called the  \emph{eternal solutions} which, as we will observe in the following definition, they do not arise from any initial data and they are actually defined for all time $t \in \mathbb{R}$.   

\begin{DefinitionP}\label{defi-eternal-sol} Let $\fe \in \mathcal{H}^{1}_{\alpha}(\Rt)$ be a stationary and divergence-free external force.  We say a couple $(\vv,q)$ is an eternal solution for the damped Navier-Stokes-Bardina  equations  with force  $\fe$, if 
\begin{equation*}
\vv \in L^{\infty}_{loc}\left(\mathbb{R},\mathcal{H}^{1}_{\alpha}(\Rt)\right)\cap L^{2}_{loc}(\mathbb{R},\dot{H}^2(\Rt)), \quad q \in  L^{2}_{loc}(\mathbb{R}, \dot{H}^{3}(\Rt)),
\end{equation*}  and if $(\vv,q)$ is a weak solution of the equations 
\begin{equation}\label{Bardina-enternal}
\partial_t \vv + \text{div}((\vv \otimes \vv)_\alpha) - \nu \Delta \vv + \vec{\nabla} q = \fe-\beta \vv, \quad  \text{div}(\vv)=0. 
\end{equation}
\end{DefinitionP}

In Proposition \ref{Prop:eternal-sol} below, we prove the existence of eternal solutions for the damped Navier-Stokes-Bardina's model. 

\medskip

Denoting by $\mathcal{E}_{\fe}$ the set of all the \emph{bounded} eternal solutions $(\vv,q)$ associated to the force $\fe$, \emph{i.e.}, the eternal solutions verifying $\vv \in L^{\infty}(\mathbb{R}, \mathcal{H}^{1}_{\alpha}(\Rt))$, by Lemma $2.18$ (page 16) in \cite{Raugel},  we have that   the global attractor $\mathcal{A}_{\fe}$   given by the Theorem \ref{Th-Atractor}  has the following structure: 

\begin{equation}\label{Internal-estructure}
 \mathcal{A}_{\fe} =   \left. \mathcal{E}_{\fe}\, \right\vert_{t=0}.  
\end{equation}
In other words,  the global attractor $\mathcal{A}_{\fe}$ of the equations (\ref{Bardina}) consists of the set of functions $\vv(0, \cdot)$, where $\vv(t,\cdot)$ is an \emph{bounded} eternal solution of the damped Navier-Stokes-Bardina  equations given in Definition  \ref{defi-eternal-sol}, and thus, its internal structure is explicitly described by the identity above. 

\medskip 

On the other hand, a particular  case of  eternal solutions for the damped Navier-Stokes-Bardina equations are the  \emph{stationary  solutions}. These solutions,  which will be  denoted by  $(\U,P)$, only depend on the spatial variable and solve the following  elliptic equation:

\begin{equation}\label{Bardina-Stat}
\left\{ \begin{array}{ll}\vspace{2mm}
-\nu \Delta \U + \text{div}((\U \otimes \U)_\alpha) + \vec{\nabla} P = \fe -\beta \U, \quad \alpha>0, \,\, \beta>0,\,\, \nu>0, \\ 
\text{div}(\U)=0.
\end{array} 
 \right.
\end{equation}

In our next result, we prove first the existence of stationary solutions, and moreover, we investigate their relation with the global attractor. 

%We are thus interested in investigating when the stationary solutions $\U$ belong to the global attractor $\mathcal{A}_{\fe}$. Before to give an answer to this question, and for the sake of completeness of this paper, in our fourth result we prove first a general result on the existence of stationary solutions which is also of independent interest.  

\begin{TheoremeP}\label{Th-Stationary-Solutions} Let $\fe \in \mathcal{H}^{1}_{\alpha}(\Rt)$ be the external force such that $div(\fe)=0$. Then, the following statements hold: 
\begin{enumerate} 
\item[$1)$] There exist $\U \in H^2(\Rt)$ and $P \in H^{1}(\Rt)$ such that the couple $(\U,P)$ is a solution of the equation (\ref{Bardina-Stat}).
\item[$2)$] All the stationary solutions verify the following energy estimate: $\ds{\Vert \U  \Vert^{2}_{\mathcal{H}^{1}_{\alpha}} + \nu \alpha^2 \Vert \U \Vert^{2}_{\dot{H}^2}\leq \frac{2}{\beta^2} \Vert \fe \Vert^{2}_{\mathcal{H}^{1}_{\alpha}}}$. 
\item[$3)$] All the stationary solutions belong to the global attractor $\mathcal{A}_{\fe}$. 
\end{enumerate}
\end{TheoremeP} 

We have the following comments. The result given in point $1)$ establishes the existence of stationary solutions for the damped Navier-Stokes-Bardina equation with \emph{any} external force, and thus, this is a general result for the elliptic equation  (\ref{Bardina-Stat})  which is also of independent interest. In Section \ref{Sec:Stationary-Sol} we comment more in details our strategy to prove this point which is based in the Scheafer's fixed point argument. Moreover, it is worth mention the uniqueness issue for the stationary solutions, in the general case of any external force, seems to be more delicate and it is matter of further investigations. 

\medskip 

On the other hand, in point $2)$ we show that \emph{all} the stationary solutions belong  to the space $H^{1}_{\alpha}(\Rt)$ and their norms are always controlled by the norm of the external force.  Finally, maybe  the most interesting feature on the stationary solutions is given in point $3)$, where we ensure that \emph{all} the stationary solutions  fallen inside the global attractor $\mathcal{A}_{\fe}$. 

\medskip 

\subsection*{Additional properties of the global attractor driven by the damping parameter} 

In this section, we study the role of the damping parameter $\beta>0$ in the description of the global attractor for the equation  (\ref{Bardina}). We start by setting some notation. For the external force $\fe \in \mathcal{H}^{1}_{\alpha}(\Rt)$, and for a numerical constant $c>0$, we introduce the following quantity depending on the damping parameter $\beta>0$: 

\begin{equation}\label{eta}
\eta(\beta)=  \frac{c}{\alpha^{5/2}\beta}\Vert \fe \Vert_{\mathcal{H}^{1}_{\alpha}} -\beta.
\end{equation}

In our next result we prove that in the case when the parameter $\beta>0$ is big enough, in the sense that this  quantity  verifies $\eta(\beta)\leq 0$, we are able to give sharp properties of the global attractor. More precisely, we will consider first the case when $\beta>0$ is such that  $\eta(\beta)=0$.  In this case, we study some kind of stability of the elements of the global attractor, also called the \emph{orbital stability}, which its definition we recall below.  For more references see the Section $1.1$ of \cite{Brieve}.  

\medskip 

We say that $\vu_0 \in \mathcal{H}^{1}_{\alpha}(\Rt)$ is orbitally stable if for all $\varepsilon>0$, there exists $\delta=\delta(\varepsilon)>0$, such that for all $\vv_0 \in \mathcal{H}^{1}_{\alpha}(\Rt)$  verifying 
\begin{equation*}
 \Vert \vu_0 - \vv_0 \Vert_{\mathcal{H}^{1}_{\alpha}} \leq \delta,    
\end{equation*}
then the solutions  $\vu(t,\cdot)$ and $\vv(t,\cdot)$ to the equation (\ref{Bardina}) arising from $\vu_0$ and $\vv_0$ respectively, satisfy 
\begin{equation*}
 \sup_{t \geq 0} \Vert \vu(t,\cdot)-\vv(t,\cdot) \Vert_{\mathcal{H}^{1}_{\alpha}} \leq \varepsilon.  
\end{equation*}
%This notion essentially asserts that if  $\vu_0$ is orbitally stable then the solution $\vv(t,\cdot)$  can stay so close to the solution $\vu(t,\cdot)$  as we want, provided that its initial datum $\vv_0$ is close enough to $\vu_0$.\\

On the other hand, when $\beta>0$ is such that $\eta(\beta)<0$, we  go further in the description of the global attractor $\mathcal{A}_{\fe}$. In this case, surprisingly, the global attractor contains a single element given by the  \emph{unique} solution of the stationary equation (\ref{Bardina-Stat}).

\medskip 

Summarizing, our result on the  role of the damping parameter $\beta>0$ in the description of the global attractor reads as follows. 

\begin{TheoremeP}\label{Th:Stability}
 Let $\mathcal{A}_{\fe}\subset \mathcal{H}^{1}_{\alpha}(\Rt)$ be global attractor of the equations (\ref{Bardina}) given by Theorem \ref{Th-Atractor}. Then, the following statements hold:
\begin{enumerate}
    \item[$1)$]  If $\beta>0$ is such that $\eta(\beta)=0$, then all the elements in the attractor $\mathcal{A}_{\fe}$ are \emph{orbitally stable}.  
    \item[$2)$]  If $\beta>0$ is such that $\eta(\beta)<0$, then the stationary solution $(\U,P)$ of equation (\ref{Bardina}) given by Theorem \ref{Th-Stationary-Solutions} is unique. Moreover,  the global attractor $\mathcal{A}_{\fe}$  only contains the unique stationary solution of equation (\ref{Bardina}).%, and we have:  
%\[ \ds{\mathcal{A}_{\fe}=\{\U\}}. \] 
\end{enumerate}
\end{TheoremeP}

The result given in point $2)$ ensures that as long as the damping parameter $\beta$ is such that  $\eta(\beta)<0$, all the weak solutions of the damped Navier-Stokes-Bardina equation (\ref{Bardina}) are attracted by the unique stationary solution of this equation when  the time is large enough. In this case, we are also able to give a sharp \emph{asymptotic profile} in time of the  solutions to the equation (\ref{Bardina}).  

\begin{PropositionP}\label{Prop:profile}
 Within the framework of point $2)$ in Theorem \ref{Th:Stability}, let $\U \in \mathcal{A}_{\fe}$   be the unique solution of the stationary problem  (\ref{Bardina-Stat}). Then, for all $x\in \Rt$ fixed, all the solutions of the equation (\ref{Bardina}) have the following \emph{asymptotic profile} in time: 
\begin{equation} 
\vu(t,x)= \U(x)+  \mathcal{R}_{\vu}(t,x) ,  \quad t>0,       
\end{equation}    
where the remainder term  $\mathcal{R}_{\vu}(t,x)$ is a vector field depending on $\vu$, which  verify the following time decaying  
\begin{equation}\label{time-decay1}
%\lim_{t\to +\infty} t^{3/4-\varepsilon} \vert  \mathcal{R}_{\vu}(t,x) \vert =0,  \quad  0<\varepsilon \ll 1.
\Vert \mathcal{R}_{\vu}(t,\cdot)\Vert_{L^{\infty}}\leq C\, t^{-3/4},  \quad t\gg 1, 
\end{equation}
with a constant $C>0$ depending on the the parameters $\alpha, \beta, \nu$, the initial value  $\vu(0,\cdot)$, and   $\fe$, $\U$. 
\end{PropositionP} 

\subsection*{The damped Navier-Stokes-Bardina's model without external force}

Finally, in this last part we consider the damped Navier-Stokes-Bardina equations in the particular case of a zero external force.  

\begin{equation}\label{Bardina-f-zero}
\left\{ \begin{array}{ll}\vspace{2mm}
\partial_t \vu + \text{div}((\vu \otimes \vu)_\alpha) - \nu \Delta \vu + \vec{\nabla} p =  -\beta \vu, \quad \alpha>0, \,\, \beta>0,\,\, \nu>0, \\ \vspace{2mm}
\text{div}(\vu)=0,\\
\vu(0,\cdot)= \vu_0,  \quad \text{div}(\vu_0)=0.
\end{array} 
 \right.
\end{equation} 

In this case, we give a sharp description of the global attractor associated with the zero force, more precisely, we show   that the global attractor of this equation only contains the zero function.  Moreover, we prove that  all the weak solutions to the equation (\ref{Bardina-f-zero})  have a fast (exponential) convergence rate to zero when the time goes to infinity.   

\begin{PropositionP}\label{Prop:force-zero} For the equation (\ref{Bardina-f-zero}), the unique global attractor verifies $\mathcal{A}_{0}=\{ 0 \}$. Moreover, all the weak  solutions $u(t,x)$  given by Theorem \ref{Th-WP} verify:
\[ \Vert \vu(t,\cdot)\Vert_{L^p} \leq C\, e^{- \frac{2\beta}{p}\, t}, \quad 2\leq p < +\infty, \quad t\gg 1,\]
where  the constant $C>0$ depends on the initial datum $\vu(0,\cdot)$ and the parameter $p$. 
\end{PropositionP} 

{\bf Organization of the paper.} In Section \ref{Sec:GWP} we prove the global well-posedness of the equation (\ref{Bardina}). Then, Section \ref{Sec:GlobalAttractor} is devoted to the existence of the global attractor associated to the equation (\ref{Bardina}), while in Section \ref{Sec:FractalDim} we derive an upper bound for its fractal dimension. Thereafter, Section \ref{Sec:Stationary-Sol} is devoted to the study of the internal structure of the global attractor, and moreover, we also construct weak solutions for the stationary problem (\ref{Bardina-Stat}). Then, in Section \ref{Sec:AdditionalProp} we study the sharp properties of the global attractor driven by the damping parameter and, finally, in Section \ref{Sec:ForceZero} we consider the case of the equations (\ref{Bardina}) without external force.

\section{Global well-posedness in the energy space}\label{Sec:GWP}
\subsection*{Proof of Theorem \ref{Th-WP}}
The proof of this theorem is rather straightforward and it follows essentially the same lines of the classical framework. The first step is to solve the following integral problem: 
\begin{equation}\label{Integral}
\vu(t,\cdot)= e^{\nu t \Delta}\vu_0 + \int_{0}^{t} e^{\nu(t-s)\Delta} \fe ds - \int_{0}^{t} e^{\nu(t-s)\Delta} \P div ((\vu \otimes \vu)_\alpha)(s,\cdot) ds - \beta \int_{0}^{t} e^{\nu (t-s) \Delta} \vu(s,\cdot) ds,
\end{equation} in the energy space  $\ds{E_T= L^{\infty}([0,T], H^1(\Rt)) \cap L^2([0,T], \dot{H}^{2}(\Rt))}$ (with $0<T<+\infty$)  with the norm  $\Vert \cdot \Vert_{T}= \Vert \cdot \Vert_{L^{\infty}_{t}H^{1}_{x}}+\Vert \cdot \Vert_{L^{2}_{t}\dot{H}^{2}_{x}}$.  We write 
\begin{equation*}
\Vert \vu(t,\cdot)\Vert_{T} \leq \underbrace{\left\Vert e^{\nu t \Delta}\vu_0 + \int_{0}^{t} e^{\nu(t-s)\Delta} \fe ds \right\Vert_T}_{(a)} + \underbrace{\left\Vert \int_{0}^{t} e^{\nu(t-s)\Delta} \P div ((\vu \otimes \vu)_\alpha)(s,\cdot) ds \right\Vert_T}_{(b)}   + \beta \underbrace{\left\Vert  \int_{0}^{t} e^{\nu (t-s) \Delta} \vu(s,\cdot) ds \right\Vert_T}_{(c)}. 
\end{equation*}
The term $(a)$ is classical to treat. For the term $\ds{e^{\nu t \Delta} \vu_0 }$ we have the following estimate (see the proof of Theorem $7.1$, page $131$ of the book \cite{PGLR1}):  
\begin{equation}\label{estim1}
\Vert e^{\nu t \Delta} \vu_0 \Vert_{E_T}\leq c(1+ 1/ \sqrt{ 2\nu}) \Vert \vu_0 \Vert_{H^1}.
\end{equation}
Moreover, to study the term $\ds{  \int_{0}^{t} e^{\nu(t-s)\Delta} \fe ds }$, we shall use the  following well-known estimates. See the Lemma $7.2$, page $129$ of the book \cite{PGLR1}.
\begin{Lemme}\label{LemaTech1} Let $g \in L^{2}([0,T], L^2(\Rt))$ and let $\ds{G(t,x)= \int_{0}^{t} e^{\nu (t-s)\Delta} g(s,x)ds}$. Then, $G(t,x)$ belongs to the space $E_T$ and we have the following estimates: 
\begin{enumerate}
\item[$1)$] $\ds{\Vert G(t,\cdot) \Vert_{L^{\infty}_{t}L^{2}_{x}} \leq c\, \sqrt{T} \, \Vert g \Vert_{L^{2}_{t}L^{2}_{x}}}.$ 
\item[$2)$] $\ds{ \Vert G(t,\cdot) \Vert_{L^{\infty}_{t}\dot{H}^{1}_{x}} \leq \frac{c}{\sqrt{2\nu}}\,  \Vert g \Vert_{L^{2}_{t}L^{2}_{x}}}$. 
\item[$3)$] $\ds{ \Vert G(t,\cdot) \Vert_{L^{2}_{t}\dot{H}^{2}_{x}} \leq \frac{c}{\nu}\, \Vert g \Vert_{L^{2}_{t}L^{2}_{x}}}$. 
\end{enumerate}   
\end{Lemme}	 
Thus, in this lemma we set $g=\fe$, and moreover, since $\fe$ is a time-independing function then we have 
\begin{eqnarray}\label{estim3}\nonumber
\left\Vert \int_{0}^{t} e^{\nu(t-s)\Delta} \fe ds \right\Vert_{E_T} &\leq & c (\sqrt{T}+ 1/\sqrt{2\nu}+ 1/\nu) \Vert \fe \Vert_{L^{2}_{t}L^{2}_{x}} \leq  c (\sqrt{T}+ 1/\sqrt{2\nu}+ 1/\nu)\sqrt{T} \Vert \fe \Vert_{L^2} \\
&\leq & c  (\sqrt{T}+ 1/\sqrt{2\nu}+ 1/\nu)\sqrt{T} \Vert \fe \Vert_{H^1}. 
\end{eqnarray}
At this point, we set the time $T\leq 1$ and then by the estimates (\ref{estim1}) and (\ref{estim3}) we can write 
\begin{equation}\label{estim-lin}
\left\Vert e^{\nu t \Delta}\vu_0 + \int_{0}^{t} e^{\nu(t-s)\Delta} \fe ds \right\Vert_T \leq c (1+ 1/\sqrt{2\nu}+ 1/\nu) (\Vert \vu_0 \Vert_{H^1}+ \Vert \fe \Vert_{H^1}) \leq c_\nu (\Vert \vu_0 \Vert_{H^1}+ \Vert \fe \Vert_{H^1}).  
\end{equation}
We study now the term $(b)$. We recall first that by (\ref{filtering}) we have $\ds{(\vu \otimes \vu)_\alpha}=(-\alpha^2 \Delta + I_d)(\vu \otimes \vu)$ and then, by well-known properties of the Bessel potential $\ds{(-\alpha^2 \Delta + I_d)^{-1}}$, we can write $\ds{(-\alpha^2 \Delta + I_d)^{-1}(\vu \otimes \vu)=K_\alpha \ast (\vu \otimes \vu)}$, where the kernel $K_\alpha(x)$ has good decaying properties. In particular we have $\Vert K_\alpha \Vert_{L^1}\leq c_\alpha$, for a constant $0<c_\alpha<+\infty$ depending on $\alpha >0$ (see the Section $6.1.2$ of the book \cite{Grafakos} for all the details).  \\

With this remark in mind, we start by estimating the quantity $\ds{\left\Vert \int_{0}^{t} e^{\nu(t-s)\Delta} \P div ((\vu \otimes \vu)_\alpha)(s,\cdot) ds \right\Vert_{L^{\infty}_{t} L^{2}_{x}}}$.  For $0<t\leq T$ fix, by the Young inequalities,   and moreover, by the continuity of the Leray projector in the Lebesgue spaces, we write 
\begin{eqnarray*}
& & \left\Vert  \int_{0}^{t} e^{\nu(t-s)\Delta} \P div ((\vu \otimes \vu)_\alpha)(s,\cdot) ds \right\Vert_{L^2}  \leq \int_{0}^{t} \Vert e^{\nu(t-s)\Delta} \P div ((\vu \otimes \vu)_\alpha)(s,\cdot) \Vert_{L^2} ds \\
& \leq & \int_{0}^{t} \Vert e^{\nu(t-s)\Delta} \P div (K_\alpha\ast  (\vu \otimes \vu))(s,\cdot) \Vert_{L^2} ds \leq \int_{0}^{t} \left\Vert  K_\alpha \ast \left( e^{\nu(t-s)\Delta} \P div (\vu \otimes \vu)\right)(s,\cdot) \right\Vert_{L^2} ds \\
&\leq & c_\alpha  \, \int_{0}^{t} \left\Vert  e^{\nu(t-s)\Delta} \P div (\vu \otimes \vu)(s,\cdot) \right\Vert_{L^2} ds \leq c_{\nu,\alpha}\,  \int_{0}^{t} \left\Vert   div (\vu \otimes \vu)(s,\cdot) \right\Vert_{L^2} ds\\
&\leq &c_{\nu,\alpha}\, \int_{0}^{t} \Vert \vu \otimes \vu(s,\cdot)\Vert_{\dot{H}^1} ds \leq c_{\nu,\alpha}\, T^{1/2} \,\Vert \vu \otimes \vu \Vert_{L^{2}_{t}\dot{H}^{1}_{x}}. 
\end{eqnarray*}
Here we shall need the following technical lemma.
\begin{Lemme}\label{LemmaTech2} We have $\Vert \vu \otimes \vu \Vert_{L^{2}_{t}\dot{H}^{1}_{x}} \leq c\, T^{1/4}\, \Vert \vu \Vert^{2}_{E_T}$. 
\end{Lemme}	
\pv Using first the product laws in the homogeneous Sobolev spaces (see the Lemma $7.3$, page $130$ of the book \cite{PGLR1}) and then, using the H\"older inequalities in the time variable  we get $\ds{\Vert \vu \otimes \vu \Vert_{L^{2}_{t}\dot{H}^{1}_{x}} \leq c\, \Vert \vu \Vert_{L^{4}_{t}\dot{H}^{3/2}_{x}}\Vert \vu \Vert_{L^{4}_{t} \dot{H}^{1}_{x}}}$.  Thereafter, to estimate  the quantity $\ds{\Vert \vu \Vert_{L^{4}_{t}\dot{H}^{3/2}_{x}}}$, we use the interpolation inequalities (first in the spatial variable and then in the temporal variable) and we have 
$\ds{\Vert \vu \Vert_{L^{4}_{t}\dot{H}^{3/2}_{x}} \leq c \Vert \vu \Vert^{1/2}_{L^{\infty}_{t} \dot{H}^{1}_{x}} \, \Vert \vu \Vert^{1/2}_{L^{2}_{t}\dot{H}^{2}_{x}}\leq c \Vert \vu \Vert_{E_T}}$. Finally, the quantity $\ds{\Vert \vu \Vert_{L^{4}_{t} \dot{H}^{1}_{x}}}$ is directly estimated as follows: $\ds{\Vert \vu \Vert_{L^{4}_{t} \dot{H}^{1}_{x}} \leq c \, T^{1/4} \Vert \vu \Vert_{L^{\infty}_{t} \dot{H}^{1}_{x}} \leq c\, T^{1/4} \Vert \vu \Vert_{E_T}}$. \finpv

With this estimate at hand, and recalling that we have assumed $T\leq 1$, we can write 
\begin{equation}\label{estim4}
\left\Vert  \int_{0}^{t} e^{\nu(t-s)\Delta} \P div ((\vu \otimes \vu)_\alpha)(s,\cdot) ds \right\Vert_{L^{\infty}_{t}L^{2}_{x}}  \leq c_{\nu,\alpha}\, T^{1/4} \Vert \vu \Vert^{2}_{E_T}.
\end{equation} 

We study now the quantity $\ds{\left\Vert \int_{0}^{t} e^{\nu(t-s)\Delta} \P div ((\vu \otimes \vu)_\alpha)(s,\cdot) ds \right\Vert_{L^{\infty}_{t} \dot{H}^{1}_{x}}}$. By point $2)$ in Lemma \ref{LemaTech1}, where we set now  $g= \P div ((\vu \otimes \vu)_\alpha) $, and moreover, by Lemma \ref{LemmaTech2} we can write
\begin{eqnarray}\label{estim5}\nonumber
\left\Vert \int_{0}^{t} e^{\nu(t-s)\Delta} \P div ((\vu \otimes \vu)_\alpha)(s,\cdot) ds \right\Vert_{L^{\infty}_{t}\dot{H}^1} &\leq & c_\nu\, \Vert \P div ((\vu \otimes \vu)_\alpha \Vert_{L^{2}_{t}L^{2}_{x}} 	\leq c_\nu \, \Vert K_\alpha \ast \left( \P div (\vu\otimes \vu)\right)\Vert_{L^{2}_{t}L^{2}_{x}} \\
&\leq  & c_{\nu,\alpha} \, \Vert \vu\otimes \vu \Vert_{L^{2}_{t}\dot{H}^{1}_{x}} \leq c_{\nu,\alpha}\, T^{1/4}\, \Vert \vu \Vert^{2}_{E_T}.
\end{eqnarray}

Finally, we study the quantity $\ds{\left\Vert \int_{0}^{t} e^{\nu(t-s)\Delta} \P div ((\vu \otimes \vu)_\alpha)(s,\cdot) ds \right\Vert_{L^{2}_{t}\dot{H}^{2}_{x}}}$. As the previous quantity, by point $3)$ in Lemma \ref{LemaTech1} and moreover by Lemma \ref{LemmaTech2} we have 
\begin{equation}\label{estim6}
\left\Vert \int_{0}^{t} e^{\nu(t-s)\Delta} \P div ((\vu \otimes \vu)_\alpha)(s,\cdot) ds \right\Vert_{L^{2}_{t}\dot{H}^2} \leq  c_\nu\, \Vert \P div ((\vu \otimes \vu)_\alpha \Vert_{L^{2}_{t}L^{2}_{x}} \leq  c_{\nu,\alpha}\, T^{1/4}\, \Vert \vu \Vert^{2}_{E_T}.
\end{equation}
Thus, gathering the estimates (\ref{estim4}), (\ref{estim5}) and (\ref{estim6}) we obtain 
\begin{equation}\label{estin-bilin}
\left\Vert \int_{0}^{t} e^{\nu(t-s)\Delta} \P div ((\vu \otimes \vu)_\alpha)(s,\cdot) ds \right\Vert_{E_T}\leq c_{\nu,\alpha}\, T^{1/4}\, \Vert \vu \Vert^{2}_{E_T}. 
\end{equation}
It remains to estimate the term $(c)$. By Lemma \ref{LemaTech1} (where we set now $g=\vu $, and in point $1)$  we recall that   $T\leq 1$)  we write 
\begin{eqnarray}\label{estil-lin2}
\left\Vert \int_{0}^{t}e^{\nu(t-s)\Delta} \vu(s,\cdot)ds \right\Vert_{E_T} \leq c_\nu \, \Vert \vu \Vert_{L^{2}_{t}L^{2}_{t}} \leq c_\nu\, T^{1/2} \, \Vert \vu \Vert_{E_T}.
\end{eqnarray} 

Once we have inequalities (\ref{estim-lin}), (\ref{estin-bilin}) and (\ref{estil-lin2}), for a time $T>0$ small enough by the Banach contraction principle we obtain a local  solution $\vu \in E_T$ of equations (\ref{Integral}).\\

The second step is to prove that this solution is global in time. Remark first that the solution $\vu$ obtained above   also solves the problem  
\begin{equation*}
\partial_t \vu + \P(\text{div}((\vu \otimes \vu)_\alpha)) - \nu \Delta \vu  = \fe -\beta \vu, 
\end{equation*} in the distributional sense, where,  as $\vu \in  E_T$ then  each term in this equation  belong to the space $L^{2}([0,T], L^2(\Rt))$.  By the identity  (\ref{filtering}) we can write 
\begin{equation}\label{Eq-Aux}
\partial_t \vu + \P(\text{div}((-\alpha^2 \Delta +I_d)^{-1} (\vu \otimes \vu)) - \nu \Delta \vu  = \fe -\beta \vu,
\end{equation} and applying the operator $(-\alpha \Delta + I_d)$ in each term we get that the solution $\vu $ also verifies  the following equation
\begin{equation*}
(-\alpha^2 \Delta + I_d) \partial_t \vu=- \P(\text{div}(\vu \otimes \vu)) + \nu (-\alpha^2\Delta +I_d) \Delta \vu+ (-\alpha^2\Delta +I_d) \fe - \beta (-\alpha^2\Delta +I_d) \vu.
\end{equation*} Here  each term belong to the space $L^{2}([0,T], H^{-2}(\Rt))$.  Now, always by the fact  $\vu \in E_T$ we get $\vu \in L^{2}([0,T], H^2(\Rt))$ and then  we can write 
\begin{equation}\label{Eq-dif-energ}
\begin{split}
\frac{1}{2} \frac{d}{dt}\left(\Vert \vu(t,\cdot) \Vert^{2}_{L^2}+ \alpha^2 \Vert \vu(t,\cdot)\Vert^{2}_{\dot{H}^1}\right)=&  \left\langle (-\alpha^2 \Delta + I_d) \partial_t \vu(t,\cdot), \vu(t,\cdot) \right\rangle_{H^{-2}\times H^2}\\ 
=& -\nu \Vert \vu(t,\cdot) \Vert^{2}_{\dot{H}^1}-\alpha^2 \Vert \vu(t,\cdot)\Vert^{2}_{\dot{H}^2} + \langle \fe, \vu(t,\cdot)\rangle_{L^2\times L^2} \\
 & +\alpha^2 \langle \vec{\nabla} \otimes \fe , \vec{\nabla} \otimes \vu(t,\cdot)\rangle_{L^2\times L^2} -\beta \left( \Vert \vu(t,\cdot)\Vert^{2}_{L^2}+\alpha^2 \Vert\vu(t,\cdot)\Vert^{2}_{\dot{H}^1}\right). 
\end{split}
\end{equation} 
As the quantity $\ds{-\beta \left( \Vert \vu(t,\cdot)\Vert^{2}_{L^2}+\alpha^2 \Vert\vu(t,\cdot)\Vert^{2}_{\dot{H}^1}\right)}$ is negative, and moreover, applying the Cauchy-Schwarz inequalities, we obtain 
\begin{equation*}
\begin{split}
\frac{1}{2} \frac{d}{dt} \Big(\Vert \vu(t,\cdot) \Vert^{2}_{L^2}+ & \alpha^2 \Vert \vu(t,\cdot)\Vert^{2}_{\dot{H}^1} \Big)  \\
\leq &-\nu \Vert \vu(t,\cdot) \Vert^{2}_{\dot{H}^1}-\alpha^2 \Vert \vu(t,\cdot)\Vert^{2}_{\dot{H}^2} + \langle \fe, \vu(t,\cdot)\rangle_{L^2\times L^2}  +\alpha^2 \langle \vec{\nabla} \otimes \fe , \vec{\nabla} \otimes \vu(t,\cdot)\rangle_{L^2\times L^2}\\
\leq & -\nu \Vert \vu(t,\cdot) \Vert^{2}_{\dot{H}^1}-\alpha^2 \Vert \vu(t,\cdot)\Vert^{2}_{\dot{H}^2}  + 
\Vert \fe \Vert_{L^2}\Vert \vu(t,\cdot)\Vert_{L^2} + \alpha^2 \Vert \fe \Vert_{\dot{H}^1}\Vert \vu(t,\cdot)\Vert_{\dot{H}^{1}} \\
\leq & -\nu \Vert \vu(t,\cdot) \Vert^{2}_{\dot{H}^1}-\alpha^2 \Vert \vu(t,\cdot)\Vert^{2}_{\dot{H}^2}  + \Vert \fe \Vert^{2}_{L^2} + \Vert \vu(t,\cdot) \Vert^{2}_{L^2}+ \alpha^2\Vert \fe \Vert^{2}_{\dot{H}^{1}}+\alpha^2\Vert \vu(t,\cdot)\Vert^{2}_{\dot{H}^1}\\
\leq & \left( \Vert\vu(t,\cdot)\Vert_{L^2}+\alpha^2\Vert \vu(t,\cdot)\Vert^{2}_{\dot{H}^{1}}\right) - \nu \Vert \vu(t,\cdot)\Vert^{2}_{\dot{H}^{1}} - \alpha^2\Vert \vu(t,\cdot)\Vert^{2}_{\dot{H}^{2}} + \Vert \fe \Vert^{2}_{L^2}+\alpha^2\Vert \fe \Vert^{2}_{\dot{H}^{1}}.
\end{split}
\end{equation*} Then, applying the Gr\"onwall inequalities,  for all $t\in [0,T]$ we have
\begin{equation*}
\begin{split}
\Vert \vu(t,\cdot) \Vert^{2}_{L^2}+ \alpha^2 \Vert \vu(t,\cdot)\Vert^{2}_{\dot{H}^1} \leq & \left( \Vert \vu_0 \Vert^{2}_{L^2}+ \alpha^2 \Vert \vu_0\Vert^{2}_{\dot{H}^1} \right)\, e^{2t} + \int_{0}^{t} e^{2(t-s)} \left(- 2\nu \Vert \vu(s,\cdot)\Vert^{2}_{\dot{H}^{1}} - 2 \alpha^2\Vert \vu(s,\cdot)\Vert^{2}_{\dot{H}^{2}}\right) ds  \\
& + \int_{0}^{t} e^{2(t-s)} \left( 2 \Vert \fe \Vert^{2}_{L^2}+2 \alpha^2\Vert \fe \Vert^{2}_{\dot{H}^{1}}\right)ds,
\end{split}
\end{equation*}
hence we get
\begin{equation*}
\begin{split}
\Vert \vu(t,\cdot) \Vert^{2}_{L^2}+\alpha^2 \Vert \vu(t,\cdot)\Vert^{2}_{\dot{H}^1} \leq & \left(\Vert \vu_0 \Vert^{2}_{L^2}+ \alpha^2 \Vert \vu_0\Vert^{2}_{\dot{H}^1}\right)\, e^{2t} - \int_{0}^{t}  \left(2\nu \Vert \vu(s,\cdot)\Vert^{2}_{\dot{H}^{1}} + 2\alpha^2\Vert \vu(s,\cdot)\Vert^{2}_{\dot{H}^{2}}\right) ds  \\
& +  e^{2t} t  \left(  \Vert \fe \Vert^{2}_{L^2}+\alpha^2\Vert \fe \Vert^{2}_{\dot{H}^{1}}\right),
\end{split}
\end{equation*}
and we  thus obtain the following control in time
\begin{equation}\label{control-aux}
\begin{split}
\Vert \vu(t,\cdot) \Vert^{2}_{L^2} + \alpha^2 \Vert \vu(t,\cdot)\Vert^{2}_{\dot{H}^1} &+ \int_{0}^{t}  \left(\nu \Vert \vu(s,\cdot)\Vert^{2}_{\dot{H}^{1}} + \alpha^2\Vert \vu(s,\cdot)\Vert^{2}_{\dot{H}^{2}}\right) ds\\
\leq & \,\, e^{2t} \left( \Vert \vu_0 \Vert^{2}_{L^2}+ \alpha^2 \Vert \vu_0\Vert^{2}_{\dot{H}^1}\right)  +e^{2t}t\left( \Vert \fe \Vert^{2}_{L^2} + t\, \alpha^2\Vert \fe \Vert^{2}_{\dot{H}^{1}} \right).
\end{split}
\end{equation} which allows us to extend the local solution $\vu$ to the  interval $[0, +\infty[$. \\

The third step is to obtain the global energy equality (\ref{energ}). It directly follows by integrating the identity (\ref{Eq-dif-energ}) on the interval of time $[0,t]$. \\

The fourth step is to recover the pressure $p$ which is always related to the velocity $\vu$. Indeed, as $\vu$ verifies the equation (\ref{Eq-Aux}), and moreover, as we have $div(\vu)=0$ and $div(\fe)=0$, then we can write 
$$ \P\left( \partial_t \vu + \text{div}((-\alpha^2 \Delta +I_d)^{-1} (\vu \otimes \vu)) - \nu \Delta \vu  -\fe +\beta \vu \right) =0,$$ hence, by the well-known properties of the Leray projector $\P$ (see the Lemma $6.3$, page $118$ of the  book \cite{PGLR1}) there exists $p \in \mathcal{D}^{'}([0,+\infty[\times \Rt)$ such that we have 
\begin{equation*}
\partial_t \vu + \text{div}((-\alpha^2 \Delta +I_d)^{-1} (\vu \otimes \vu)) - \nu \Delta \vu  -\fe +\beta \vu=\vec{\nabla} p. 
\end{equation*}

Applying the divergence operator in each term of this identity, and moreover, denoting by  $\mathcal{R}_i = \frac{\partial_i}{\sqrt{-\Delta}}$ the  Riesz transforms,  we get the following identity 
\begin{equation}\label{Caract-pression}
p = \sum_{i=1}^{3}\sum_{j=1}^{3}  \mathcal{R}_{i}\mathcal{R}_i \left((-\alpha^2 \Delta + I_d)^{-1}(u_i u_j)\right).
\end{equation}
By this identity, the pressure term verifies $p \in L^{2}_{loc}([0,+\infty[, 
H^3(\Rt))$. Indeed, we shall quickly verify that  $ \vu \otimes \vu \in (L^{2}_{t})_{loc}\,H_{x}^{1}$. By Lemma \ref{LemmaTech2} we have $\vu \otimes \vu  \in (L^{2}_{t})_{loc}\dot{H}^{1}_{x}$. On the other hand, as $\vu \in (L^{\infty}_{t})_{loc}\,H^{1}_{x}$ then by the product laws in the non-homogeneous Sobolev spaces  (see always  the Lemma $7.3$, page $130$ of the book \cite{PGLR1}) we get $\vu \otimes \vu \in (L^{\infty}_{t})_{loc} \, H^{1/2}_{x}$, hence we have $\vu \otimes \vu \in (L^{2}_{t})_{loc} L^{2}_{x}$ and thereafter, joint with the fact that $\vu \otimes \vu  \in (L^{2}_{t})_{loc}\, \dot{H}^{1}_{x}$ we finally obtain $ \vu \otimes \vu \in (L^{2}_{t})_{loc}H^{1}_{x}$. \\   

The fifth and last step is to prove the uniqueness of solutions. So, let 
$(\vu_1, p_1) \in (L^{\infty}_{t})_{loc}H^{1}_{x} \cap (L^{2}_{t})_{loc}\dot{H}^{2}_{x} \times (L^{2}_{t})_{loc}H^{3}_{x}$ and  $(\vu_2, p_2) \in (L^{\infty}_{t})_{loc}H^{1}_{x} \cap (L^{2}_{t})_{loc}\dot{H}^{2}_{x} \times (L^{2}_{t})_{loc}H^{3}_{x}$  be two solutions of equation (\ref{Bardina}) arising from the initial data $\vu_{0,1}$ and $\vu_{0,2}$ respectively. We define $\vw=\vu_1-\vu_2$ and $q=p_1-p_2$, and then we get that the couple $(\vw, q)$ verifies the equation  
\begin{equation*}
\partial_t \vw + \left((\vw \cdot \vec{\nabla}) \vu_1+ (\vu_2\cdot \vec{\nabla})\vw\right)_\alpha- \nu \Delta \vw + \vec{\nabla} q= -\beta \vw, \quad \vw(0,\cdot)= \vu_{1}(0,\cdot)-\vu_2(0,\cdot).  
\end{equation*} Following the computations done in (\ref{Eq-dif-energ})  we have  
\begin{equation}\label{estim7}  
\begin{split}
 \frac{1}{2}\frac{d}{dt}\left(\Vert \vw(t,\cdot)\Vert^{2}_{L^2}+\alpha^2 \Vert \vw(t,\cdot)\Vert^{2}_{\dot{H}^1}\right)	+\nu \Vert \vw(t,\cdot)\Vert^{2}_{\dot{H}^1}+ &\alpha^2 \Vert \vw(t,\cdot)\Vert^{2}_{\dot{H}^2} =- \beta \left( \Vert \vw(t,\cdot)\Vert^{2}_{L^2}+\alpha^2 \Vert \vw(t,\cdot)\Vert^{2}_{\dot{H}^1} \right) \\
 & - \left\langle (\vw \cdot \vec{\nabla}) \vu_1(t,\cdot) , \vw(t,\cdot)\right\rangle_{\dot{H}^{-1}\times \dot{H}^1},
 \end{split}
\end{equation}	where we must estimate the term $\ds{\left\langle (\vw \cdot \vec{\nabla}) \vu_1(t,\cdot) , \vw(t,\cdot)\right\rangle_{\dot{H}^{-1}\times \dot{H}^1}}$. For this, we use the Hardy-Littlewood-Sobolev inequalities and the H\"older inequalities to write 
\begin{equation*}
\begin{split}
 \left\vert \left\langle (\vw \cdot \vec{\nabla}) \vu_1(t,\cdot) , \vw(t,\cdot)\right\rangle_{\dot{H}^{-1}\times \dot{H}^1} \right\vert & \leq c \Vert (\vw \cdot \vec{\nabla}) \vu_1(t,\cdot)\Vert_{\dot{H}^{-1}} \Vert \vw(t,\cdot)\Vert_{\dot{H}^1} \leq c \Vert (\vw \cdot \vec{\nabla}) \vu_1(t,\cdot)\Vert_{L^{6/5}} \Vert \vw(t,\cdot)\Vert_{\dot{H}^1}  \\
 & \leq c \Vert \vw(t,\cdot) \Vert_{L^2} \Vert \vec{\nabla} \otimes \vu_1 (t,\cdot) \Vert_{L^3} \, \Vert \vw(t,\cdot)\Vert_{\dot{H}^1} =(a).
\end{split} 
\end{equation*}  To treat the term $\ds{\Vert \vec{\nabla} \otimes \vu_1 (t,\cdot) \Vert_{L^3} }$ in expression (a),  we apply first the interpolation inequalities, and thereafter, by the Hardy-Littlewood-Sobolev we obtain
\begin{equation*}
\begin{split}
(a) \leq & c \,  \Vert \vw(t,\cdot) \Vert_{L^2}  \, \Vert \vec{\nabla} \otimes \vu_1(t,\cdot) \Vert^{1/2}_{L^2}\Vert \vec{\nabla} \otimes \vu_1(t,\cdot) \Vert^{1/2}_{L^6}\, \Vert \vw(t,\cdot)\Vert_{\dot{H}^1}  \\
\leq &  c\,  \Vert \vw(t,\cdot) \Vert_{L^2}\, \Vert \vec{\nabla} \otimes \vu_1(t,\cdot)\Vert^{1/2}_{L^2} \Vert \vec{\nabla} \otimes \vu_1(t,\cdot) \Vert^{1/2}_{\dot{H}^1}\, \Vert \vw(t,\cdot)\Vert_{\dot{H}^1}  \\
 \leq & c \, \Vert \vw(t,\cdot) \Vert_{L^2}\, \Vert \vu_1(t,\cdot) \Vert^{1/2}_{\dot{H}^1} \Vert \vu_1(t,\cdot) \Vert^{1/2}_{\dot{H}^2} \, \Vert \vw(t,\cdot) \Vert_{\dot{H}^1}\\
 \leq & \, \frac{c}{\nu} \Vert \vw(t,\cdot) \Vert^{2}_{L^2} \, \Vert \vu_1(t,\cdot) \Vert_{\dot{H}^1} \Vert \vu_1(t,\cdot) \Vert_{\dot{H}^2} + \frac{\nu}{2}  \Vert \vw(t,\cdot) \Vert^{2}_{\dot{H}^1} \\
 \leq & \,  \frac{c}{\nu} \left( \Vert \vw(t,\cdot) \Vert^{2}_{L^2} +\alpha^2 \Vert \vw(t,\cdot) \Vert^{2}_{\dot{H}^1}\right) \, \Vert \vu_1(t,\cdot) \Vert_{\dot{H}^1} \Vert \vu_1(t,\cdot) \Vert_{\dot{H}^2}+ \frac{\nu}{2}  \Vert \vw(t,\cdot) \Vert^{2}_{\dot{H}^1}. 
\end{split}
\end{equation*}
With this estimate at hand, we get back to (\ref{estim7}) where we have 
\begin{equation*}
\begin{split}
 \frac{1}{2}\frac{d}{dt}\left(\Vert \vw(t,\cdot)\Vert^{2}_{L^2}+\alpha^2 \Vert \vw(t,\cdot)\Vert^{2}_{\dot{H}^1}\right)&	+\nu \Vert \vw(t,\cdot)\Vert^{2}_{\dot{H}^1}+\alpha^2 \Vert \vw(t,\cdot)\Vert^{2}_{\dot{H}^2} \leq - \beta \left( \Vert \vw(t,\cdot)\Vert^{2}_{L^2}+\alpha^2 \Vert \vw(t,\cdot)\Vert^{2}_{\dot{H}^1} \right) \\
&+  \frac{c}{\nu} \left( \Vert \vw(t,\cdot) \Vert^{2}_{L^2} +\alpha^2 \Vert \vw(t,\cdot) \Vert^{2}_{\dot{H}^1}\right) \, \Vert \vu_1(t,\cdot) \Vert_{\dot{H}^1} \Vert \vu_1(t,\cdot) \Vert_{\dot{H}^2}+ \frac{\nu}{2}  \Vert \vw(t,\cdot) \Vert^{2}_{\dot{H}^1}
\end{split},
\end{equation*} hence we get 
\begin{equation*}
 \frac{1}{2}\frac{d}{dt}\left(\Vert \vw(t,\cdot)\Vert^{2}_{L^2}+\alpha^2 \Vert \vw(t,\cdot)\Vert^{2}_{\dot{H}^1}\right) \leq \frac{c}{\nu} \left( \Vert \vw(t,\cdot) \Vert^{2}_{L^2} +\alpha^2 \Vert \vw(t,\cdot) \Vert^{2}_{\dot{H}^1}\right) \, \Vert \vu_1(t,\cdot) \Vert_{\dot{H}^1} \Vert \vu_1(t,\cdot) \Vert_{\dot{H}^2}.
\end{equation*}
Thus, applying the Gr\"onwall inequalities we have 
\begin{equation*}
\Vert \vw(t,\cdot)\Vert^{2}_{L^2}+\alpha^2 \Vert \vw(t,\cdot)\Vert^{2}_{\dot{H}^1} \leq \left( \Vert \vw(0,\cdot) \Vert^{2}_{L^2}+\alpha^2 \Vert \vw(0,\cdot)\Vert^{2}_{\dot{H}^1} \right) \, e^{ \ds{\frac{c}{\nu} \int_{0}^{t} \Vert \vu_1(s,\cdot) \Vert_{\dot{H}^1} \Vert \vu_1(s,\cdot) \Vert_{\dot{H}^2}\, ds}}.
\end{equation*}
Moreover, using (\ref{control-aux}) the term $\ds{\frac{c}{\nu}\int_{0}^{t} \Vert \vu_1(s,\cdot) \Vert_{\dot{H}^1} \Vert \vu_1(s,\cdot) \Vert_{\dot{H}^2}\, ds}$ is estimated as follows: 
\begin{equation*}
\begin{split}
& \frac{c}{\nu}\int_{0}^{t} \Vert \vu_1(s,\cdot) \Vert_{\dot{H}^1} \Vert \vu_1(s,\cdot) \Vert_{\dot{H}^2}\, ds \leq  c_{\nu, \alpha} \left( \nu \int_{0}^{t}\Vert \vu_1(s,\cdot)\Vert^{2}_{\dot{H}^1}ds+ \alpha^2 \int_{0}^{t} \Vert \vu_1(s,\cdot)\Vert^{2}_{\dot{H}^2} ds \right)\\
\leq & c_{\nu,\alpha} \left(e^{2t} \left( \Vert \vu_{0,1} \Vert^{2}+ \alpha^2 \Vert \vu_{0,1}\Vert^{2}_{\dot{H}^1}\right)  +e^{2t}t\left( \Vert \fe \Vert^{2}_{L^2} + t\, \alpha^2\Vert \fe \Vert^{2}_{\dot{H}^{1}} \right)\right)=c_2(\alpha,\nu,\fe,\vu_{0,1},t).  
 \end{split}
\end{equation*}
With this estimate, we can write the following inequality: 
\begin{equation}\label{Estim-Unicidad}
\Vert \vw(t,\cdot)\Vert^{2}_{L^2}+\alpha^2 \Vert \vw(t,\cdot)\Vert^{2}_{\dot{H}^1} \leq \left( \Vert \vw(0,\cdot) \Vert^{2}_{L^2}+\alpha^2 \Vert \vw(0,\cdot)\Vert^{2}_{\dot{H}^1} \right) \, e^{c_2(\alpha,\nu,\fe, \vu_{0,1},t)}. 
\end{equation}   Uniqueness of solutions directly follows from this inequality. Theorem \ref{Th-WP} is now proven. \finpv

\subsection*{Proof of Proposition \ref{Prop:time-control}} 
\begin{enumerate}
    \item[1)] The control (\ref{control-tiempo}) directly follows from the identity (\ref{Eq-dif-energ}). Indeed, by the Cauchy-Schwarz inequalities, and as $\ds{-\nu \Vert \vu(t,\cdot) \Vert^{2}_{\dot{H}^1}-\alpha^2 \Vert \vu(t,\cdot)\Vert^{2}_{\dot{H}^2}}$ is negative quantity,   we write 
\begin{equation*}
\begin{split}
&\frac{1}{2} \frac{d}{dt}\left(\Vert \vu(t,\cdot) \Vert^{2}+ \alpha^2 \Vert \vu(t,\cdot)\Vert^{2}_{\dot{H}^1}\right) \\
\leq & -\nu \Vert \vu(t,\cdot) \Vert^{2}_{\dot{H}^1}-\alpha^2 \Vert \vu(t,\cdot)\Vert^{2}_{\dot{H}^2} +  \Vert \fe \Vert_{L^{2}} \Vert \vu(t,\cdot)\Vert_{L^2} +\alpha^2 \Vert \fe \Vert_{\dot{H}^{1}} \Vert \vu(t,\cdot)\Vert_{\dot{H}^{1}} -\beta \left( \Vert \vu(t,\cdot)\Vert^{2}_{L^2}+\alpha^2 \Vert\vu(t,\cdot)\Vert^{2}_{\dot{H}^1}\right)\\
\leq & -\nu \Vert \vu(t,\cdot) \Vert^{2}_{\dot{H}^1}-\alpha^2 \Vert \vu(t,\cdot)\Vert^{2}_{\dot{H}^2} + \frac{2}{\beta} \Vert \fe \Vert^{2}_{L^2}+ \frac{\beta}{2}\Vert \vu(t,\cdot)\Vert^{2}_{L^2} + \frac{2 \alpha^2}{\beta}\Vert \fe \Vert^{2}_{\dot{H}^{1}}+ \frac{\beta \alpha^2}{2} \Vert \vu(t,\cdot)\Vert^{2}_{\dot{H}^1} \\
& - \beta \left( \Vert \vu(t,\cdot)\Vert^{2}_{L^2}+\alpha^2 \Vert\vu(t,\cdot)\Vert^{2}_{\dot{H}^1}\right)\\
\leq & \frac{2}{\beta}\left(\Vert \fe \Vert^{2}_{L^2}+\alpha^2\Vert \fe \Vert^{2}_{\dot{H}^1} \right) -\frac{\beta}{2} \left( \Vert \vu(t,\cdot)\Vert^{2}_{L^2}+\alpha^2 \Vert\vu(t,\cdot)\Vert^{2}_{\dot{H}^1}\right).
\end{split}
\end{equation*} Then, applying the Gr\"onwall inequalities
we have 
\begin{equation*}
 \Vert \vu(t,\cdot) \Vert^{2}+ \alpha^2 \Vert \vu(t,\cdot)\Vert^{2}_{\dot{H}^1} \leq\left( \Vert \vu_0 \Vert^{2}+ \alpha^2 \Vert \vu_0\Vert^{2}_{\dot{H}^1} \right)e^{-\beta t}+ \frac{4}{\beta} \int_{0}^{t}e^{-\beta(t-s)}\left( \Vert \fe \Vert^{2}_{L^2}+\alpha^2\Vert \fe \Vert^{2}_{\dot{H}^1} \right) ds,  
\end{equation*} hence, the desired control (\ref{control-tiempo}) follows. 
\item[2)] Always by the identity (\ref{Eq-dif-energ}) and   the Cauchy-Schwarz inequalities, we have 
\begin{equation*}
\begin{split}
& \frac{1}{2} \frac{d}{dt}\left(\Vert \vu(t,\cdot) \Vert^{2}+ \alpha^2 \Vert \vu(t,\cdot)\Vert^{2}_{\dot{H}^1}\right)\\
\leq & -\nu \Vert \vu(t,\cdot) \Vert^{2}_{\dot{H}^1}-\alpha^2 \Vert \vu(t,\cdot)\Vert^{2}_{\dot{H}^2} +  \Vert \fe \Vert_{L^{2}} \Vert \vu(t,\cdot)\Vert_{L^2}  +\alpha^2 \Vert \fe \Vert_{\dot{H}^{1}} \Vert \vu(t,\cdot)\Vert_{\dot{H}^{1}} -\beta \left( \Vert \vu(t,\cdot)\Vert^{2}_{L^2}+\alpha^2 \Vert\vu(t,\cdot)\Vert^{2}_{\dot{H}^1}\right)\\
\leq & -\nu \Vert \vu(t,\cdot) \Vert^{2}_{\dot{H}^1}-\alpha^2 \Vert \vu(t,\cdot)\Vert^{2}_{\dot{H}^2} + \frac{2}{\beta}\Vert \fe \Vert^{2}_{L^2}+ \frac{\beta}{2}\Vert \vu(t,\cdot)\Vert^{2}_{L^2}  + \frac{2 \alpha^2}{\beta} \Vert \fe \Vert^{2}_{\dot{H}^1}+ \frac{\beta \alpha^2}{2}\Vert \vu(t,\cdot)\Vert^{2}_{\dot{H}^1}\\
& -\beta \left( \Vert \vu(t,\cdot)\Vert^{2}_{L^2}+\alpha^2 \Vert\vu(t,\cdot)\Vert^{2}_{\dot{H}^1}\right)\\
\leq &-\nu \Vert \vu(t,\cdot) \Vert^{2}_{\dot{H}^1}-\alpha^2 \Vert \vu(t,\cdot)\Vert^{2}_{\dot{H}^2} + \frac{2}{\beta}\left( \Vert \fe \Vert^{2}_{L^2}+\alpha^2 \Vert \fe \Vert^{2}_{\dot{H}^1}\right) - \frac{\beta}{2}\left( \Vert \vu(t,\cdot) \Vert^{2}+ \alpha^2 \Vert \vu(t,\cdot)\Vert^{2}_{\dot{H}^1} \right) \\
 \leq &-\nu \Vert \vu(t,\cdot) \Vert^{2}_{\dot{H}^1}-\alpha^2 \Vert \vu(t,\cdot)\Vert^{2}_{\dot{H}^2} + \frac{2}{\beta}\left( \Vert \fe \Vert^{2}_{L^2}+\alpha^2 \Vert \fe \Vert^{2}_{\dot{H}^1}\right).
\end{split}
\end{equation*}
Then, integrating in the interval $[t,t+T]$ we get 
\begin{equation*}
\begin{split}
\Vert  \vu(t+T,\cdot) \Vert^{2}_{L^2} + \alpha^2 \Vert \vu(t+T,\cdot)\Vert^{2}_{\dot{H}^1} \leq &-\nu \int_{t}^{t+T} \Vert \vu(s,\cdot) \Vert^{2}_{\dot{H}^1}ds - \alpha^2 \int_{t}^{t+T}  \Vert \vu(s,\cdot)\Vert^{2}_{\dot{H}^2} ds \\
&+ \frac{2 T}{\beta}\left( \Vert \fe \Vert^{2}_{L^2}+\alpha^2 \Vert \fe \Vert^{2}_{\dot{H}^1}\right) + \Vert  \vu(t,\cdot) \Vert^{2}_{L^2} + \alpha^2 \Vert \vu(t,\cdot)\Vert^{2}_{\dot{H}^1},
\end{split}    
\end{equation*}hence we have the desired control (\ref{control-tiempo2}) \finpv
\end{enumerate}

\section{The global attractor}\label{Sec:GlobalAttractor}
%%%
To prove the existence of a global attractor associated to the equation (\ref{Bardina}), stated in Theorem \ref{Th-Atractor}, we will use some  results arising from the theory of dynamical systems which, for completeness of the paper and the reader's convenience, we shall state below.\\

We will start by recalling that  for the parameter $\alpha >0$ fixed, we consider the Banach space
$\mathcal{H}^{1}_{\alpha}(\Rt)=\left\{ \vg  \in H^1(\Rt): \,\, div(\vg)=0, \quad \Vert \vg \Vert^{2}_{\mathcal{H}^{1}_{\alpha}}=\Vert \vg \Vert^{2}_{L^2}+\alpha^{2} \Vert \vg \Vert^{2}_{\dot{H}^1} <+\infty\right\}$. Thereafter, for $t\geq 0$ let $S(t): \mathcal{H}^{1}_{\alpha}(\Rt)\to \mathcal{H}^{1}_{\alpha}(\Rt)$ be the semigruop associated with equation (\ref{Bardina}) and defined in (\ref{def-semigrupo}). \\

Before to state a general result leading the existence of a global attractor for the semigroup $S(t)$, we introduce the following definitions  we shall need later. 

\begin{Definition}[Absorbig set]\label{def-Abosorbing-set}
A closed set $\mathcal{B} \subset \mathcal{H}^{1}_{\alpha}(\Rt)$ is an absorbing set for the semigroup $S(t)$ if for every bounded set $B \subset \mathcal{H}^{1}_{\alpha}(\Rt)$,
there exists a time $T = T(B)>0$ such that, for all $t>T$ we have 
$S(t)B \subset \mathcal{B}$.
\end{Definition}

In this definition, it is worth mention  we use the notation 
$\ds{S(t)B= \left\{ S(t)\vu_0: \vu_0 \in B \right\}}$. 
\begin{Definition}[Semigroup asymptotically compact]\label{def-asymptotically-compact}
The semigroup $S(t)$ is  asymptotically compact 
if for any bounded sequence $(\vu_{0,n})_{n\in \mathbb{N}}$ in $\mathcal{H}^{1}_{\alpha}(\Rt)$, and moreover,  for any sequence of times $(t_n)_{n \in \mathbb{N}}$ such that  $t_n \to \infty$, the sequence $(S(t_n)\vu_{0,n})_{n\in \mathbb{N}}$  is precompact in $\mathcal{H}^{1}_{\alpha}(\Rt)$.
\end{Definition}

Once we have introduced these definitions, we are able to state the following result on the existence of a global attractor. For a proof of this fact see \cite{Temam}. 
\begin{Proposition}\label{Prop:Existence-Attractor}
 Suppose that:
 \begin{enumerate}
     \item The semigroup $S(t)$ has a bounded and closed absorbing set $\mathcal{B} \subset \mathcal{H}^{1}_{\alpha}(\Rt)$ given in Definition \ref{def-Abosorbing-set}.
     \item The semigroup $S(t)$ is asymptotically compact in the sense of Definition \ref{def-asymptotically-compact}.
\item  For every fixed $t \geq 0$ the map $S(t) : \mathcal{B} \to \mathcal{H}^{1}_{\alpha}(\Rt)$ is continuous.
 \end{enumerate}
Then, the semigroup $S(t)$ has a  global attractor $\mathcal{A}_{\fe} \subset \mathcal{H}^{1}_{\alpha}(\Rt)$ given in Definition \ref{Def:atractor}. 
\end{Proposition} 
%%
%% Esto deberia ir en la presentacion de los resultados (a ver lego) 
%Moreover, the attractor $\mathcal{A}$ has the following structure:
%\begin{equation*}
%\mathcal{A}=\mathcal{K}\Big|_{t=0}, 
%\end{equation*}
%where $\mathcal{K}\subset L^{\infty}$
%R, H) is the set of complete trajectories u : R → H of
%semigroup S(t) which are defined for all t ∈ R and bounded.
%%

\subsection*{Proof of Theorem \ref{Th-Atractor}} 
We will prove that the semigroup $S(t)$ associated to equation (\ref{Bardina}) and defined in (\ref{def-semigrupo}) verify  the points $1.$, $2.$ and $3.$ in Proposition \ref{Prop:Existence-Attractor}. We start by verifying the point $1$ with the following result.
\begin{Lemme}\label{LemaTech2-Absorbing set}  Let $\ds{\mathcal{B}=\left\{ \vu_0 \in \mathcal{H}^{1}_{\alpha}(\Rt): \Vert \vu_0 \Vert^{2}_{\mathcal{H}^{1}_{\alpha}} \leq \frac{8}{\beta^2} \Vert \fe \Vert^{2}_{\mathcal{H}^{1}_{\alpha}}\right\}}$. Then, $\mathcal{B}$ is a absorbing set in the sense of Definition \ref{def-Abosorbing-set}.
\end{Lemme}
\pv We observe that $\mathcal{B}$ is a bounded and closed set in $\mathcal{H}^{1}_{\alpha}(\Rt)$ and we will prove that $\mathcal{B}$ is moreover an absorbing set. Indeed, let $\ds{B \subset \mathcal{H}^{1}_{\alpha}(\Rt)}$ be a bounded set. Then, for $R>0$ (large enough) we have $\Vert \vu_0 \Vert^{2}_{\mathcal{H}^{1}_{\alpha}}\leq R^2$ for all $\vu_0 \in B$. On the other hand, by point $1)$ in Proposition  \ref{Prop:time-control}, for all  $\vu_0 \in \mathcal{B}$ we have 
$$ \Vert S(t)\vu_0 \Vert^{2}_{\mathcal{H}^{1}_{\alpha}} \leq \Vert \vu_0 \Vert^{2}_{\mathcal{H}^{1}_{\alpha}}e^{-\beta t} + \frac{4}{\beta^2} \Vert \fe \Vert^{2}_{\mathcal{H}^{1}_{\alpha}}\leq R^2e^{-\beta t} + \frac{4}{\beta^2} \Vert \fe \Vert^{2}_{\mathcal{H}^{1}_{\alpha}}.$$
Here, we set a time $T=T(B)>0$ such that for all $t>T$ we have  $\ds{R^2e^{-\beta t} \leq \frac{4}{\beta^2} \Vert \fe \Vert^{2}_{\mathcal{H}^{1}_{\alpha}}}$, and then,  and for all $\vu_0 \in \mathcal{B}$ we get $S(t)\vu_0 \in \mathcal{B}$. \finpv

We verify now the point $2$.   Let $(\vu_{0,n})_{n\in \mathbb{N}}$ be a bounded sequence in $\mathcal{H}^{1}_{\alpha}(\Rt)$, and moreover, let $(t_n)_{n\in \mathbb{N}}$ be a sequence of positive times such that  $t_n \to +\infty$. We must show that the sequence $(S(t_n)\vu_{0,n})_{n\in \mathbb{N}}$ is precompact in $\mathcal{H}^{1}_{\alpha}(\Rt)$ and for this, we will use a energy method. For each $n \in \mathbb{N}$, we consider the following Cauchy problem for the equation (\ref{Bardina}): 
\begin{equation}\label{Bardina-n}
\left\{ \begin{array}{ll}\vspace{2mm}
\partial_t \vu_n + \text{div} \,( (\vu_n \otimes \vu_n)_\alpha) - \nu \Delta \vu_n + \vec{\nabla} p_n = \fe -\beta \vu_n, \qquad \text{div}(\vu_n)=0, \\ 
\vu_n(-t_n,\cdot)= \vu_{0,n}.
\end{array} 
 \right.
\end{equation}
Let  $\vu_n : [-t_n, +\infty[\times \Rt \to \Rt$ and $p_n:  [-t_n, +\infty[\times \Rt \to \R$, be the unique solution of this equation given by Theorem \ref{Th-WP}, where it verifies $\ds{\vu_n \in L^{\infty}([-t_n,+\infty[, H^1(\Rt)) \cap L^{2}_{loc}([-t_n, +\infty[, \dot{H}^2(\Rt))}$ and $p_n \in  L^{2}_{loc}([-t_n,+\infty[, 
H^3(\Rt))$. \\        

Now, by uniqueness of solution $\vu_n$, and moreover, by definition of the semigroup $S(t)$ associated with the equation (\ref{Bardina-n}),  defined in (\ref{def-semigrupo}), for all $n\in \mathbb{N}$ we have the identity  $\ds{S(t_n)\vu_{0,n}= \vu_{n}(0, \cdot)}$ and thus, it is enough to  verify  that the sequence $(\vu_{n}(0,\cdot))_{n\in \mathbb{N}}$ is precompact in $\mathcal{H}^{1}_{\alpha}(\Rt)$. For this, our general strategy is the following one.  First, we   prove the existence of  a solution  $(\vv,q)$  to the equation (\ref{Bardina-enternal}), called the  eternal  given in Definition \ref{defi-eternal-sol}. Then, we will show that the sequence $(\vu_{n}(0,\cdot))_{n\in \mathbb{N}}$ strongly converges (via a sub-sequence) to $v(0,\cdot)$   in the space $H^{1}_{\alpha}(\Rt)$. \\

We start by the construction of an eternal  solution $(\vv, q)$. 
\begin{Proposition}\label{Prop:eternal-sol} There exists a couple of functions $(\vv,q)$, with $\vv \in \ds{L^{\infty}_{loc}(\R, H^1(\Rt)) \cap L^{2}_{loc}(\R, \dot{H}^2(\Rt))}$ and  $q \in L^{2}_{loc}(\R, \dot{H}^3(\Rt))$, which is a weak solution of the equation (\ref{Bardina-enternal}). 
\end{Proposition}
\pv  This solution will be obtained as the limit of the solutions $\vu_n : [-t_n, +\infty[\times \Rt \to \Rt$  and $p_n:[-t_n,+\infty[\times \Rt \to \R$ of equations (\ref{Bardina-n}) when $n \to +\infty$. \\

We observe  that by point $1)$ in Proposition  \ref{Prop:time-control}, for all $n\in \mathbb{N}$  and for all  $t \geq -t_n$, we have
\begin{equation}\label{control-tiempo-n}
\Vert \vu_n(t,\cdot)\Vert^{2}_{L^2}+\alpha^2\Vert \vu_n(t,\cdot)\Vert^{2}_{\dot{H}^1} \leq \left( \Vert \vu_{0,n} \Vert^{2}_{L^2}+\alpha^2\Vert \vu_{0,n} \Vert^{2}_{\dot{H}^1}\right)e^{-\beta\, (t+t_n)} + \frac{4}{\beta^2} \left( \Vert \fe \Vert^{2}_{L^2}+\alpha^2\Vert \fe \Vert^{2}_{\dot{H}^1}\right).
\end{equation}
Moreover, by point $2)$ in Proposition \ref{Prop:time-control}, for all  $t\geq -t_n$ and $T\geq 0$, we have 
\begin{equation}\label{control-tiempo2-n}
\begin{split}
\nu \int_{t}^{t+T} \Vert \vu_n(s,\cdot)\Vert^{2}_{\dot{H}^1} ds + \alpha^2 \int_{t}^{t+T} \Vert \vu_n(s,\cdot)\Vert^{2}_{\dot{H}^2}ds \leq & \frac{2 T}{\beta}\left( \Vert \fe \Vert^{2}_{L^2}+\alpha^2 \Vert \fe \Vert^{2}_{\dot{H}^1}\right) \\
& + \Vert \vu_n(t,\cdot)\Vert^{2}_{L^2}+\alpha^2\Vert \vu_n(t,\cdot)\Vert^{2}_{\dot{H}^1}. 
\end{split} 
\end{equation}

By estimate (\ref{control-tiempo-n}), and recalling that the sequence $(\vu_{0,n})_{n\in \mathbb{N}}$ is bounded in $\mathcal{H}^{1}_{\alpha}(\Rt)$,  we can write 
\begin{equation}\label{estim1-n}
\begin{split}
\sup_{n \in \mathbb{N}} \  \sup_{t \geq -t_n} \left( \Vert \vu_n(t,\cdot)\Vert^{2}_{L^2}+\alpha^2\Vert \vu_n(t,\cdot)\Vert^{2}_{\dot{H}^1}  \right) \leq  & R^2  e^{-\beta\, (t+t_n)} + \frac{4}{\beta^2} \left( \Vert \fe \Vert^{2}_{L^2}+\alpha^2\Vert \fe \Vert^{2}_{\dot{H}^1}\right)  \\
\leq & R^2 + \frac{4}{\beta^2} \left( \Vert \fe \Vert^{2}_{L^2}+\alpha^2\Vert \fe \Vert^{2}_{\dot{H}^1}\right). 
\end{split}
\end{equation}
Moreover, by estimate (\ref{control-tiempo2-n}) (with $T=1$) we get
\begin{equation}\label{estim2-n}
\begin{split}
& \sup_{n \in \mathbb{N}} \ 
\sup_{ t \geq -t_n} \left( \nu \int_{t}^{t+1} \Vert \vu_n(s,\cdot)\Vert^{2}_{\dot{H}^1} ds + \alpha^2 \int_{t}^{t+1} \Vert \vu_n(s,\cdot)\Vert^{2}_{\dot{H}^2}ds \right) \\
\leq &   \frac{2}{\beta}\left( \Vert \fe \Vert^{2}_{L^2}+\alpha^2 \Vert \fe \Vert^{2}_{\dot{H}^1}\right)  + \Vert \vu_n(t,\cdot)\Vert^{2}_{L^2}+\alpha^2\Vert \vu_n(t,\cdot)\Vert^{2}_{\dot{H}^1}\\
 \leq &  \frac{2}{\beta}\left( \Vert \fe \Vert^{2}_{L^2}+\alpha^2 \Vert \fe \Vert^{2}_{\dot{H}^1}\right) + R^2 + \frac{4}{\beta^2} \left( \Vert \fe \Vert^{2}_{L^2}+\alpha^2\Vert \fe \Vert^{2}_{\dot{H}^1}\right). 
\end{split}    
\end{equation} 

Then, by estimates (\ref{estim1-n}) and (\ref{estim2-n}) and the Banach-Alaoglu theorem, there exists $\vv \in L^{\infty}(\R, H^1(\Rt)) \cap L^{2}_{loc}(\R, \dot{H}^2(\Rt))$ such that the sequence $(\vu_n)_{n\in \mathbb{N}}$ converges (via s sub-sequence) to $\vv$ in the weak-* topology of the spaces $L^{\infty}([-\tau,\tau], H^1(\Rt))$, $L^{2}([-\tau, \tau], \dot{H}^1(\Rt))$  and $L^{2}([-\tau, \tau], \dot{H}^2(\Rt))$, for all $\tau>0$. \\

It remains to prove that the sequence of pressure terms $(p_n)_{n \in \mathbb{N}}$ converges to a limit $q: \R\times \Rt \to \R$. Recall that  by identity  (\ref{Caract-pression}), for all $n \in \mathbb{N}$, we have:  $\ds{p_n =\sum_{i=1}^{3}\sum_{j=1}^{3}  \mathcal{R}_{i}\mathcal{R}_i \left((-\alpha^2 \Delta + I_d)^{-1}(u_{n,i} u_{n,j})\right)}$, hence, for all $t\geq -t_n$ and $T\geq 0$, we get
\begin{equation*}
 \int_{t}^{t+T}\Vert p_n(s,\cdot) \Vert^{2}_{H^3} ds \leq  c_\alpha  \int_{t}^{t+T} \Vert \vu_n \otimes \vu_n (s,\cdot)\Vert^{2}_{H^1}ds \leq   c_\alpha  \int_{t}^{t+T}  \Vert  \vu_n \otimes \vu_n (s,\cdot)\Vert^{2}_{L^2} ds\, + c_\alpha   \int_{t}^{t+T} \Vert  \vu_n \otimes \vu_n (s,\cdot)\Vert^{2}_{\dot{H}^1}ds,
\end{equation*}
where we must study each term in the right side. For the first term, by the product laws in the non-homogeneous Sobolev spaces (Lemma $7.3$, page $130$ of \cite{PGLR1}), and moreover, by estimate (\ref{control-tiempo-n}) we write
\begin{equation*}
\begin{split}
c_\alpha  \int_{t}^{t+T} \Vert  \vu_n \otimes \vu_n (s,\cdot)\Vert^{2}_{L^2} ds\leq & c_\alpha  \int_{t}^{t+T} \Vert \vu_n \otimes \vu_n(s,\cdot) \Vert^{2}_{H^{1/2}} ds \leq c_\alpha  \int_{t}^{t+T}\Vert \vu_{n}(s,\cdot) \Vert^{4}_{H^{1}} ds\\
\leq & c_\alpha  \int_{t}^{t+T}\left( \Vert \vu_n(s,\cdot)\Vert^{2}_{L^2}+ \Vert \vu_n(s,\cdot)\Vert^{2}_{\dot{H}^1}\right)^{2} ds \\
\leq & c_\alpha  \left(\max\left(1, 1/\alpha^2 \right) \right)^2  \int_{t}^{t+T}\left( \Vert \vu_n(s,\cdot)\Vert^{2}_{L^2}+ \alpha^2\Vert \vu_n(s,\cdot)\Vert^{2}_{\dot{H}^1} \right)^2  ds\\
\leq & C_\alpha\, T\left( R^2 + \frac{4}{\beta^2} \left(\Vert \fe \Vert^{2}_{L^2}+ \alpha^2\Vert \fe \Vert^{2}_{\dot{H}^1}\right) \right)^2.
\end{split}
\end{equation*}
We estimate now the second term. Recall first  that we have $\Vert \cdot \Vert_{E_T}= \Vert \cdot \Vert_{L^{\infty}_{t}H^{1}_{x}}+\Vert \cdot \Vert_{L^{2}_{t}\dot{H}^{2}_{x}}$, and then,  by  Lemma \ref{LemmaTech2} we get
\begin{equation*}
\begin{split}
 c_\alpha   \int_{t}^{t+T} & \Vert  \vu_n \otimes \vu_n (s,\cdot)\Vert^{2}_{\dot{H}^1}ds  \leq  c_\alpha T^{1/2} \Vert \vu_n \Vert^{2}_{E_T} \leq c_\alpha\, T^{1/2} \left(  \sup_{t\leq s \leq t+T} \Vert \vu_n(s,\cdot)\Vert^{2}_{H^1}+ \int_{t}^{t+T} \Vert \vu_n(s,\cdot)\Vert^{2}_{\dot{H}^2} ds\right) \\
 \leq & c_{\alpha} \left(\max\left(1, 1/\alpha^2 \right) \right)^2  \, T^{1/2} \left(  \sup_{t\leq s \leq t+T} \left( \Vert \vu_n(s,\cdot)\Vert^{2}_{L^2}+\alpha^2\Vert \vu_n(s,\cdot)\Vert^{2}_{\dot{H^1}} \right)+ \alpha^2 \int_{t}^{t+T} \Vert \vu_n(s,\cdot)\Vert^{2}_{\dot{H}^2} ds\right)=(a). 
\end{split}    
\end{equation*}
Thereafter, by estimates (\ref{control-tiempo-n}) and (\ref{control-tiempo2-n}) we can write 
\begin{equation*}
\begin{split}
 (a)\leq C_\alpha T^{1/2}\left( 2R^2 + \frac{8}{\beta^2} \left( \Vert \fe \Vert^{2}_{L^2}+\alpha^2\Vert \fe \Vert^{2}_{\dot{H}^1}\right)+ \frac{2}{\beta}\left( \Vert \fe \Vert^{2}_{L^2}+\alpha^2 \Vert \fe \Vert^{2}_{\dot{H}^1}\right) \right).   
\end{split}    
\end{equation*}

Once we have studied each term,  gathering  the estimates obtained  we finally get 
\begin{equation*}
\begin{split}
 \int_{t}^{t+T}\Vert p_n(s,\cdot) \Vert^{2}_{H^3} ds \leq & C_\alpha \, T\left( R^2 + \frac{4}{\beta^2} \left(\Vert \fe \Vert^{2}_{L^2}+ \alpha^2\Vert \fe \Vert^{2}_{\dot{H}^1}\right) \right)^2\\
 & + C_\alpha \, T^{1/2}\left( 2R^2 + \frac{8}{\beta^2} \left( \Vert \fe \Vert^{2}_{L^2}+\alpha^2\Vert \fe \Vert^{2}_{\dot{H}^1}\right)+ \frac{2}{\beta}\left( \Vert \fe \Vert^{2}_{L^2}+\alpha^2 \Vert \fe \Vert^{2}_{\dot{H}^1}\right) \right).
\end{split}    
\end{equation*}
Hence, taking $T=1$ we can write 
\begin{equation*}
\begin{split}
\sup_{n\in \mathbb{N}} \sup_{ t \geq -t_n} \,  \int_{t}^{t+T}\Vert p_n(s,\cdot) \Vert^{2}_{H^3} ds \leq C_\alpha \left( 3R^2 + \left(\frac{12}{\beta^2}+ \frac{2}{\beta} \right) \left( \Vert \fe \Vert^{2}_{L^2}+\alpha^2\Vert \fe \Vert^{2}_{\dot{H}^1}\right)\right),   
\end{split}    
\end{equation*}
and always by the Banach-Alaoglu theorem there exists $q \in L^{2}_{loc}(\R, H^3(\Rt))$, such that the sequence $(p_n)_{n\in \mathbb{N}}$ converges (through a  a subsequence) to a limit $q$ in the weak-* topology of the spaces $L^{2}([-\tau, \tau], H^3(\Rt))$, for all $\tau>0$.\\ 

We want to show that the limit $(\vv,q)$ verifies the equation \eqref{Bardina-enternal} in the distributional  sense.  From equation (\ref{Bardina-n}) we can write 
\begin{equation*}
\begin{split}
&   \partial_t (-\alpha^2 \Delta + Id ) \vu_n +   \text{div} \,(\vu_n \otimes \vu_n) -   \nu \Delta (-\alpha^2 \Delta + Id )\vu_n +  \,  \vec{\nabla}(-\alpha^2 \Delta + Id )  p_n \\
= &  (-\alpha^2 \Delta + Id )\, \fe -\beta (-\alpha^2 \Delta + Id )\, \vu_n,  
\end{split}   \qquad div(\vu_n)=0,
\end{equation*} 
and then, it is enough to prove that the non-linear terms $\text{div} \,(\vu_n \otimes \vu_n)$  converge to $\text{div} \,(\vv \otimes \vv)$ in the sense of distributions. For this, we will use the following Rellich-Lions  lemma. For a proof see the Theorem $12.1$, page $349$, of \cite{PGLR1}.
\begin{Lemme}\label{Lemma:RL}
 Let  $(g_n)_{n\in \mathbb{N}}$ be a sequence of measurable functions on $\R\times \Rt$, such that, for every $\varphi \in \mathcal{C}^{\infty}_{0}(\R\times \Rt)$, we have:
\begin{enumerate}
\item  for some  positive $\gamma>0$,  
$\ds{\sup_{n\in \mathbb{N}} \Vert  \varphi g_n\Vert_{L^{2}_{t} H^{\gamma}_{x}}<+\infty}$, and 
\item for some negative $\sigma <0$, 
$\ds{\sup_{n \in \mathbb{N}} \Vert \varphi \partial_t  g_n \Vert_{L^{2}_{t} H^{\sigma}_{x}}<+\infty}$.
\end{enumerate}
Then, the sequence $(g_n)_{n\in \mathbb{N}}$ strongly converges (via a sub-sequence) to a limit $g$ in the space $L^{2}_{loc}(\R \times \Rt)$. 
\end{Lemme}

We will prove now that the sequence $(\vu_n)_{n\in \mathbb{N}}$ verifies the points $1$ and $2$ in this lemma, and for this, we will extend first these functions to whole real line by setting $\vu_n(t,\cdot)=0$ for all $t<-t_n$. Then, let $\varphi \in  \mathcal{C}^{\infty}_{0}(\R \times \Rt)$ be a test function, and let $\tau>0$ be such that $\varphi(t,\cdot)=0$ for all $\vert t \vert > \tau$.  To verify the point $1$, we set $\gamma=1$ to get 
\begin{equation*}
\begin{split}
\Vert  \varphi \vu_n\Vert^{2}_{L^{2}_{t} H^{1}_{x}}=& \int_{-\tau}^{\tau} \Vert \varphi(t,\cdot)\vu_n(t,\cdot)\Vert^{2}_{H^1}dt = \int_{-\tau}^{\tau} \Vert \varphi(t,\cdot)\vu_n(t,\cdot)\Vert^{2}_{L^2}dt + \int_{-\tau}^{\tau} \Vert \varphi(t,\cdot)\vu_n(t,\cdot)\Vert^{2}_{\dot{{H}^{1}}}dt \\
\leq & c_\alpha \, 2\tau(\Vert \varphi \Vert_{L^{\infty}}+\Vert \vec{\nabla}\varphi \Vert_{L^{\infty}}) \left( \sup_{-\tau \leq t \leq \tau}  \Vert \vu_n(t,\cdot)\Vert^{2}_{L^2}+ \sup_{-\tau \leq t  \leq \tau} \alpha^2  \Vert \vu_n(t,\cdot)\Vert^{2}_{\dot{H}^{1}}\right) \\
\leq & C(\tau, \varphi, \alpha) \left( \sup_{-\tau \leq t \leq \tau}  \Vert \vu_n(t,\cdot)\Vert^{2}_{L^2}+ \sup_{-\tau \leq t  \leq \tau} \alpha^2  \Vert \vu_n(t,\cdot)\Vert^{2}_{\dot{H}^{1}}\right).
\end{split}    
\end{equation*}
Moreover, by estimate (\ref{control-tiempo-n}) we can write 
\begin{equation}\label{estim-unif-u-H1}
\begin{split}
 C(\tau, \varphi, \alpha) \left( \sup_{-\tau \leq t \leq \tau}  \Vert \vu_n(t,\cdot)\Vert^{2}_{L^2}+ \sup_{-\tau \leq t  \leq \tau} \alpha^2  \Vert \vu_n(t,\cdot)\Vert^{2}_{\dot{H}^{1}}\right)   \leq    C(\tau, \varphi, \alpha)\left( R^2 + \frac{4}{\beta^2} \left( \Vert \fe \Vert^{2}_{L^2}+\alpha^2\Vert \fe \Vert^{2}_{\dot{H}^1}\right)\right).
 \end{split}
\end{equation}
Thus, for all $n\in \mathbb{N}$ we have 
\begin{equation*}
 \Vert  \varphi \vu_n\Vert^{2}_{L^{2}_{t} H^{1}_{x}} \leq    C(\tau, \varphi, \alpha)\left( R^2 + \frac{4}{\beta^2} \left( \Vert \fe \Vert^{2}_{L^2}+\alpha^2\Vert \fe \Vert^{2}_{\dot{H}^1}\right)\right).
\end{equation*}
We verify now the point $2$. We observe first that for all $n\in \mathbb{N}$ we have 
\begin{equation*}
\varphi\, \partial_t \vu_n=  - \varphi\, \P\,  \text{div} \,( (\vu_n \otimes \vu_n)_\alpha) + \nu \varphi\, \Delta \vu_n +  \varphi\,\fe -\beta \varphi\,\vu_n,  
\end{equation*}
where, for $\tau>0$ given above, we will show that each term  in the right side is uniformly bounded in the space $L^{2}([-\tau,\tau], L^2(\Rt))$. Indeed,  for the term $\ds{\varphi\, \P\,  \text{div} \,( (\vu_n \otimes \vu_n)_\alpha)}$, using the product laws in the Sobolev spaces, and moreover, by estimate (\ref{estim1-n}) we have 
\begin{equation*}
\begin{split}
\Vert \varphi\, \P\,  \text{div} \,( (\vu_n \otimes \vu_n)_\alpha) \Vert_{L^{2}_{t}L^{2}_{x}} & \leq   \Vert \varphi \Vert_{L^{\infty}_{t,x}} \Vert  \P\,  \text{div} \,( (\vu_n \otimes \vu_n)_\alpha) \Vert_{L^{2}_{t}L^{2}_{x}} \leq C(\alpha, \varphi) \Vert \vu_n \otimes \vu_n \Vert_{L^{2}_{t}H^{-1}_{x}} \\
& \leq  C(\alpha, \varphi) \Vert \vu_n \otimes \vu_n \Vert_{L^{2}_{t}H^{1/2}_{x}}
\leq  C(\tau,\alpha, \varphi) \Vert \vu_n \otimes \vu_n \Vert_{L^{\infty}_{t}H^{1/2}_{x}}  \\
&\leq C(\tau,\alpha, \varphi) \Vert \vu_n \Vert^{2}_{L^{\infty}_{t}H^{1}_{x}}\\
& \leq  C(\tau,\alpha, \varphi) \left(R^2 + \frac{4}{\beta^2} \left( \Vert \fe \Vert^{2}_{L^2}+\alpha^2\Vert \fe \Vert^{2}_{\dot{H}^1}\right) \right).
\end{split}
\end{equation*}
For the term $\nu \varphi\, \Delta \vu_n$, using the estimate (\ref{estim2-n}) we can write 
\begin{equation*}
\begin{split}
\nu \Vert  \varphi\, \Delta \vu_n \Vert_{L^{2}_{t}L^{2}_{x}} & \leq \frac{\nu}{\alpha} \Vert \varphi  \Vert_{L^{\infty}_{t,x}}  \alpha \Vert   \Delta \vu_n \Vert_{L^{2}_{t}L^{2}_{x}} \\
&\leq C(\tau, \alpha, \nu, \varphi)\left(\frac{2}{\beta}\left( \Vert \fe \Vert^{2}_{L^2}+\alpha^2 \Vert \fe \Vert^{2}_{\dot{H}^1}\right) + R^2 + \frac{4}{\beta^2} \left( \Vert \fe \Vert^{2}_{L^2}+\alpha^2\Vert \fe \Vert^{2}_{\dot{H}^1} \right)\right)^{1/2}.
\end{split}    
\end{equation*}
Finally, we may observe  that the third and fourth terms: $\varphi\,\fe$ and $\beta \varphi\,\vu_n$, are uniformly bounded in the space  $L^{2}([-\tau,\tau], L^2(\Rt))$, and then, the point $2$ hols true for any $\sigma <0$. \\

We can use now the Lemma \ref{Lemma:RL} to conclude that the sequence $(\vu_n)_{n\in \mathbb{N}}$ converges to the limit $\vv$ in the strong topology of the space $L^{2}_{loc}(\R \times \Rt)$. Then, we have that $\vu_n \otimes \vu_n$ strongly converges to $\vv \otimes \vv $ in the space  $L^{1}_{loc}(\R \times \Rt)$, hence we get that $div(\vu_n \otimes \vu_n)$ converges to $div(\vv \otimes \vv )$ in the distributional sense. \\

At this point, we have proven that $(\vv, q)$ verifies (in the sense of distributions) the equation 
\begin{equation*}
\begin{split}
&   \partial_t (-\alpha^2 \Delta + Id ) \vv +   \text{div} \,(\vv \otimes \vv) -   \nu \Delta (-\alpha^2 \Delta + Id )\vv +  \,  \vec{\nabla}(-\alpha^2 \Delta + Id ) q \\
= &  (-\alpha^2 \Delta + Id )\, \fe -\beta (-\alpha^2 \Delta + Id )\, \vv,
\end{split}   \qquad div(\vv)=0,
\end{equation*} 
but, recalling that $\vv \in \ds{L^{\infty}_{loc}(\R, H^1(\Rt)) \cap L^{2}_{loc}(\R, \dot{H}^2(\Rt))}$ and $q \in L^{2}_{loc}(\R, \dot{H}^3(\Rt))$, we can apply the filtering operator $(\cdot)_\alpha=(-\alpha \Delta + Id)^{-1}$ to each term in this equation to finally obtain that $(\vv,q)$ verifies the equation (\ref{Bardina-enternal}).  Proposition \ref{Prop:eternal-sol} is proven.  \finpv

Once we have constructed a solution $(\vv,q)$ to the equation (\ref{Bardina-enternal}), we will prove now that the sequence $(\vu_n(0,\cdot))_{n\in \mathbb{N}}$ strongly converges (through a sub-sequence) to $\vv(0,\cdot)$ in  $\mathcal{H}^{1}_{\alpha}(\Rt)$. Recall first that this space is equipped with the norm  $\ds{\Vert \cdot \Vert^{2}_{\mathcal{H}^{1}_{\alpha}}=\Vert \cdot \Vert^{2}_{L^2}+\alpha^2 \Vert \cdot \Vert^{2}_{\dot{H}^1}}$. Then, recall moreover that,  for all $n \in \mathbb{N}$ and for all $t \geq -t_n$, the solution $\ds{\vu_n \in L^{\infty}([-t_n,+\infty[, H^1(\Rt)) \cap L^{2}_{loc}([-t_n, +\infty[, \dot{H}^2(\Rt))}$ of equation (\ref{Bardina-n}) verifies the identity (\ref{Eq-dif-energ}): 

\begin{equation*}
\begin{split}
\frac{1}{2} \frac{d}{dt}\Vert \vu_n(t,\cdot)\Vert^{2}_{\mathcal{H}^{1}_{\alpha}} =& -\nu \Vert \vu_n(t,\cdot) \Vert^{2}_{\dot{H}^1}-\alpha^2 \Vert \vu_n(t,\cdot)\Vert^{2}_{\dot{H}^2} + \langle \fe, \vu_n(t,\cdot)\rangle_{L^2\times L^2} \\
 & +\alpha^2 \langle \vec{\nabla} \otimes \fe , \vec{\nabla} \otimes \vu_n(t,\cdot)\rangle_{L^2\times L^2} -\beta \Vert \vu_n(t,\cdot)\Vert^{2}_{\mathcal{H}^{1}_{\alpha}}.
\end{split}
\end{equation*} 
We multiply each term in this identity by $e^{2\beta t}$, and moreover,  we integrate  in the interval $[-t_n, 0]$ to get:
\begin{equation*}
\begin{split}
\frac{1}{2} \Vert \vu_n(0,\cdot)\Vert^{2}_{\mathcal{H}^{1}_{\alpha}} - & \frac{1}{2} e^{- 2 \beta t_n }\Vert \vu_{0,n} \Vert^{2}_{\mathcal{H}^{1}_{\alpha}}-  \beta \int_{-t_n}^{0} e^{2\beta t} \Vert \vu_n(t,\cdot)\Vert^{2}_{\mathcal{H}^{1}_{\alpha}}  dt = -\nu \int_{-t_n}^{0} e^{2 \beta t} \Vert \vu_n(t,\cdot) \Vert^{2}_{\dot{H}^1}  dt \\
&-\alpha^2 \int_{-t_n}^{0} e^{2 \beta t} \Vert \vu_n(t,\cdot)\Vert^{2}_{\dot{H}^2}  dt  + \int_{-t_n}^{0} e^{2 \beta t} \langle \fe, \vu_n(t,\cdot)\rangle_{L^2\times L^2}  dt  \\
& +\alpha^2 \int_{-t_n}^{0} e^{2 \beta t} \langle \vec{\nabla} \otimes \fe , \vec{\nabla} \otimes \vu_n(t,\cdot)\rangle_{L^2\times L^2} -\beta \int_{-t_n}^{0} e^{2\beta t} \Vert \vu_n(t,\cdot)\Vert^{2}_{\mathcal{H}^{1}_{\alpha}}  dt,
\end{split}
\end{equation*}
hence we obtain 
\begin{equation*}
\begin{split}
 \Vert \vu_n(0,\cdot)\Vert^{2}_{\mathcal{H}^{1}_{\alpha}}=&  e^{- 2 \beta t_n }\Vert \vu_{0,n} \Vert^{2}_{\mathcal{H}^{1}_{\alpha}}   -2\nu \int_{-t_n}^{0} e^{2 \beta t} \Vert \vu_n(t,\cdot) \Vert^{2}_{\dot{H}^1}  dt -2 \alpha^2 \int_{-t_n}^{0} e^{2 \beta t} \Vert \vu_n(t,\cdot)\Vert^{2}_{\dot{H}^2}  dt \\
& +2 \int_{-t_n}^{0} e^{2 \beta t} \langle \fe, \vu_n(t,\cdot)\rangle_{L^2\times L^2}  dt  +2 \alpha^2 \int_{-t_n}^{0} e^{2 \beta t} \langle \vec{\nabla} \otimes \fe , \vec{\nabla} \otimes \vu_n(t,\cdot)\rangle_{L^2\times L^2}.
\end{split}
\end{equation*}
In each term of this identity, we take now the $\ds{\limsup}$ when $n\to +\infty$ to write: 
\begin{equation}\label{Estim-energ-limsup}
\begin{split}
\limsup_{n\to +\infty} \Vert \vu_n(0,\cdot)\Vert^{2}_{\mathcal{H}^{1}_{\alpha}}\leq & \limsup_{n\to +\infty}   e^{- 2 \beta t_n }\Vert \vu_{0,n} \Vert^{2}_{\mathcal{H}^{1}_{\alpha}}   +  \limsup_{n\to +\infty} \left( -2\nu \int_{-t_n}^{0} e^{2 \beta t} \Vert \vu_n(t,\cdot) \Vert^{2}_{\dot{H}^1}  dt \right) \\
&+ \limsup_{n\to +\infty} \left( -2 \alpha^2 \int_{-t_n}^{0} e^{2 \beta t} \Vert \vu_n(t,\cdot)\Vert^{2}_{\dot{H}^2}  dt \right)\\
&+\limsup_{n\to +\infty} \left(2 \int_{-t_n}^{0} e^{2 \beta t} \langle \fe, \vu_n(t,\cdot)\rangle_{L^2\times L^2}  dt\right) \\
& + \limsup_{n\to +\infty} \left( 2 \alpha^2 \int_{-t_n}^{0} e^{2 \beta t} \langle \vec{\nabla} \otimes \fe , \vec{\nabla} \otimes \vu_n(t,\cdot)\rangle_{L^2\times L^2} dt\right),
\end{split}
\end{equation}
where we must study each term in the right side. For the first term, recalling that the sequence $(\vu_{0,n})_{n\in \mathbb{N}}$ is bounded in $\mathcal{H}^{1}_{\alpha}(\Rt)$, we have
\begin{equation}\label{E1}
  \limsup_{n\to +\infty}   e^{- 2 \beta t_n }\Vert \vu_{0,n} \Vert^{2}_{\mathcal{H}^{1}_{\alpha}} =0.    
\end{equation}
For the second  term,  since  by estimate (\ref{control-tiempo2-n}) we have that the sequence $(\vu_n)_{n \in \mathbb{N}}$ converges to $\vv$ in the weak-* topology of the space $L^{2}_{loc}(\R, \dot{H}^1(\Rt))$, then we can write 
\begin{equation*}
\liminf_{n\to + \infty} \left( 2\nu \int_{-t_n}^{0} e^{2 \beta t} \Vert \vu_n(t,\cdot) \Vert^{2}_{\dot{H}^1}  dt \right)  \geq 2 \nu \int_{-\infty}^{0} e^{2\beta t} \Vert \vv(t,\cdot)\Vert^{2}_{\dot{H}^1} dt,  
\end{equation*} hence we have
\begin{equation}\label{E2}
\begin{split}
\limsup_{n\to +\infty} \left( -2\nu \int_{-t_n}^{0} e^{2 \beta t} \Vert \vu_n(t,\cdot) \Vert^{2}_{\dot{H}^1}  dt \right) \leq  & - \liminf_{n\to +\infty} \left( 2\nu \int_{-t_n}^{0} e^{2 \beta t} \Vert \vu_n(t,\cdot) \Vert^{2}_{\dot{H}^1}  dt \right)   \\ 
\leq & -2\nu  \int_{-\infty}^{0} e^{2\beta t} \Vert \vv(t,\cdot)\Vert^{2}_{\dot{H}^1} dt. 
\end{split}
\end{equation}
For the third term, always by estimate (\ref{control-tiempo2-n}), we have that the sequence $(\vu_n)_{n \in \mathbb{N}}$ converges to $\vv$ in the weak-* topology of the space $L^{2}_{loc}(\R, \dot{H}^2(\Rt))$, and we write \begin{equation*}
\liminf_{n\to + \infty} \left( 2\alpha^2 \int_{-t_n}^{0} e^{2 \beta t} \Vert \vu_n(t,\cdot) \Vert^{2}_{\dot{H}^2}  dt \right)  \geq 2 \alpha^2 \int_{-\infty}^{0} e^{2\beta t} \Vert \vv(t,\cdot)\Vert^{2}_{\dot{H}^1} dt. 
\end{equation*} Then, we have 
\begin{equation}\label{E3}
\begin{split}
\limsup_{n\to +\infty} \left( -2\alpha^2 \int_{-t_n}^{0} e^{2 \beta t} \Vert \vu_n(t,\cdot) \Vert^{2}_{\dot{H}^2}  dt \right) \leq  & - \liminf_{n\to +\infty} \left( 2\alpha^2 \int_{-t_n}^{0} e^{2 \beta t} \Vert \vu_n(t,\cdot) \Vert^{2}_{\dot{H}^2}  dt \right)   \\ 
\leq & -2\alpha^2  \int_{-\infty}^{0} e^{2\beta t} \Vert \vv(t,\cdot)\Vert^{2}_{\dot{H}^2} dt. 
\end{split}
\end{equation}
Similarly,  by the estimate (\ref{control-tiempo-n}), the sequence $(\vu_n)_{n \in \mathbb{N}}$ converges to $\vv$ in the weak-* topology of the space $L^{2}_{loc}(\R,L^2(\Rt))$, and moreover, as we also have the weak-* convergence in the space $L^{2}_{loc}(\R, \dot{H}^1(\Rt))$. Then, for the fourth and fifth terms we have  
\begin{equation}\label{E4}
\begin{split}
 \limsup_{n\to +\infty} \left(2 \int_{-t_n}^{0} e^{2 \beta t} \langle \fe, \vu_n(t,\cdot)\rangle_{L^2\times L^2}  dt\right)= 2 \int_{-\infty}^{0}  e^{2 \beta t} \langle \fe, \vv(t,\cdot)\rangle_{L^2\times L^2} \, dt, \\ \vspace{3mm}
 \limsup_{n\to +\infty} \left( 2 \alpha^2 \int_{-t_n}^{0} e^{2 \beta t} \langle \vec{\nabla} \otimes \fe , \vec{\nabla} \otimes \vu_n(t,\cdot)\rangle_{L^2\times L^2}\right)= 2 \alpha^2 \int_{-\infty}^{0} e^{2 \beta t}\, \langle \vec{\nabla} \otimes \fe , \vec{\nabla} \otimes \vv(t,\cdot)\rangle_{L^2\times L^2} \, dt.
 \end{split}
\end{equation} 
Thus, gathering the estimates  (\ref{E1}), (\ref{E2}), (\ref{E3}) and (\ref{E4}), we get back to (\ref{Estim-energ-limsup}) to write: 
\begin{equation*}
\begin{split}
\limsup_{n\to +\infty} \Vert \vu_n(0,\cdot)\Vert^{2}_{\mathcal{H}^{1}_{\alpha}}\leq & -2\nu  \int_{-\infty}^{0} e^{2\beta t} \Vert \vv(t,\cdot)\Vert^{2}_{\dot{H}^1} dt -2\alpha^2  \int_{-\infty}^{0} e^{2\beta t} \Vert \vv(t,\cdot)\Vert^{2}_{\dot{H}^2} dt \\
&+ 2 \int_{-\infty}^{0}  e^{2 \beta t} \langle \fe, \vv(t,\cdot)\rangle_{L^2\times L^2} \, dt + 2 \alpha^2 \int_{-\infty}^{0} e^{2 \beta t}\, \langle \vec{\nabla} \otimes \fe , \vec{\nabla} \otimes \vv(t,\cdot)\rangle_{L^2\times L^2} \, dt=(b).
\end{split}
\end{equation*} 

We shall study now the term $(b)$ above. Since the solution $(\vv, q)$ of the equation (\ref{Bardina-enternal}) verifies $\vv \in L^{\infty}(\R, H^1(\Rt)) \cap L^{2}_{loc}(\R, \dot{H}^2(\Rt))$ and $q \in L^{2}_{loc}(\R, H^3(\Rt))$ then, following the same computations done in (\ref{Eq-dif-energ}), we have  the following energy equality:  
\begin{equation*}
\begin{split}
  \frac{1}{2} \frac{d}{dt} \Vert \vv(t,\cdot)\Vert^{2}_{\mathcal{H}^{1}_{\alpha}} =& -\nu \Vert \vv(t,\cdot) \Vert^{2}_{\dot{H}^1}-\alpha^2 \Vert \vv(t,\cdot)\Vert^{2}_{\dot{H}^2} + \langle \fe, \vv(t,\cdot)\rangle_{L^2\times L^2}  +\alpha^2 \langle \vec{\nabla} \otimes \fe , \vec{\nabla} \otimes \vv(t,\cdot)\rangle_{L^2\times L^2} -\beta \Vert \vv(t,\cdot)\Vert^{2}_{\mathcal{H}^{1}_{\alpha}}.
\end{split}  
\end{equation*}

We multiply each term by $e^{2 \beta t}$, and integrating in the interval $]-\infty, 0]$ we get: 
\begin{equation*}
\begin{split}
\frac{1}{2} \Vert \vv(0,\cdot)\Vert^{2}_{\mathcal{H}^{1}_{\alpha}} - & \beta \int_{-\infty}^{0} e^{2\beta t} \Vert \vv(t,\cdot)\Vert^{2}_{\mathcal{H}^{1}_{\alpha}}  dt  = -\nu \int_{-\infty}^{0} e^{2 \beta t} \Vert \vv(t,\cdot) \Vert^{2}_{\dot{H}^1}  dt -\alpha^2 \int_{-\infty}^{0} e^{2 \beta t} \Vert \vv(t,\cdot)\Vert^{2}_{\dot{H}^2}  dt \\
&+ \int_{-\infty}^{0} e^{2 \beta t} \langle \fe, \vv(t,\cdot)\rangle_{L^2\times L^2}  dt  +\alpha^2 \int_{-\infty}^{0} e^{2 \beta t} \langle \vec{\nabla} \otimes \fe , \vec{\nabla} \otimes \vv(t,\cdot)\rangle_{L^2\times L^2} -\beta \int_{-\infty}^{0} e^{2\beta t} \Vert \vv(t,\cdot)\Vert^{2}_{\mathcal{H}^{1}_{\alpha}}  dt.
\end{split}
\end{equation*}

From this identity  we have $\ds{(b)= \Vert \vv(0,\cdot)\Vert^{2}_{\mathcal{H}^{1}_{\alpha}}}$, and then, getting back to the previous estimate, we obtain $\ds{\limsup_{n\to +\infty} \Vert \vu_n(0,\cdot)\Vert^{2}_{\mathcal{H}^{1}_{\alpha}}\leq \Vert \vv(0,\cdot)\Vert^{2}_{\mathcal{H}^{1}_{\alpha}}}$. On the other hand, since by the estimate (\ref{control-tiempo-n}) we know that the sequence  $(\vu_n)_{n\in \mathbb{N}}$ converges (via a sub-sequence) to $\vv$ in the weak-* topology of the space $L^{\infty}(\R, H^1(\Rt))$, we also have the inequality $\ds{\Vert \vv(0,\cdot)\Vert^{2}_{\mathcal{H}^{1}_{\alpha}} \leq \liminf_{n\to +\infty} \Vert \vu_{n}(0,\cdot)\Vert^{2}_{\mathcal{H}^{1}_{\alpha}}}$. Then, we obtain the desired strong convergence: $\ds{\lim_{n\to +\infty} \Vert \vu_n(0,\cdot)\Vert^{2}_{\mathcal{H}^{1}_{\alpha}}= \Vert \vv(0,\cdot)\Vert^{2}_{\mathcal{H}^{1}_{\alpha}}}$,  and the point $2$ in Proposition \ref{Prop:Existence-Attractor} is now verified. \\ 

To  verify now the point $3$ in Proposition \ref{Prop:Existence-Attractor}, we just observe that by estimate (\ref{Estim-Unicidad}), where we have $ \vw(t,\cdot)=\vu_1(t,\cdot)-\vu_2(t,\cdot)= S(t)\vu_{0,1}-S(t)\vu_{0,2}$, and $\vw(0,\cdot)=\vu_{0,1}-\vu_{0,2}$, then the continuity of the map $S(t) : \mathcal{B} \to \mathcal{H}^{1}_{\alpha}(\Rt)$ follows directly. \\

Thus, by Proposition \ref{Prop:Existence-Attractor}, the semigroup $S(t)$ has a  global attractor $\mathcal{A}_{\fe} \subset \mathcal{H}^{1}_{\alpha}(\Rt)$. Theorem \ref{Th-Atractor} is now proven. \finpv 
%%%%%%%%%
\section{Fractal box counting  dimension of the attractor}\label{Sec:FractalDim} 

In this section, we prove that the global attractor $\mathcal{A}_{\fe}\subset \mathcal{H}^{1}_{\alpha}(\Rt)$, constructed in Theorem \ref{Th-Atractor}, 
has finite fractal  box counting dimension and we give an explicit upper bound. In order to estimate the fractal dimension of the attractor, we will use the following volume contraction method adapted from \cite{Ilyn}. See also \cite{Constantin1} and  \cite{Temam} for more details. We start by introducing some definitions  that we shall use later.\\

The first definition concerns the following quasi-differential operator. Let $t\geq 0$ be a fixed time and let $\vu_0 \in \mathcal{A}_{\fe}$ be an initial datum. Moreover, let $u(t,\cdot)$ be the solution  of the equation (\ref{Bardina}) arising from the initial datum $\vu_0$ and given by Theorem \ref{Th-WP}. Thus, for  $u(t,\cdot)$ fixed, let $\vv \in L^{\infty}([0,+\infty[,\mathcal{H}^{1}_{\alpha}(\Rt)) \cap L^{2}_{loc}([0, +\infty[, \dot{H}^2(\Rt))$ be the solution of the following linearized version of the equation (\ref{Bardina}): 

\begin{equation}\label{Bardina-Lin}
\left\{ \begin{array}{ll}\vspace{2mm}
\partial_t \vv + \P \left(\left((\vv \cdot \vec{\nabla}) \vu + (\vu \cdot \vec{\nabla})\vv \right)_\alpha\right) - \nu \Delta \vv  = -\beta \vv, \quad \text{div}(\vv)=0,  \\ 
\vv(0,\cdot)= \vv_0 \in \mathcal{H}^{1}_{\alpha}(\Rt),
\end{array} 
 \right.
\end{equation}
where $\vv_0$ denotes the initial datum. The existence and uniqueness of the solution $\vv$ for this equation is straightforward (since we have a linear equation) and follows  the main ideas and  estimates in the proof of Theorem \ref{Th-WP}. So, we will omit the prove of this fact.

\begin{Definition}[Quasi-differential operator]\label{def-quasi-diff-op} 

The quasi-differential operator $DS(t,\vu_0)$, depending on the time $t\geq 0$ and the datum $\vu_0 \in \mathcal{A}_{\fe}$, is the linear and bounded operator $DS(t,\vu_0): \mathcal{H}^{1}_{\alpha}(\Rt) \to \mathcal{H}^{1}_{\alpha}(\Rt)$ defined as 
\begin{equation*}
DS(t,\vu_0) \vv_0= \vv(t,\cdot),    
\end{equation*}
where  $\vv(t,\cdot)$ is the solution of the linearized  equation (\ref{Bardina-Lin}).
\end{Definition}

Once we have defined this operator, our second definition lies with the notion a semigroup uniformly quasi-differentiable. 

\begin{Definition}[Semigroup uniformly quasi-differentiable]\label{def-unif-quasi-diff} For $t\geq 0$ fixed, let $S(t)$ be the semi-group associated to equation (\ref{Bardina}) and defined in (\ref{def-semigrupo}). We say that this semigroup is uniformly quasi-differentiable on the global attractor $\mathcal{A}_{\fe}\subset \mathcal{H}^{1}_{\alpha}(\Rt)$,  if for all $\vu_{0,1}, \vu_{0,2}\in \mathcal{A}_{\fe}$ we have 
\begin{equation*}
\left\Vert S(t)\vu_{0,2}-S(t)\vu_{0,1}-DS(t,\vu_{0,1})(\vu_{0,2}-\vu_{0,1}) \right\Vert_{\mathcal{H}^{1}_{\alpha}}\leq \mathfrak{o} \left( \Vert \vu_{0,2}-\vu_{0,1}\Vert_{\mathcal{H}^{1}_{\alpha}}\right),    
\end{equation*}
where the quasi-differential operator $DS(t,\vu_{0,1})$ is given in Definition \ref{def-quasi-diff-op}, and moreover, the quantity $\mathfrak{o}(\cdot)$ verifies: $\ds{\lim_{h \to 0^{+}} \mathfrak{o}(h)/ h=0}$. 
\end{Definition}

Finally, in our last definition, we  need to introduce the notion the $m-$ global Lyapunov exponents for a  $m \in \mathbb{N}$ given. For this,  we shall need to precise first  some notation.  On one hand, we denote by $\mathcal{O}_m$, the set of all the orthonormal families  $(\vw_i)_{1\leq i \leq m}$ in the space $\mathcal{H}^{1}_{\alpha}(\Rt)$ dotted with its natural scalar  product:
\begin{equation}\label{def-scalar-product}
[\vw_i,\vw_j ]_\alpha= \left(\vw_i, \vw_j \right)_{L^2\times L^{2}}+ \alpha^2 \left( \vec{\nabla} \otimes \vw_i, \vec{\nabla} \otimes  \vw_j\right)_{L^2\times L^2}.     
\end{equation}
On the other hand,  getting back to the linearized equation (\ref{Bardina-Lin}), we can write 
\begin{equation*}
\partial_t \vv=  - \P \left(\left((\vv \cdot \vec{\nabla}) \vu + (\vu \cdot \vec{\nabla})\vv \right)_\alpha\right) + \nu \Delta \vv  -\beta \vv.
\end{equation*} and then, from the right side of this identity, and for all $\vw \in \mathcal{H}^{1}_{\alpha}(\Rt)$,  we define now the linear operator
\begin{equation}\label{defi-op-L}
\mathcal{L}(t,\vu_0)\vw=  - \P \left(\left((\vw \cdot \vec{\nabla}) \vu + (\vu \cdot \vec{\nabla})\vw \right)_\alpha\right) + \nu \Delta \vw  -\beta \vw.   
\end{equation}
Once we have introduced the set $\mathcal{O}_m$ and the linear operator $\mathcal{L}(t,\vu_0)(\cdot)$ above, we have the following definition.
\begin{Definition}[$m-$ global Lyapunov exponents]\label{def-global-lyapunov-expo} Let $m\in \mathbb{N}$ fixed. We define the $m-$ global Lyapunov exponent $\ell (m)$ as the quantity: 
\begin{equation*}
\ell (m)= \limsup_{T\to+\infty} \left(  \sup_{\vu_0 \in \mathcal{A}_{\fe}}\, \,  \sup_{(\vw_i)_{1\leq i \leq m} \in \mathcal{O}_m} \left( \frac{1}{T}\int_{0}^{T}\sum_{i=1}^{m} \left[ \mathcal{L}(t,\vu_0)\vw_i, \vw_i \right]_{\alpha} dt \right)\right).   
\end{equation*}
\end{Definition} 

We have now all the tools to state the following technical result that we shall use to derive an upper bound of the fractal dimension of the attractor $\mathcal{A}_{\fe}$. For a proof of this result see \cite{Chepyzhov}.

\begin{Theoreme}[Upper bound of the fractal dimension]\label{Th-tec-upper-bound} Let $S(t)$ be the semigruop associated to equation (\ref{Bardina}) and defined in (\ref{def-semigrupo}). Moreover, let $\mathcal{A}_{\fe}\subset \mathcal{H}^{1}_{\alpha}(\Rt)$ be the global attractor of the semigroup $S(t)$ given by Theorem \ref{Th-Atractor}. Finally, let $dim \left(\mathcal{A}_{\fe}\right)$ be the fractal box counting dimension of the attractor $\mathcal{A}_{\fe}$ given in Definition \ref{def-fractal-dim-attractor}. \\

If the following statements hold: 
\begin{enumerate}
    \item The semigroup $S(t)$   is uniformly quasi-differentiable on the attractor $\mathcal{A}_{\fe}$ in the sense of Definition \ref{def-unif-quasi-diff}.
    \item The quasi-differential operator $DS(t,\vu_0)(\cdot)$ given in Definition \ref{def-quasi-diff-op}, depends continuously  on the initial datum $\vu_0 \in \mathcal{A}_{\fe}$. 
    \item There exists $\gamma\geq 1$, and moreover, there exist two constants $c_1,c_2>0$ such that, for all $m\in \mathbb{N}$, the $m-$ global Lyapunov exponent $\ell(m)$ given in Definition \ref{def-global-lyapunov-expo} verifies:
    \begin{equation}\label{estim-lyapunov-exp}
     \ell(m)\leq -c_1 \, m^{\gamma}+c_2.   
    \end{equation}
\end{enumerate}
Then, we have the following upper bound:  $\ds{dim\left(\mathcal{A}_{\fe}\right) \leq \left(\frac{c_2}{c_1}\right)^{1/\gamma}}$.
\end{Theoreme}

\subsection*{Proof of Theorem  \ref{Th-dim-attractor}}
We must verify that the points $1,2$ and $3$ in Theorem \ref{Th-tec-upper-bound} hold. However, the points $1$ and $2$ are given in \cite{Babin} where it is proven that the semigroup $S(t)$ is even differentiable for all $\vu_0 \in \mathcal{A}_{\fe}$ and the differential operator $DS(t,\vu_0)$ depends continuously on $\vu_0 \in \mathcal{A}_{\fe}$. So, we will focus on the point $3$. \\

To estimate the $m-$ global Lyapunov exponent $\ell(m)$ according to the desired estimate (\ref{estim-lyapunov-exp}), we shall prove the following technical estimates. In the expression of the quantity $\ell(m)$ given in Definition \ref{def-global-lyapunov-expo}, we derive first  an upper bound for the term $\ds{\sum_{i=1}^{m} \left[ \mathcal{L}(t,\vu_0) \vw_{i}, \vw_i\right]_\alpha}$ as follows:   

\begin{Proposition}\label{Prop:estim-lyapunov} Let $m\in \mathbb{N}$ fixed and let $(\vw_i)_{1\leq i \leq m} \in \mathcal{O}_m$. Moreover,  let $\mathcal{L}(t,\vu_0)(\cdot)$ be the linear operator given in (\ref{defi-op-L}), and let $[\cdot, \cdot]_\alpha$ be the scalar product defined in (\ref{def-scalar-product}). Then,   we have: 
\begin{equation}\label{estim-tech}
\sum_{i=1}^{m} \left[ \mathcal{L}(t,\vu_0)\vw_i, \vw_i\right]_{\alpha} \leq  -\beta m + 2\frac{C^{4}_{LT}}{\nu^{12/5}\alpha^{6/5}} \Vert \vu(t,\cdot) \Vert^{14/5}_{\dot{H}^{1}}+ \frac{3}{8} \alpha^{2}\Vert \vu(t,\cdot) \Vert^{2}_{\dot{H}^2},
\end{equation}
where $C_{LT}>0$ is a numerical constant given in (\ref{CLT}).
\end{Proposition} 
\pv  By definition of the operator $\mathcal{L}(t,\vu_0)(\cdot)$, we write 
\begin{equation}\label{estil-L}
\sum_{i=1}^{m} \left[ \mathcal{L}(t,\vu_0)\vw_i, \vw_i\right]_{\alpha}= \sum_{i=1}^{m} \left[- \P \left(\left((\vw_i \cdot \vec{\nabla}) \vu + (\vu \cdot \vec{\nabla})\vw_i \right)_\alpha\right), \vw_i\right]_{\alpha}      +  \sum_{i=1}^{m} \left[ \nu \Delta \vw_i  -\beta \vw_i, \vw_i\right]_{\alpha}= I_1+I_2,  
\end{equation}
where we shall study each term $I_1$ and $I_2$ separately. \\

For the term $I_1$, as we have $div(\vw_i)=0$ (for $1\leq i \leq m$), and moreover,  by the well-known properties of the Leray's projector $\P$,  we can write 
\begin{equation*}
I_1=  \sum_{i=1}^{m} \left[-  \left((\vw_i \cdot \vec{\nabla}) \vu + (\vu \cdot \vec{\nabla})\vw_i \right)_\alpha, \vw_i\right]_{\alpha}   
\end{equation*}
Then, we will use the following identity. 
\begin{Lemme} By definition of the filtering operator $(\cdot)_\alpha= (-\alpha^2 \Delta + I_d)^{-1}$, and moreover, by definition of the scalar product $[\cdot , \cdot]_\alpha$ given in (\ref{def-scalar-product}), for $\vg_1, \vg_2 \in H^1(\Rt)$ we have:  $\ds{[ (\vg_1)_\alpha, \vg_2 ]_\alpha = (\vg_1, \vg_2)_{L^2\times L^2}}$.
\end{Lemme}
\pv We  write $\ds{ [ (\vg_1)_\alpha, \vg_2 ]_\alpha =  [ \vg_1, (\vg_2)_\alpha ]_\alpha = \left(\vg_1, (\vg_2)_\alpha   \right)_{L^2\times L^2} +\alpha^2 \left(\vec{\nabla} \otimes \vg_1, \vec{\nabla} \otimes (\vg_2)_\alpha   \right)_{L^2\times L^2}}$, 
and integrating by parts the last term we have 
\begin{equation*}
\begin{split}
 [ (\vg_1)_\alpha, \vg_2 ]_\alpha =&  \left(\vg_1, (\vg_2)_\alpha   \right)_{L^2\times L^2} -\alpha^2 \left(\vg_1, \Delta (\vg_2)_\alpha   \right)_{L^2\times L^2} = \left( \vg_1, (\vg_2)_\alpha - \alpha^2 \Delta((\vg_2)_{\alpha}) \right)_{L^{2}\times L^2} \\
 = &  \left( \vg_1, (-\alpha^2 \Delta + I_d)(\vg_2)_\alpha \right)_{L^{2}\times L^2}=  \left( \vg_1, (-\alpha^2 \Delta + I_d) (-\alpha^2 \Delta + I_d)^{-1}\vg_2 \right)_{L^{2}\times L^2}=  \left( \vg_1,\vg_2\right)_{L^{2}\times L^2}.
 \end{split}
\end{equation*}
\finpv 
Thus, applying  by this identity, and morerover, as we have $div(\vw_i)=0$,  we can write 
\begin{equation*} 
\begin{split}
I_1 = &- \sum_{i=1}^{m} \left((\vw_i \cdot \vec{\nabla}) \vu + (\vu \cdot \vec{\nabla})\vw_i, \vw_i\right)_{L^2 \times L^2}=  -  \sum_{i=1}^{m} \left((\vw_i \cdot \vec{\nabla}) \vu, \vw_i\right)_{L^2 \times L^2} - \sum_{i=1}^{m} \left( (\vu \cdot \vec{\nabla})\vw_i, \vw_i\right)_{L^2 \times L^2}\\
= & \sum_{i=1}^{m} \left((\vw_i \cdot \vec{\nabla}) \vu, \vw_i\right)_{L^2 \times L^2}= \ds{ \sum_{i=1}^{m} \int_{\Rt} \left( \sum_{j,k=1}^{3} w_{i,k} (\partial_k u_j)  \, w_{i,j} \right) dx} \leq \ds{ \sum_{i=1}^{m} \int_{\Rt} \left\vert \sum_{j,k=1}^{3} w_{i,k} (\partial_k u_j)  \, w_{i,j} \right\vert dx}.
\end{split}
\end{equation*}
Here, we need to estimate  the term $\ds{\left\vert \sum_{j,k=1}^{3} w_{i,k} (\partial_k u_j)  \, w_{i,j} \right\vert}$, and following the same computations done in the estimate (3.5), page 16 in \cite{Ilyn}, we have  $\ds{ \left\vert \sum_{j,k=1}^{3} w_{i,k} (\partial_k u_j)  \, w_{i,j} \right\vert \leq \vert \vec{\nabla}\otimes \vu \vert \, \vert \vw_i \vert^2}$.\\

With this pointwise inequality at hand, we get back to the estimate of the term $I_1$ above, where we can write 
\begin{equation*}
I_1 \leq \sum_{i=1}^{m} \int_{\Rt} \vert \vec{\nabla}\otimes \vu \vert \vert \vw_i \vert^2 dx \leq \int_{\Rt} \vert \vec{\nabla}\otimes \vu \vert \left( \sum_{i=1}^{m} \vert \vw_i \vert^2 \right) dx.
\end{equation*} Then, applying the H\"older inequalities (with $2/5+3/5=1$) we have 
\begin{equation*}
\int_{\Rt} \vert \vec{\nabla}\otimes \vu \vert \left( \sum_{i=1}^{m} \vert \vw_i \vert^2 \right) dx \leq \Vert \vec{\nabla}\otimes \vu \Vert_{L^{5/2}} \left\Vert \sum_{i=1}^{m} \vert \vw_i \vert^2 \right\Vert_{L^{5/3}}.    
\end{equation*}
In order to estimate the last term in the right we will use the following Lieb-Thirring inequality, for a  proof see the equation $(6)$, page $2$ in \cite{Lieb}:
\begin{equation*}
 \left\Vert \sum_{i=1}^{m} \vert \vw_i \vert^2 \right\Vert_{L^{5/3}} \leq C_{LT}\, \left( \sum_{i=1}^{m} \Vert \vw_i \Vert^{2}_{\dot{H^1}} \right)^{3/5},  
\end{equation*} where, for the function Gamma $\Gamma(\cdot)$, the constant $C_{LT}>0$ writes down as: 
\begin{equation}\label{CLT}
C_{LT}= \frac{3}{5^{5/3}} \left(16 \pi^{3/2} \frac{\Gamma(7/2)}{\Gamma(5)} \right)^{2/3}.     
\end{equation}
Thus, by the Lieb-Thirring inequality above, and moreover, by the Young inequalities (always with $2/5+3/5=1$), we have 
\begin{equation*}
\begin{split}
\Vert \vec{\nabla}\otimes \vu \Vert_{L^{5/2}} \left\Vert \sum_{i=1}^{m} \vert \vw_i \vert^2 \right\Vert_{L^{5/3}} & \leq C_{LT}\Vert \vec{\nabla}\otimes \vu \Vert_{L^{5/2}} \left( \sum_{i=1}^{m} \Vert \vw_i \Vert^{2}_{\dot{H^1}} \right)^{3/5}  \leq \frac{C_{LT}}{\nu^{3/5}}\Vert \vec{\nabla}\otimes \vu \Vert_{L^{5/2}} \left( \nu \sum_{i=1}^{m} \Vert \vw_i \Vert^{2}_{\dot{H^1}} \right)^{3/5}\\
&\leq \frac{2}{5}\frac{C^{5/2}_{LT}}{\nu^{3/2}} \Vert \vec{\nabla}\otimes \vu \Vert^{5/2}_{L^{5/2}} + \frac{3 \nu}{5} \sum_{i=1}^{m} \Vert \vw_i \Vert^{2}_{\dot{H^1}}.
\end{split} 
\end{equation*}
Now, in the first term to the right hand side, by the interpolation inequalities (with $2/5= \theta/2+ (1-\theta)/6$, and $\theta=7/10$), and moreover, by the Sobolev embedding,  we write 
\begin{equation*}
\begin{split}
 \frac{2}{5}\frac{C^{5/2}_{LT}}{\nu^{3/2}} \Vert \vec{\nabla}\otimes \vu \Vert^{5/2}_{L^{5/2}} + \frac{3 \nu}{5} \sum_{i=1}^{m} \Vert \vw_i \Vert^{2}_{\dot{H^1}} & \leq \frac{2}{5}\frac{C^{5/2}_{LT}}{\nu^{3/2}} \Vert \vec{\nabla}\otimes \vu \Vert^{7/4}_{L^{2}}\,\Vert \vec{\nabla}\otimes \vu \Vert^{3/4}_{L^{6}} + \frac{3 \nu}{5} \sum_{i=1}^{m} \Vert \vw_i \Vert^{2}_{\dot{H^1}}  \\
 & \leq \frac{2}{5}\frac{C^{5/2}_{LT}}{\nu^{3/2}} \Vert  \vu \Vert^{7/4}_{\dot{H}^{1}}\,\left( 4 \Vert \vu \Vert_{\dot{H}^2}\right)^{3/4} + \frac{3 \nu}{5} \sum_{i=1}^{m} \Vert \vw_i \Vert^{2}_{\dot{H^1}} \\
 &\leq \frac{2^{5/2}}{5}\frac{C^{5/2}_{LT}}{\nu^{3/2}} \Vert  \vu \Vert^{7/4}_{\dot{H}^{1}}\, \Vert \vu \Vert^{3/4}_{\dot{H}^2} + \frac{3 \nu}{5} \sum_{i=1}^{m} \Vert \vw_i \Vert^{2}_{\dot{H^1}}\\
 &\leq 2 \frac{C^{5/2}_{LT}}{\nu^{3/2} \alpha^{3/4}} \Vert  \vu \Vert^{7/4}_{\dot{H}^{1}}\, (\alpha \Vert \vu \Vert_{\dot{H}^2})^{3/4} + \frac{3 \nu}{5} \sum_{i=1}^{m} \Vert \vw_i \Vert^{2}_{\dot{H^1}}
\end{split}    
\end{equation*}
At this point, always in the first term to the right hand side, we apply the Young inequalities (with $1=5/8+3/8$) to get
\begin{equation*}
\begin{split}
2 \frac{C^{5/2}_{LT}}{\nu^{3/2} \alpha^{3/4}} \Vert  \vu \Vert^{7/4}_{\dot{H}^{1}}\, (\alpha \Vert \vu \Vert_{\dot{H}^2})^{3/4} + \frac{3 \nu}{5} \sum_{i=1}^{m} \Vert \vw_i \Vert^{2}_{\dot{H^1}}& \leq \frac{5}{8}\left(2 \frac{C^{5/2}_{LT}}{\nu^{3/2} \alpha^{3/4}} \Vert  \vu \Vert^{7/4}_{\dot{H}^{1}}\right)^{8/5} + \frac{3}{8} \alpha^{2}\Vert \vu \Vert^{2}_{\dot{H}^2}+ \frac{3 \nu}{5} \sum_{i=1}^{m} \Vert \vw_i \Vert^{2}_{\dot{H^1}} \\
&\leq 2 \frac{C^{4}_{LT}}{\nu^{12/5}\alpha^{6/5}} \Vert \vu \Vert^{14/5}_{\dot{H}^{1}}+ \frac{3}{8} \alpha^{2}\Vert \vu \Vert^{2}_{\dot{H}^2}+  \frac{3 \nu}{5} \sum_{i=1}^{m} \Vert \vw_i \Vert^{2}_{\dot{H^1}}.
\end{split}    
\end{equation*}
Finally, for the term $I_1$ in (\ref{estil-L}) we have the estimate
\begin{equation}\label{estim-I1}
I_1 \leq  2\frac{C^{4}_{LT}}{\nu^{12/5}\alpha^{6/5}} \Vert \vu(t,\cdot) \Vert^{14/5}_{\dot{H}^{1}}+ \frac{3}{8} \alpha^{2}\Vert \vu(t,\cdot) \Vert^{2}_{\dot{H}^2}+  \frac{3 \nu}{5} \sum_{i=1}^{m} \Vert \vw_i \Vert^{2}_{\dot{H^1}}.   
\end{equation}
We estimate now the term $I_2$ in (\ref{estil-L}). By definition of the scalar product $[\cdot , \cdot]_\alpha$ given in (\ref{def-scalar-product}), and integrating by parts,  we write
\begin{equation*}
\begin{split}
I_2 =& \sum_{i=1}^{m}\int_{\Rt} (\nu \Delta \vw_i - \beta \vw_i)\cdot \vw_i \, dx + \alpha^2 \sum_{i=1}^{m} \int_{\Rt} \vec{\nabla}\otimes (\nu \Delta \vw_i -\beta \vw_i) \cdot \vec{\nabla} \otimes \vw_i \, dx \\
=& \nu \sum_{i=1}^{m}\Vert \vw_i \Vert^{2}_{\dot{H}^1} - \beta \sum_{i=1}^{m} \Vert \vw_i \Vert^{2}_{L^2}-\alpha^2 \nu \sum_{i=1}^{m} \Vert \vw_i \Vert^{2}_{\dot{H}^{2}} - \alpha^2 \beta \sum_{i=1}^{m}\Vert \vw_i \Vert^{2}_{\dot{H}^1}\\
=& -\beta \sum_{i=1}^{m} \left( \Vert \vw_i\Vert^{2}_{L^2}+\alpha^2\Vert \vw_i \Vert^{2}_{\dot{H}^1}\right) -\nu \sum_{i=1}^{m} \Vert \vw_i \Vert^{2}_{\dot{H}^1}-\alpha^2 \nu \sum_{i=1}^{m}\Vert \vw_i \Vert^{2}_{\dot{H}^2}\\
=& -\beta \sum_{i=1}^{m} \Vert \vw_i \Vert^{2}_{\mathcal{H}^{1}_{\alpha}}  -\nu \sum_{i=1}^{m} \Vert \vw_i \Vert^{2}_{\dot{H}^1}-\alpha^2 \nu \sum_{i=1}^{m}\Vert \vw_i \Vert^{2}_{\dot{H}^2},
\end{split}
\end{equation*}
and recalling that $(\vw_i)_{1\leq i \leq m}$ is an orthonormal family in $\mathcal{H}^{1}_{\alpha}(\Rt)$, we finally get
\begin{equation}\label{estim-I2}
I_2= -\beta m   -\nu \sum_{i=1}^{m} \Vert \vw_i \Vert^{2}_{\dot{H}^1}-\alpha^2 \nu \sum_{i=1}^{m}\Vert \vw_i \Vert^{2}_{\dot{H}^2} \leq   -\beta m   -\nu \sum_{i=1}^{m} \Vert \vw_i \Vert^{2}_{\dot{H}^1}. 
\end{equation}
Once we have the estimates (\ref{estim-I1}) and (\ref{estim-I2}), we get back to (\ref{estil-L}) to write
\begin{equation*}
\begin{split}
\sum_{i=1}^{m} \left[ \mathcal{L}(t,\vu_0)\vw_i, \vw_i\right]_{\alpha} & \leq   2\frac{C^{4}_{LT}}{\nu^{12/5}\alpha^{6/5}} \Vert \vu(t,\cdot) \Vert^{14/5}_{\dot{H}^{1}}+ \frac{3}{8} \alpha^{2}\Vert \vu(t,\cdot) \Vert^{2}_{\dot{H}^2}+  \frac{3 \nu}{5} \sum_{i=1}^{m} \Vert \vw_i \Vert^{2}_{\dot{H^1}}    -\beta m   -\nu \sum_{i=1}^{m} \Vert \vw_i \Vert^{2}_{\dot{H}^1} \\
&\leq -\beta m + 2\frac{C^{4}_{LT}}{\nu^{12/5}\alpha^{6/5}} \Vert \vu(t,\cdot) \Vert^{14/5}_{\dot{H}^{1}}+ \frac{3}{8} \alpha^{2}\Vert \vu(t,\cdot) \Vert^{2}_{\dot{H}^2},
\end{split}
\end{equation*}
and then we obtain  the desired estimate. Proposition \ref{Prop:estim-lyapunov} is proven. \finpv 

Once we have the estimate on the quantity $\ds{\sum_{i=1}^{m} \left[ \mathcal{L}(t,\vu_0)\vw_i, \vw_i\right]_{\alpha}}$ given in Proposition \ref{Prop:estim-lyapunov}, we shall continue estimating the quantity $\ell(m)$ given in (\ref{def-global-lyapunov-expo}). So, for $T>0$, we take the time-average $\frac{1}{T} \int_{0}^{T}(\cdot)dt$ in each term of (\ref{estim-tech}) to get 
\begin{equation}\label{estim-tech2}
\frac{1}{T} \int_{0}^{T} \sum_{i=1}^{m} \left[ \mathcal{L}(t,\vu_0)\vw_i, \vw_i\right]_{\alpha} \leq   -\beta m + 2\frac{C^{4}_{LT}}{\nu^{12/5}\alpha^{6/5}}\, \underbrace{\frac{1}{T} \int_{0}^{T} \Vert \vu(t,\cdot) \Vert^{14/5}_{\dot{H}^{1}} dt}_{(a)} + \frac{3}{8}  \underbrace{\alpha^{2} \frac{1}{T} \int_{0}^{T}\Vert \vu(t,\cdot) \Vert^{2}_{\dot{H}^2}dt}_{(b)},  
\end{equation} where we must study now the terms $(a)$ and $(b)$. For the term $(a)$, recalling that $\Vert \cdot \Vert^{2}_{\mathcal{H}^{1}_{\alpha}}=\Vert \cdot \Vert^{2}_{L^2}+\alpha^2 \Vert \cdot \Vert^{2}_{\dot{H}^1}$,  by estimate (\ref{control-tiempo}) we have 
\begin{equation*}
\begin{split}
& \frac{1}{T} \int_{0}^{T} \Vert \vu(t,\cdot) \Vert^{14/5}_{\dot{H}^{1}} dt \leq \frac{1}{\alpha^{14/5}}  \frac{1}{T} \int_{0}^{T} \left(\alpha^2 \Vert \vu(t,\cdot) \Vert^{2}_{\dot{H}^{1}}\right)^{7/5} dt  \leq \frac{1}{\alpha^{14/5}}  \frac{1}{T} \int_{0}^{T} \left( \Vert \vu_0 \Vert^{2}_{\mathcal{H}^{1}_{\alpha}}\, e^{-\beta\, t} + \frac{4}{\beta^2} \Vert \fe \Vert^{2}_{\mathcal{H}^{1}_{\alpha}}\right)^{7/5} dt \\
 &\leq  \frac{2^{2/5}}{\alpha^{14/5}}\left( \Vert \vu_0 \Vert^{14/5}_{\mathcal{H}^{1}_{\alpha}}  \frac{1}{T} \int_{0}^{T} e^{-\frac{7 \beta}{5}t} dt + \frac{2^{14 /5}}{\beta^{14/5}} \Vert \fe \Vert^{14/5}_{\mathcal{H}^{1}_{\alpha}}\right) \leq  \frac{2^{2/5}}{\alpha^{14/5}} \Vert \vu_0 \Vert^{14/5}_{\mathcal{H}^{1}_{\alpha}} \frac{5}{7 \beta \, T} (1-e^{-\frac{7 \beta}{5}T})+ \frac{2^{16/5}}{\alpha^{14/5}\beta^{14/5}} \Vert \fe \Vert^{14/5}_{\mathcal{H}^{1}_{\alpha}}.   
\end{split} 
\end{equation*}
Recall also that by definition of the quantity $\ell(m)$ (see always the expression (\ref{def-global-lyapunov-expo})) we have $\vu_0 \in \mathcal{A}_{\fe}$, and moreover,  as the global attractor $\mathcal{A}_{\fe}$ is a compact set in $\mathcal{H}^{1}_{\alpha}(\Rt)$ (see Definition \ref{Def:atractor}) then there exists $M>0$ such that $\Vert \vu_{0} \Vert_{\mathcal{H}^{1}_{\alpha}} \leq M$ for all $\vu_0 \in \mathcal{H}^{1}_{\alpha}(\Rt)$. Then, by the previous estimate we can write
\begin{equation}\label{estim-H1}
\frac{1}{T} \int_{0}^{T} \Vert \vu(t,\cdot) \Vert^{14/5}_{\dot{H}^{1}} dt \leq \frac{2^{2/5}}{\alpha^{14/5}}  \frac{5 M^{14/5}}{7 \beta \, T} (1-e^{-\frac{7 \beta}{5}T})+ \frac{2^{16/5}}{\alpha^{14/5}\beta^{14/5}} \Vert \fe \Vert^{14/5}_{\mathcal{H}^{1}_{\alpha}}.
\end{equation}
In order to estimate the term $(b)$, we will use now the inequality (\ref{control-tiempo2}) (with $t=0$ and $T>0$) to write 
\begin{equation}\label{estim-H2}
\alpha^{2} \frac{1}{T} \int_{0}^{T}\Vert \vu(t,\cdot) \Vert^{2}_{\dot{H}^2}dt \leq   \frac{2}{\beta} \Vert \fe \Vert^{2}_{\mathcal{H}^{1}_{\alpha}}+ \frac{1}{T} \Vert \vu_{0} \Vert^{2}_{\mathcal{H}^{1}_{\alpha}}\leq  \frac{2}{\beta} \Vert \fe \Vert^{2}_{\mathcal{H}^{1}_{\alpha}}+ \frac{1}{T}M^2. 
\end{equation}
Gathering the estimates (\ref{estim-H1}) and (\ref{estim-H2}), in (\ref{estim-tech2}) we have 
\begin{equation*}
\begin{split}
 \frac{1}{T} \int_{0}^{T} \sum_{i=1}^{m} \left[ \mathcal{L}(t,\vu_0)\vw_i, \vw_i\right]_{\alpha} dt \leq &   -\beta m + 2\frac{C^{4}_{LT}}{\nu^{12/5}\alpha^{6/5}} \left( \frac{2^{2/5}}{\alpha^{14/5}}  \frac{5 M^{14/5}}{7 \beta \, T} (1-e^{-\frac{7 \beta}{5}T})+ \frac{2^{16/5}}{\alpha^{14/5}\beta^{14/5}} \Vert \fe \Vert^{14/5}_{\mathcal{H}^{1}_{\alpha}} \right)   \\
  &+ \frac{3}{8} \left(   \frac{2}{\beta} \Vert \fe \Vert^{2}_{\mathcal{H}^{1}_{\alpha}}+ \frac{1}{T}M^2 \right)\\
  \leq &  -\beta m + \frac{1}{T} \underbrace{\left( 2\frac{C^{4}_{LT}}{\nu^{12/5}\alpha^{6/5}} \frac{2^{2/5}}{\alpha^{14/5}}  \frac{5 M^{14/5}}{7 \beta} (1-e^{-\frac{7 \beta}{5}T})  + \frac{3}{8}M^2 \right)}_{(c)}\\
  &+ 2\frac{C^{4}_{LT}}{\nu^{12/5}\alpha^{6/5}} \frac{2^{16/5}}{\alpha^{14/5}\beta^{14/5}} \Vert \fe \Vert^{14/5}_{\mathcal{H}^{1}_{\alpha}} + \frac{3}{4 \beta} \Vert \fe \Vert^{2}_{\mathcal{H}^{1}_{\alpha}},
\end{split}    
\end{equation*}
and moreover, we set now the constant 
\begin{equation}\label{defi-const}
C(\alpha,\beta,\nu)=  2\frac{C^{4}_{LT}}{\nu^{12/5}\alpha^{6/5}} \frac{2^{16/5}}{\alpha^{14/5}}+ \frac{3}{4 \beta}, 
\end{equation} to write
\begin{equation*}
 \frac{1}{T} \int_{0}^{T} \sum_{i=1}^{m} \left[ \mathcal{L}(t,\vu_0)\vw_i, \vw_i\right]_{\alpha} dt \leq   -\beta m +  \frac{1}{T}(c)+ C(\alpha,\beta,\nu)\max\left(\Vert \fe \Vert^{14/5}_{\mathcal{H}^{1}_{\alpha}}, \Vert \fe \Vert^{2}_{\mathcal{H}^{1}_{\alpha}} \right).
\end{equation*}
Finally, by (\ref{def-global-lyapunov-expo}) we obtain the estimate 
\begin{equation}
\ell(m)\leq - \beta m   + c(\alpha,\beta,\nu)\max\left(\Vert \fe \Vert^{14/5}_{\mathcal{H}^{1}_{\alpha}}, \Vert \fe \Vert^{2}_{\mathcal{H}^{1}_{\alpha}} \right), 
\end{equation} which is the desired estimate (\ref{estim-lyapunov-exp}) with $c_1=\beta$, $\gamma=1$ and $c_2= C(\alpha,\beta,\nu)\max\left(\Vert \fe \Vert^{14/5}_{\mathcal{H}^{1}_{\alpha}}, \Vert \fe \Vert^{2}_{\mathcal{H}^{1}_{\alpha}} \right)$. Thus, appliying the Theorem \ref{Th-tec-upper-bound} we get the upper bound for the fractal dimension of the attractor $\mathcal{A}_{\fe}$ given in (\ref{dimfrac}), where the constant $c(\alpha, \beta,\nu)$ is defined through the constant $C(\alpha,\beta, \nu)$, given in (\ref{defi-const}) as 
\begin{equation}\label{defi-constant-Th}
c(\alpha,\beta,\nu)= \frac{1}{\beta} C(\alpha, \beta,\nu).
\end{equation}
Theorem \ref{Th-dim-attractor} is now proven. \finpv

\section{Internal estructure of the global attractor $\mathcal{A}_{\fe}$}\label{Sec:Stationary-Sol} 
\subsection*{Proof of Theorem \ref{Th-Stationary-Solutions}}
In order to prove this theorem, we will verify each point stated in this result separately. 
\subsubsection*{1) The existence of stationary solutions} 
Given an external force $\fe \in \mathcal{H}^{1}_{\alpha}(\Rt)$, we will construct a solution $(\U,P)$ for the stationary problem (\ref{Bardina-Stat}).  We shall start by explaining the strategy of our proof. Since we work on the whole space $\Rt$, and moreover, since we have $div(\fe)=0$ and $div(\U)=0$, then we can apply the Leray's projector $\P$ to the equation (\ref{Bardina-Stat}) to obtain the following equation:
\begin{equation}
-\nu \Delta \U + \P \text{div}((\U \otimes \U)_\alpha)  = \fe -\beta \U.
\end{equation}
From this equation, as $div(\U)=0$ we may observe that the velocity $\U$ formally solves the following \emph{equivalent} fixed point problem: 
\begin{equation}\label{Fixed-Problem}
\begin{split}
\U&=-\frac{1}{-\nu \Delta + \beta I_d} \left[ \P \text{div}((\U \otimes \U)_\alpha)\right] + \frac{1}{-\nu \Delta + \beta I_d} \left[ \fe\right] \\
& = -\frac{1}{-\nu \Delta + \beta I_d} \left[ \P \left( (\U \cdot \vec{\nabla}) \U \right)_\alpha \right] + \frac{1}{-\nu \Delta + \beta I_d} \left[ \fe\right]= T(\U),
\end{split}
\end{equation} where the non-local operator $\ds{\frac{1}{-\nu \Delta + \beta I_d}}$ is easily defined in the Fourier variable by its symbol $\ds{\frac{1}{\nu \vert \xi \vert^2+\beta}}$. Thus, the main idea to construct a solution $(\U,P)$ to the equation (\ref{Bardina-Stat}) is to solve first this fixed point problem, and  then, the pressure term $P$ is calculated through the velocity $\U$ by the formula  (\ref{Caract-pression}). \\

To solve the fixed point problem (\ref{Fixed-Problem}), we could use the (classical) Picard's iterative scheme in the space $H^2(\Rt)$, however, this scheme needs a control on the quantity $\ds{\left\Vert \frac{1}{-\nu \Delta + \beta I_d} \left[ \fe\right] \right\Vert_{H^2}}$. Moreover,  as we have the estimate $\ds{\left\Vert \frac{1}{-\nu \Delta + \beta I_d} \left[ \fe\right] \right\Vert_{H^2} \lesssim \Vert \fe \Vert_{H^1}}$, we observe that the  Picard's iterative scheme  ultimately needs a control on the quantity $\ds{\Vert \fe \Vert_{H^1}}$. Consequently, this scheme  works for the case of  external forces small enough: $\ds{\Vert \fe \Vert_{H^1} \lesssim 1}$.\\

Since we want to construct stationary solutions $(\U,P)$ for any external force $\fe \in \mathcal{H}^{1}_{\alpha}(\Rt)$, we shall use here a different approach. Instead of the Picard's fixed point theorem, the key tool to construct a solution $(\U,P)$ will be the Sheafer's fixed point theorem. For a proof see the Theorem $16.1$, page $529$ of \cite{PGLR1}. 
\begin{Theoreme}[Sheafer's fixed point]\label{Th-Sheafer} Let $E$ be a Banach space and let $T:E \to E$ such that: 
\begin{enumerate}
\item[1)] $T$ is a continuous operator.
\item[2)] $T$ is a compact operator.
\item[3)] There exists a constant $M>0$, such that for all $\lambda \in [0,1]$, if $e\in E$ verifies $e=\lambda T(e)$ then we have $\Vert e \Vert_E \leq M$.
\end{enumerate}
Then, there exists $\overline{e} \in E$ such that $\overline{e}   =T(\overline{e})$. 
\end{Theoreme}
It is worth mention  this result asserts the existence of at least a fixed point of the operator $T$, but there is not any supplementary information on its uniqueness. \\

We would like to apply this result to the Banach space $E=\left\{ \U \in H^2: \,\, div(\U)=0\right\}$, and the operator $T$ given in (\ref{Fixed-Problem}), however, there is a technical problem to overcome: to the best of our knowledge, we do not how to verify that  the compactness of  the operator  $T$ in the space $E$ above. To solve this problem, we will consider  a family of \emph{modified} operators $(T_r)_{r>0}$ which verify all the points of the  Sheafer's fixed theorem above.  Then, for all $r>0$ fixed, by this theorem we will get a solution $\U_r \in E$ of the equation $\U_r =T_r(\U_r)$. Finally, using a Rellich-Lions lemma we will show that the sequence $(\U_r)_{r>0}$ converges to a solution of the stationary problem (\ref{Bardina-Stat}). \\  

For $r>0$ fixed, we define now the modified operator $T_r:E \to E$ as follows. Let $\theta \in \mathcal{C}^{\infty}_{0}(\Rt)$ be a test function such that $0\leq \theta(x) \leq 1$, $\theta(x)=1$ for $\vert x \vert <1$, and $\theta(x)=0$ for $\vert x \vert >2$. Then,  we define the cut-off function $\theta_r(x)= \theta(x / r)$, and moreover, we define  $T_r$ by
\begin{equation}\label{def-Op-Tr}
T_r(\U)= -\frac{1}{-\nu \Delta + \beta I_d} \left[ \P \left( (\theta_r \U \cdot \vec{\nabla}) (\theta_r\U) \right)_\alpha \right] + \frac{1}{-\nu \Delta + \beta I_d} \left[ \fe\right].
\end{equation}

We will prove that this operator verifies the points stated in the Theorem \ref{Th-Sheafer}.  For the point $1)$, for $\U, \V \in E$ we write 
\begin{equation}\label{estim-tech-Tr}
\begin{split}
\Vert T_r(\U) - T_r(\V) \Vert_{H^2}=& \left\Vert -\frac{1}{-\nu \Delta + \beta I_d} \left[ \P \left( (\theta_r \U \cdot \vec{\nabla}) (\theta_r\U) -  (\theta_r \V \cdot \vec{\nabla}) (\theta_r\V) \right)_\alpha \right]  \right\Vert_{H^2}  \\
=&  \left\Vert -\frac{1}{-\nu \Delta + \beta I_d} \left[ \P \left(  ((\theta_r(\U-\V)) \vec{\nabla})(\theta_r \U) + ((\theta_r \V)\cdot \vec{\nabla})(\theta_r(\U-\V))\right)_\alpha \right]  \right\Vert_{H^2}. 
\end{split}    
\end{equation}
Using the fact that the symbol of the operator $\ds{\frac{1}{-\nu \Delta + \beta I_d}}$ is a  bounded function in the frequency variable, and moreover, by the well-known properties of the Leray's projector $\P$,  we can write
\begin{equation*}
\begin{split}
 & \left\Vert -\frac{1}{-\nu \Delta + \beta I_d} \left[ \P \left(  ((\theta_r(\U-\V)) \vec{\nabla})(\theta_r \U) + ((\theta_r \V)\cdot \vec{\nabla})(\theta_r(\U-\V))\right)_\alpha \right]  \right\Vert_{H^2} \\
 \leq & c(\beta,\nu)\left( \left\Vert \left(  ((\theta_r(\U-\V)) \vec{\nabla})(\theta_r \U) + ((\theta_r \V)\cdot \vec{\nabla})(\theta_r(\U-\V))\right)_\alpha \right\Vert_{H^2}\right),
\end{split}    
\end{equation*}
and then, by definition of the filterig operator $(\cdot)_{\alpha}$ given in formula (\ref{filtering}) we can write 
\begin{equation*}
\begin{split}
\Vert T_r(\U) - T_r(\V) \Vert_{H^2}  \leq & c(\beta,\nu)\left( \left\Vert \left(  ((\theta_r(\U-\V)) \vec{\nabla})(\theta_r \U) + ((\theta_r \V)\cdot \vec{\nabla})(\theta_r(\U-\V))\right)_\alpha \right\Vert_{H^2}\right) \\
\leq & c(\alpha,\beta,\nu) \left\Vert  ((\theta_r(\U-\V)) \vec{\nabla})(\theta_r \U) + ((\theta_r \V)\cdot \vec{\nabla})(\theta_r(\U-\V)) \right\Vert_{L^2}.  
\end{split}    
\end{equation*}
Moreover, applying the H\"older inequalities and recalling that the space $H^2(\Rt)$ embeds in the spaces $L^{\infty}(\Rt)$ and $H^{1}(\Rt)$,  we obtain 
\begin{equation*}
\begin{split}
&  c(\alpha,\beta,\nu) \left\Vert  ((\theta_r(\U-\V)) \vec{\nabla})(\theta_r \U) + ((\theta_r \V)\cdot \vec{\nabla})(\theta_r(\U-\V)) \right\Vert_{L^2} \\
\leq & c(\alpha,\beta,\nu) \left( \left\Vert \theta_r(\U-\V) \right\Vert_{L^{\infty}} \left\Vert \vec{\nabla}\otimes (\theta_r \U)\right\Vert_{L^2} + \left\Vert \theta_r \V\right\Vert_{L^{\infty}} \left\Vert \vec{\nabla}\otimes (\theta_r(\U-\V))\right\Vert_{L^2}  \right)\\
\leq & c(\alpha,\beta,\nu,\theta_r) \left( \left\Vert \U-\V \right\Vert_{L^{\infty}} \left\Vert  \U\right\Vert_{H^1} + \left\Vert  \V\right\Vert_{L^{\infty}} \left\Vert \U-\V\right\Vert_{H^1}  \right) \\
\leq & c(\alpha,\beta,\nu,\theta_r)\left( \left\Vert \U-\V \right\Vert_{H^2}\left\Vert  \U\right\Vert_{H^2} + \left\Vert  \V\right\Vert_{H^2} \left\Vert \U-\V\right\Vert_{H^2}\right). 
\end{split}    
\end{equation*}
Thus, we are able to write 
\begin{equation*}
\left\Vert T_r(\U) - T_r(\V) \right\Vert_{H^2} \leq c(\alpha,\beta,\nu,\theta_r) \left( \left\Vert  \U\right\Vert_{H^2}+\left\Vert  \V\right\Vert_{H^2}\right) \, \left\Vert \U-\V \right\Vert_{H^2},   
\end{equation*}
hence we have the continuity of the operator $T_r:E \to E$. \\

We verify now the point $2)$.  Let $(\V_n)_{n\in \mathbb{N}}$ be a bounded sequence in the space $E$. Then, this sequence converges to a limit $\V \in H^2(\Rt)$ in the weak  topology of the space $H^2(\Rt)$. Moreover, as we have $div(\V_n)=0$, for all $n\in \mathbb{N}$, then we get that $div(\V)=0$, and thus we have $\V \in E$. \\

We must show that the sequence $(T_r(\V_n))_{n\in \mathbb{N}}$ strongly converges (via a sub-sequence) to $T_r(\V)$ in the space $E$.  For this, by (\ref{estim-tech-Tr}), and moreover, by the well-known properties of operators  $\ds{\frac{1}{-\nu \Delta + \beta I_d}}$, $\P$ and $\ds{(\cdot)_\alpha= \frac{1}{-\alpha^2 \Delta +  I_d}}$, we can write 
\begin{equation}\label{estim-tech-comp}
\begin{split}
 \Vert T_r(\V_n) - T_r(\V) \Vert_{H^2}  = & \left\Vert\frac{1}{-\nu \Delta + \beta I_d} \left[ \P \left(  ((\theta_r(\V_n-\V)) \vec{\nabla})(\theta_r \V_n) + ((\theta_r \V)\cdot \vec{\nabla})(\theta_r(\V_n-\V))\right)_\alpha \right]  \right\Vert_{H^2} \\
 \leq & c_{\alpha,\beta,\nu}\,  \left\Vert   ((\theta_r(\V_n-\V)) \vec{\nabla})(\theta_r \V_n) + ((\theta_r \V)\cdot \vec{\nabla})(\theta_r(\V_n-\V)) \right\Vert_{H^{-2}} \\
  \leq & c_{\alpha,\beta,\nu}\,\left\Vert   ((\theta_r(\V_n-\V)) \vec{\nabla})(\theta_r \V_n) + ((\theta_r \V)\cdot \vec{\nabla})(\theta_r(\V_n-\V)) \right\Vert_{H^{-1}} \\
   \leq & c_{\alpha,\beta,\nu}\, \left\Vert   ((\theta_r(\V_n-\V)) \vec{\nabla})(\theta_r \V_n) + ((\theta_r \V)\cdot \vec{\nabla})(\theta_r(\V_n-\V)) \right\Vert_{\dot{H}^{-1}}\\
 \leq & c_{\alpha,\beta,\nu}\,\left( \underbrace{\left\Vert  ((\theta_r(\V_n-\V)) \vec{\nabla})(\theta_r \V_n)  \right\Vert_{\dot{H}^{-1}}}_{(a)} +  \underbrace{ \left\Vert ((\theta_r \V)\cdot \vec{\nabla})(\theta_r(\V_n-\V)) \right\Vert_{\dot{H}^{-1}}}_{(b)}\right),
\end{split}    
\end{equation}
where, we will prove that the terms $(a)$ and $(b)$ converge to zero when $n$ goes to $+\infty$.  \\

For the term $(a)$, by the Hardy-Littlewood-Sobolev inequalities,  by the H\"older inequalities, and moreover, recalling that the sequence $(\V_n)_{n\in \mathbb{N}}$ is bounded in $H^2(\Rt)$,  we can write 
\begin{equation*}
\begin{split}
(a) & \leq c \left\Vert  ((\theta_r(\V_n-\V)) \vec{\nabla})(\theta_r \V_n)  \right\Vert_{L^{6/5}}  \leq c \left\Vert  (\theta_r(\V_n-\V) \right\Vert_{L^{3}}\, \left\Vert \vec{\nabla} \otimes (\theta_r \V_n) \right\Vert_{L^{2}}  \\
&\leq c(\theta_r) \left\Vert  (\theta_r(\V_n-\V) \right\Vert_{L^{3}}\, \Vert \V_n \Vert_{H^1}  \leq c(\theta_r) \left\Vert  (\theta_r(\V_n-\V) \right\Vert_{L^{3}}\Vert \V_n \Vert_{H^2} \leq  c(\theta_r) \left\Vert  (\theta_r(\V_n-\V) \right\Vert_{L^{3}}. 
\end{split}    
\end{equation*}
We will prove now the term in the right side converges to zero when $n$ goes to $+\infty$. For this we have the following technical result.

\begin{Lemme}\label{Lema-tech-conv} For $r>0$ fixed, the sequence $(\theta_r\V_n)_{n\in \mathbb{N}}$  strongly converges (through a sub-sequence) to $\theta_r \V$ in the space $L^p(\Rt)$ with $2\leq p < +\infty$.
\end{Lemme}
\pv 
 We observe first that the sequence $(\theta_{2r}\V_n)_{n\in \mathbb{N}}$ is also bounded in the space $H^2(\Rt)$, and moreover, by definition of the cut-off function $\theta_{2r}$ we have $\ds{supp\left( \theta_{2r} \V_n \right)\subset B(0, 4r)}$, for all $n\in \mathbb{N}$. Then, by the Rellich-Lions lemma there exists $\W \in H^2(\Rt)$ such that the sequence $(\theta_{2r}\V_n)_{n\in \mathbb{N}}$ strongly converges (through a sub-sequence) to $\W$ in $L^2(\Rt)$. Moreover, the sequence $(\theta_{2r}\V_n)_{n\in \mathbb{N}}$ strongly converges  to $\W$ in $L^p(\Rt)$ with $2\leq p \leq +\infty$. On the other hand,  as we have $\theta_{2r}=1$ on $supp(\theta_r)$ then we obtain that the sequence $(\theta_{r}\V_n)_{n\in \mathbb{N}}$ strongly converges  to $\theta_r \W$ in $L^p(\Rt)$ with $2\leq p < +\infty$.  Finally, we will prove the identity $\theta_r\W= \theta_r \V$. For this, just recall that $\V_n$ converges to $\V$ in the weak topology of $H^2(\Rt)$ and then it also converges to $\V$ in the weak topology of $L^2(\Rt)$. Hence we get that the sequence  $\theta_r \V_n$ weakly converges to $\theta_r \V$ in $L^2(\Rt)$. Thus, since $\theta_r \V_n$ strongly converges $\theta_r \W$ in $L^2(\Rt)$ then we have $\theta_r\W= \theta_r \V$. \finpv 

We study now the term $(b)$. We will need here the following identity. 
\begin{Lemme} Let $\vec{A}=(A_1, A_2, A_3), \vec{B}=(B_1, B_2, B_3)$ be two vector fields. If $div(\vec{A})=0$ then we have  $$\ds{((\theta_r\vec{A})\cdot \vec{\nabla})(\theta_r \vec{B})= div (\theta^{2}_{r} \vec{A} \otimes \vec{B})-\theta_r (\vec{\nabla}\theta_r \cdot \vec{A})\vec{B}}.$$
\end{Lemme}
\pv  We write 
\begin{equation*}
\begin{split}
((\theta_r\vec{A})\cdot \vec{\nabla})(\theta_r \vec{B})= &\sum_{i=1}^{3} \theta_r A_i \partial_i (\theta_r \vec{B})= \sum_{i=1}^{3} \partial_i \left(  \theta^{2}_{r} A_i \vec{B}\right)-  \sum_{i=1}^{3} \theta_r \partial_i (\theta_r A_i) \vec{B} = div (\theta^{2}_{r} \vec{A} \otimes \vec{B})-   \sum_{i=1}^{3} \theta_r \partial_i (\theta_r A_i) \vec{B}\\
= & div (\theta^{2}_{r} \vec{A} \otimes \vec{B}) -  \sum_{i=1}^{3} \theta_r (\partial_i \theta_r) A_i- \sum_{i=1}^{3} \theta^{2}_{r} (\partial_i A_i) \vec{B}  = div (\theta^{2}_{r} \vec{A} \otimes \vec{B})-\theta_r (\vec{\nabla}\theta_r \cdot \vec{A})\vec{B}.\qquad\,\,\,\,\,  \hfill \blacksquare
\end{split}    
\end{equation*}

By this lemma, setting $\vec{A}=\V$ and $\vec{B}=\V_n-\V$, in the term $(b)$ we have  
\begin{equation*}
\begin{split}
(b)=&  \left\Vert div (\theta^{2}_{r} \V \otimes (\V_n-\V)) \right\Vert_{\dot{H}^{-1}} + \left\Vert \theta_r (\vec{\nabla}\theta_r \cdot \V)(\V_n-\V) \right\Vert_{\dot{H}^{-1}} =(b_1)+(b_2),  
\end{split}    
\end{equation*}
and then, we must study the terms $(b_1)$ and $(b_2)$ separately. For the term $(b_1)$ we write 
\begin{equation*}
\begin{split}
(b_1) \leq c  \left\Vert \theta^{2}_{r} \V \otimes (\V_n-\V)) \right\Vert_{L^2} \leq  \left\Vert \theta_r \V \right\Vert_{L^{\infty}}\, \left\Vert \theta_r (\V_n-\V)\right\Vert_{L^2} \leq c \left\Vert \V \right\Vert_{H^2} \, \left\Vert \theta_r (\V_n-\V)\right\Vert_{L^2}   
\end{split}    
\end{equation*}
where, by Lemma \ref{Lema-tech-conv} we conclude that this term tends to zero when $n$ goes to $+\infty$. On the other hand, to estimate the term $(b_2)$, we use again the Hardy-Littlewood-Sobolev inequalities and the H\"older inequalities to write 
\begin{equation*}
\begin{split}
(b_2) \leq & c \left\Vert \theta_r (\vec{\nabla}\theta_r \cdot \V)(\V_n-\V) \right\Vert_{L^{6/5}} \leq  \left\Vert  \vec{\nabla}\theta_r \cdot \V \right\Vert_{L^{3}}\, \left\Vert \theta_r (\V_n-\V)\right\Vert_{L^2},   
\end{split}    
\end{equation*}
where, always by  Lemma \ref{Lema-tech-conv}, we also conclude that this term tends to zero when $n$ goes to $+\infty$. \\

Once we have shown that the $(a)$ and $(b)$ converge to zero, by (\ref{estim-tech-comp}) we obtain that $T_r(\V_n)$ strongly converges to $T_r(\V)$ in the space $E$, and then, the operator $T_r$ is compact.\\

We verify now the point $3)$. Let $\lambda \in [0,1]$, and $\U \in E$ such that it verifies $\U= \lambda T_r(\U)$. Then, by definition of the operator $T_r$ given in (\ref{def-Op-Tr}), we get that $\U$ solves the equation: 
\begin{equation*}
 -\nu \Delta \U +\beta \U + \lambda\, \P \left( (\theta_r \U \cdot \vec{\nabla}) (\theta_r\U) \right)_\alpha =  \lambda\, \fe. 
\end{equation*}
Moreover, applying the operator $(-\alpha^2 \Delta + Id)$ to each term in this equation, and recalling that we have $(\cdot)_\alpha=(-\alpha^2 \Delta + Id)^{-1}$, we obtain that $\U$ also solves the equation: 
\begin{equation*}
 \nu\alpha^2 \Delta^{2} \U -(\nu+\beta \alpha^2) \Delta \U +\beta \U +\lambda\, \P  (\theta_r \U \cdot \vec{\nabla}) (\theta_r\U) = -\alpha^2 \lambda\, \Delta\fe + \lambda\,\fe. 
\end{equation*} 
In this equation, we multiply first by $\U$ and then we integrate on the whole space $\Rt$. After some  integration by parts, and moreover, as we have $div(\U)=0$, we finally get 
\begin{equation*}
 \nu\alpha^2  \Vert \U \Vert^{2}_{\dot{H}^2}+ (\nu+\beta \alpha^2) \Vert \U \Vert^{2}_{\dot{H}^{1}}+ \beta \Vert \U \Vert^{2}_{L^2} = \alpha^2 \lambda\,\int_{\Rt} \vec{\nabla}\otimes \fe \cdot \vec{\nabla} \otimes\U \, dx + \lambda\, \int_{\Rt} \fe \cdot \U\, dx.    
\end{equation*}
In the term in the right side, by the Cauchy-Schwarz inequalities, using the fact that $\lambda \in [0,1]$, and moreover, by the Young inequalities, we write
\begin{equation}\label{estim-term-force}
\begin{split}
-\alpha^2\, \lambda \int_{\Rt} \fe \cdot \Delta \U \, dx + \lambda \int_{\Rt} \fe \cdot \U\, dx \leq &\alpha^2 \lambda\, \Vert \fe \Vert_{\dot{H^1}}\, \Vert \U \Vert_{\dot{H}^1} + \lambda\, \Vert \fe \Vert_{L^2}\, \Vert \U \Vert_{L^2} \\
\leq & \alpha^2  \Vert \fe \Vert_{\dot{H}^1}\, \Vert \U \Vert_{\dot{H}^1} +  \Vert \fe \Vert_{L^2}\, \Vert \U \Vert_{L^2}\\
\leq & \frac{\alpha^2}{2 \beta} \Vert \fe \Vert^{2}_{\dot{H}^1} + \frac{\beta \alpha^2}{2} \Vert \U \Vert^{2}_{\dot{H}^1}+ \frac{1}{2 \beta} \Vert \fe \Vert^{2}_{L^2}+ \frac{\beta}{2} \Vert \U \Vert^{2}_{L^2}.
\end{split}
\end{equation}
Getting back to the last identity, we obtain 
\begin{equation}\label{eq-aux-1}
 \frac{\nu\alpha^2}{2}  \Vert \U \Vert^{2}_{\dot{H}^2}+ (\nu+\frac{\beta \alpha^2}{2}) \Vert \U \Vert^{2}_{\dot{H}^{1}}+ \frac{\beta}{2}\beta \Vert \U \Vert^{2}_{L^2} \leq \frac{1}{2 \beta}\left(\Vert \fe \Vert^{2}_{L^2}+\alpha^2\Vert \fe \Vert^{2}_{\dot{H}^1}\right),
\end{equation}
hence we have
\begin{equation*}
\beta \Vert \U \Vert^{2}_{L^2} + \nu\alpha^2 \Vert \U \Vert^{2}_{\dot{H}^2} \leq  \frac{1}{\beta}\left(\Vert \fe \Vert^{2}_{L^2}+\alpha^2\Vert \fe \Vert^{2}_{\dot{H}^1}\right). 
\end{equation*}
Finally,  we can write  
\begin{equation}\label{unif-bound-Ur}
\Vert \U \Vert^{2}_{H^2} \leq \frac{1}{\min(\beta,\nu \alpha^2)\beta} \left(\Vert \fe \Vert^{2}_{L^2}+\alpha^2\Vert \fe \Vert^{2}_{\dot{H}^1}\right),     
\end{equation}
where, setting the constant $\ds{M^2= \frac{1}{\min(\beta,\nu \alpha^2)\beta} \left(\Vert \fe \Vert^{2}_{L^2}+\alpha^2\Vert \fe \Vert^{2}_{\dot{H}^1}\right)}$,  the point $3)$ is now verified.\\

We may apply the Theorem \ref{Th-Sheafer} to obtain a solution $\U_r \in E$ of the problem $\U_r = T_r(\U_r)$, for all $r>0$ fixed.  We will show now that the family  $(\U_r)_{r>0}$ converges (in the distributional sense) to a solution of the equation (\ref{Bardina-Stat}) when the parameter $r$ goes (via a sub-sequence) to infinity. \\

We observe first that by estimate (\ref{unif-bound-Ur}) the family $(\U_r)_{r>0}$ in uniformly bounded in the space $H^2(\Rt)$. Thus, for all $\varphi \in \mathcal{C}^{\infty}_{0}(\Rt)$ we have: $\ds{ \sup_{r>0} \Vert \varphi\, \U_r \Vert_{H^2} <+\infty}$, and then, by the Rellich-Lions lemma  there exists $\U \in H^{2}_{loc}(\Rt)$ and a sub-sequence $(r_n)_{n\in \mathbb{N}}$, such that $r_n \to +\infty$ when $n\to +\infty$, and such  that $\U_{r_n}$ strongly converges to $\U$ in the space $L^{p}_{loc}(\Rt)$ for $2\leq p < +\infty$. Moreover, as  $div(\U_{r_n})=0$ for all $n\in \mathbb{N}$ we have $div(\U)=0$, and always by the uniform boundness of the sequence $\U_{r_n}$ in the space $H^2(\Rt)$, by the Banach-Alaoglu theorem we obtain that $\U \in H^2(\Rt)$.\\

In order to prove that  $\U \in H^2(\Rt)$  is a solution of the equation (\ref{Bardina-Stat}), we recall  that as we have  $\U_{r_n}=T_{r_n}(\U_{r_n})$, then $\U_{r_n}$ solves the equation 
\begin{equation*}
 \nu\alpha^2 \Delta^{2} \U_{r_n} -(\nu+\beta \alpha^2) \Delta \U_{r_n} +\beta \U_{r_n} + \P  (\theta_{r_n} \U_{r_n} \cdot \vec{\nabla}) (\theta_{r_n}\U_{r_n}) = -\alpha^2  \Delta\fe + \fe. 
\end{equation*}
Moreover, since we have $div(\U_{r_n})=0$ and $div(\fe)=0$, then $\U_{r_n}$ verifies the equation  
\begin{equation*}
\P\left(  \nu\alpha^2 \Delta^{2} \U_{r_n} -(\nu+\beta \alpha^2) \Delta \U_{r_n} +\beta \U_{r_n} +  (\theta_{r_n} \U_{r_n} \cdot \vec{\nabla}) (\theta_{r_n}\U_{r_n}) +\alpha^2  \Delta\fe - \fe \right)=0, 
\end{equation*}
where, by well-known properties of the Leray's projector, for all $n \in \mathbb{N}$ we have 
\begin{equation}\label{eq-rot}
\vec{\nabla}\wedge \left(   \nu\alpha^2 \Delta^{2} \U_{r_n} -(\nu+\beta \alpha^2) \Delta \U_{r_n} +\beta \U_{r_n} +  (\theta_{r_n} \U_{r_n} \cdot \vec{\nabla}) (\theta_{r_n}\U_{r_n}) +\alpha^2  \Delta\fe - \fe \right)=0. 
\end{equation}

Here,  we will prove that the non-linear term $\ds{(\theta_{r_n} \U_{r_n} \cdot \vec{\nabla}) (\theta_{r_n}\U_{r_n})}$ converges in the distributional sense to the non-linear term $\ds{( \U \cdot \vec{\nabla})\U}$ when $n\to +\infty$. Indeed, we observe first that as $div(\U_{r_n})=0$, then we can write 
\begin{equation*}
 (\theta_{r_n} \U_{r_n} \cdot \vec{\nabla}) (\theta_{r_n}\U_{r_n})= \theta_{r_n} \left(  ( \U_{r_n} \cdot \vec{\nabla}) (\theta_{r_n}\U_{r_n})\right)=\theta_{r_n} div \left( \U_{r_n} \otimes (\theta_{r_n}\U_{r_n})\right).   
\end{equation*}
Thereafter, as $\U_{r_n}$ strongly converges to $\U$ in the space $L^{4}_{loc}(\Rt)$, and moreover, as we have $\theta_{r_n}(x)=1$ when $\vert x \vert < r_n$, then $\theta_{r_n}\U_{r_n}$ also strongly converges to $\U$ in $L^{4}_{loc}(\Rt)$. We get that $\ds{\U_{r_n} \otimes (\theta_{r_n}\U_{r_n})}$ converges to $\ds{\U \otimes \U}$ in the strong topology of the space $L^{2}_{loc}(\Rt)$, hence we conclude the desired convergence. We thus have the following limit in the distributional sense: 
\begin{equation*}
\begin{split}
\lim_{r_{n}\to +\infty} & \left(  \nu\alpha^2 \Delta^{2} \U_{r_n} -(\nu+\beta \alpha^2) \Delta \U_{r_n} +\beta \U_{r_n} +  (\theta_{r_n} \U_{r_n} \cdot \vec{\nabla}) (\theta_{r_n}\U_{r_n}) +\alpha^2  \Delta\fe - \fe \right)\\
& =  \nu\alpha^2 \Delta^{2} \U -(\nu+\beta \alpha^2) \Delta \U +\beta \U +  (\U \cdot \vec{\nabla}) \U +\alpha^2  \Delta\fe - \fe. 
\end{split}
\end{equation*}

Thus, by (\ref{eq-rot}) we get 
\begin{equation*}
\vec{\nabla}\wedge \left(  \nu\alpha^2 \Delta^{2} \U -(\nu+\beta \alpha^2) \Delta \U +\beta \U +  (\U \cdot \vec{\nabla}) \U +\alpha^2  \Delta\fe - \fe \right)=0,
\end{equation*}
and then, there exists $Q \in \mathcal{D}'(\Rt)$ such that 
\begin{equation*}
\nu\alpha^2 \Delta^{2} \U -(\nu+\beta \alpha^2) \Delta \U +\beta \U +  (\U \cdot \vec{\nabla}) \U +\alpha^2  \Delta\fe - \fe = \vec{\nabla}Q.
\end{equation*}
Moreover,  as $\U \in H^2(\Rt)$ and $\fe \in H^1(\Rt)$, we observe that the term in the left side of this equation belongs to the space $H^{-2}(\Rt)$, hence $\vec{\nabla}Q \in H^{-2}(\Rt)$. \\

From this equation we can write 
\begin{equation}\label{eq-stat-ext}
-\nu (-\alpha^2 \Delta+I_d)  \Delta \U +\beta (-\alpha^2 \Delta+I_d) \U  +  (\U \cdot \vec{\nabla}) \U - \vec{\nabla}Q = (-\alpha^2 \Delta^{2}+I_d)\fe,
\end{equation}
then, applying the filtering operator $\ds{(\cdot )_\alpha = (-\alpha^2 \Delta+I_d)^{-1}}$ to each term, and defining the pressure $P$ as $P= -(-\alpha^2 \Delta^{2}+I_d)^{-1} Q \in H^1(\Rt)$, we finally obtain that the couple $(\U, P)$ is a solution of equation (\ref{Bardina-Stat}).

\subsubsection*{2) The energy estimate} 

As we have $\U \in H^2(\Rt)$ and $P\in H^1(\Rt)$, then we can multiply each term in the equation (\ref{eq-stat-ext})  by $\U$. Thereafter,  we integrate  on $\Rt$ and after some integration by parts  we obtain

\begin{equation*}
 \nu\alpha^2  \Vert \U \Vert^{2}_{\dot{H}^2}+ (\nu+\beta \alpha^2) \Vert \U \Vert^{2}_{\dot{H}^{1}}+ \beta \Vert \U \Vert^{2}_{L^2} = \alpha^2\,\int_{\Rt} \vec{\nabla}\otimes \fe \cdot \vec{\nabla} \otimes\U \, dx + \, \int_{\Rt} \fe \cdot \U\, dx.     
\end{equation*}
Moreover, the term in the right side was estimated in (\ref{estim-term-force}) (where we set $\lambda=1$) and we have (\ref{eq-aux-1}). From this estimate we can write now 
\begin{equation*}
\frac{\beta}{2} \Vert \U \Vert^{2}_{L^2} +  \frac{\beta \alpha^2}{2} \Vert \U \Vert^{2}_{\dot{H}^{1}} + \frac{\nu\alpha^2}{2}  \Vert \U \Vert^{2}_{\dot{H}^2} \leq \frac{1}{2 \beta}\left(\Vert \fe \Vert^{2}_{L^2}+\alpha^2\Vert \fe \Vert^{2}_{\dot{H}^1}\right),
\end{equation*}
hence, recalling that $\Vert \cdot \Vert^{2}_{\mathcal{H}^{1}_{\alpha}}=\Vert \cdot \Vert^{2}_{L^2}+\alpha^2 \Vert \cdot \Vert^{2}_{\dot{H}^1}$, we finally obtain the desired energy estimate \[ \Vert \U \Vert^{2}_{\mathcal{H}^{1}_{\alpha}} + \nu \alpha^2 \Vert \U \Vert^{2}_{\dot{H}^2}\leq \frac{2}{\beta^2} \Vert \fe \Vert^{2}_{\mathcal{H}^{1}_{\alpha}}.\]  

\subsubsection*{3) All the stationary solutions belong to the global attractor} 

To prove this point, we just recall that by (\ref{Internal-estructure}) the global attractor $\mathcal{A}_{\fe}\subset \mathcal{H}^{1}_{\alpha}(\Rt)$ is the set of functions $\vv(0,\cdot)$, where $\vv \in L^{\infty}\Big(\R, \mathcal{H}^{1}_{\alpha}(\Rt)\Big)$ is  an eternal solution of the equation (\ref{Bardina-enternal}). On the other hand, by the energy  estimate proven in the point $2)$, we know that all the  solutions $\U \in H^2(\Rt)$ of the stationary problem (\ref{Bardina-Stat}) verify $\U \in \mathcal{H}^{1}_{\alpha}(\Rt)$, and moreover,  as they do not depend on the time variable, they also are a particular case of  bounded in time eternal solutions, and we can write $\U(0,\cdot)=\U$.  Consequently, all the  solutions $\U \in H^2(\Rt)$ of the equation  (\ref{Bardina-Stat}) belong to the global attractor $\mathcal{A}_{\fe}$.  Theorem \ref{Th-Stationary-Solutions} is now proven. \finpv
 
\section{Additional properties of the global attractor driven by the damping parameter}\label{Sec:AdditionalProp}

\subsection*{Proof of Theorem \ref{Th:Stability}}

The proof of this theorem bases on the following result concerning the long-time behavior of two solutions of the equation (\ref{Bardina}).

\begin{Proposition}\label{Prop:long-time-behavior}
Let $\fe_{1}, \fe_{2} \in \mathcal{H}^{1}_{\alpha}(\Rt)$ be two  divergence-free external forces and  let  $\vu_{0,1}, \vu_{0,2} \in \mathcal{H}^{1}_{\alpha}(\Rt)$ be two initial data.  Moreover,  let $(\vu_1, p_1)\, (\vu_2, p_2)\in L^{\infty}_{t}\,\mathcal{H}^{1}_{\alpha} \cap (L^{2}_{loc})_{t}\,\dot{H}^{2}_{x} \times (L^{2}_{loc})_{t}\,H^{1}_{x}$, be two global in time  solutions of the equation (\ref{Bardina})  arising from  the data $(\vu_{0,1},\fe_1)$ and $(\vu_{0,2},\fe_2)$ respectively.\\ 

Moreover, for a numerical constant $c>0$, and the parameters $\alpha>0$ and $\beta>0$, we define the quantity $\eta(\beta)$ as: 
\begin{equation*}
\eta(\beta)= 2\left( -\frac{\beta}{2}+ \frac{c}{\alpha^{5/2}\beta}\Vert \fe_1 \Vert_{\mathcal{H}^{1}_{\alpha}}\right),
\end{equation*}
Then, for all time $t\geq 0$, the following estimate holds: 
\begin{equation}\label{Estim-Stab}
\Vert \vu_1(t,\cdot) - \vu_{2}(t,\cdot)\Vert^{2}_{\mathcal{H}^{1}_{\alpha}} \leq \Vert \vu_{0,1}-\vu_{0,2}\Vert^{2}_{\mathcal{H}^{1}_{\alpha}}e^{\eta(\beta)\, t}+ \frac{1}{\beta} \Vert \fe_1 - \fe_2 \Vert^{2}_{\mathcal{H}^{1}_{\alpha}}\,\frac{1}{\eta(\beta)}\Big(e^{\eta(\beta)\, t} -1\Big). 
\end{equation}
\end{Proposition}
\pv From the solutions $(\vu_1, p_1) \in L^{\infty}_{t}\,\mathcal{H}^{1}_{\alpha} \cap (L^{2}_{t})_{loc}\,\dot{H}^{2}_{x} \times (L^{2}_{t})_{loc}H^{1}_{x}$ and $(\vu_2, p_2) \in L^{\infty}_{t}\,\mathcal{H}^{1}_{\alpha} \cap (L^{2}_{t})_{loc}\,\dot{H}^{2}_{x} \times (L^{2}_{t})_{loc}H^{1}_{x}$, we define  $\vw=\vu_1-\vu_2$ and $q=p_1-p_2$. We observe that $(\vw,q) \in L^{\infty}_{t}\,\mathcal{H}^{1}_{\alpha} \cap (L^{2}_{t})_{loc}\dot{H}^{2}_{x} \times (L^{2}_{t})_{loc}H^{1}_{x}$, solves the equation: 
\begin{equation*}
\partial_t \vw + \left((\vw \cdot \vec{\nabla}) \vu_1+ (\vu_2\cdot \vec{\nabla})\vw\right)_\alpha- \nu \Delta \vw + \vec{\nabla} q= \fe_1 - \fe_2 -\beta \vw, \quad div(\vw)=0, \quad \vw(0,\cdot)= \vu_{1}(0,\cdot)-\vu_2(0,\cdot).  
\end{equation*}
Moreover, performing  the same computations done in (\ref{Eq-dif-energ}), for  $t \geq 0$, we have the following energy equality: 
\begin{equation}\label{estim-energ-stab} 
\begin{split}
 \frac{1}{2}\frac{d}{dt} \, \Vert \vw(t,\cdot)\Vert^{2}_{\mathcal{H}^{1}_{\alpha}} =& 	-\nu \Vert \vw(t,\cdot)\Vert^{2}_{\dot{H}^1} - \alpha^2 \Vert \vw(t,\cdot)\Vert^{2}_{\dot{H}^2} - \beta  \Vert \vw(t,\cdot)\Vert^{2}_{\mathcal{H}^{1}_{\alpha}} - \left\langle \fe_1-\fe_2, \vw \right\rangle_{L^2 \times L^2}\\
 & + \alpha^2 \left\langle \vec{\nabla}\otimes (\fe_1-\fe_2), \vec{\nabla}\otimes \vw \right\rangle_{L^2 \times L^2}  - \left\langle (\vw \cdot \vec{\nabla}) \vu_1(t,\cdot) , \vw(t,\cdot)\right\rangle_{\dot{H}^{-1}\times \dot{H}^1}.
 \end{split}
\end{equation} We study now  the term $\ds{\left\langle (\vw \cdot \vec{\nabla}) \vu_1(t,\cdot) , \vw(t,\cdot)\right\rangle_{\dot{H}^{-1}\times \dot{H}^1}}$.  More precisely, as $div(\vw)=0$ and integrating by parts, we can write 
\begin{equation*}
\begin{split}
\left\langle (\vw \cdot \vec{\nabla}) \vu_1(t,\cdot) , \vw(t,\cdot)\right\rangle_{\dot{H}^{-1}\times \dot{H}^1}= & \sum_{i,j=1}^{3} \left\langle w_j (\partial_j u_{1,i}), w_{i})  \right\rangle_{\dot{H}^{-1}\times \dot{H}^1}= \sum_{i,j=1}^{3} \left\langle (\partial_j w_j  u_{1,i}), w_{i})  \right\rangle_{\dot{H}^{-1}\times \dot{H}^1} \\
=& - \sum_{i,j=1}^{3} \left\langle  w_j  u_{1,i}, \partial_j w_{i})  \right\rangle_{L^2 \times L^2} =  - \sum_{i,j=1}^{3} \int_{\Rt}  w_j \, u_{1,i} \, \partial_j w_{i} \, dx. 
\end{split}    
\end{equation*}
By the Parseval's identity and the Cauchy-Schwarz inequality we obtain 
\begin{equation*}
\begin{split}
 - \sum_{i,j=1}^{3} \int_{\Rt}  w_j  u_{1,i} \, \partial_j w_{i} \, dx =& - \sum_{i,j=1}^{3} \int_{\Rt}  \widehat{w_j} \,\, \widehat{u_{1,i}  \partial_j w_{i}} \, d\xi =  - \sum_{i,j=1}^{3} \int_{\Rt}  \vert \xi \vert \widehat{w_j} \,\, \vert \xi \vert^{-1} ( \widehat{u_{1,i}  \partial_j w_{i}})  \, d\xi \\
\leq & \Vert \vw(t,\cdot) \Vert_{\dot{H}^1} \Vert \vu_1(t,\cdot) (\vec{\nabla}\otimes \vw)(t,\cdot) \Vert_{\dot{H}^{-1}}.  
\end{split}    
\end{equation*}
Then, applying the Hardy-Littlewood-Sobolev inequalities,  the H\"older inequalities (with $5/6=1/3+1/2$), and moreover, recalling that $\Vert \cdot \Vert^{2}_{\mathcal{H}^{1}_{\alpha}}=\Vert \cdot \Vert^{2}_{L^2}+\alpha^2 \Vert \cdot \Vert_{\dot{H}^1}$, we have 
\begin{equation*}
\begin{split}
\Vert \vw(t,\cdot) \Vert_{\dot{H}^1} \Vert \vu_1 (\vec{\nabla}\otimes \vw) \Vert_{\dot{H}^{-1}} \leq & c \Vert \vw(t,\cdot) \Vert_{\dot{H}^1} \Vert \vu_1 (\vec{\nabla}\otimes \vw) \Vert_{L^{6/5}} \leq c  \Vert \vw(t,\cdot) \Vert_{\dot{H}^1}  \Vert \vu_1(t,\cdot)\Vert_{L^3} \Vert  \vec{\nabla}\otimes \vw (t,\cdot)\Vert_{L^2}\\
\leq & c\Vert \vw(t,\cdot) \Vert^{2}_{\dot{H}^1}  \Vert \vu_1(t,\cdot)\Vert_{L^3} \leq c\, \alpha^2 \Vert \vw(t,\cdot) \Vert^{2}_{\dot{H}^1} \, \frac{1}{\alpha^2}  \Vert \vu_1(t,\cdot)\Vert_{L^3} \\
\leq &c\, \Vert \vw(t,\cdot) \Vert^{2}_{\mathcal{H}^{1}_{\alpha}} \, \frac{1}{\alpha^2}  \Vert \vu_1(t,\cdot)\Vert_{L^3}. 
\end{split}    
\end{equation*}
We still need to estimate the term $\ds{\frac{1}{\alpha^2}  \Vert \vu_1(t,\cdot)\Vert_{L^3}}$. By the interpolation inequalities (with $1/3=\theta/2+ (1-\theta)/ 6$, and $\theta=1/2$) and applying again the Hardy-Littlewood-Sobolev inequalities, we can write  
\begin{equation*}
\begin{split}
\frac{1}{\alpha^2}  \Vert \vu_1(t,\cdot)\Vert_{L^3} \leq & \frac{c}{\alpha^2} \Vert \vu_1(t,\cdot)\Vert^{1/2}_{L^2}\, \Vert \vu_1(t,\cdot)\Vert^{1/2}_{L^6} \leq  \frac{c}{\alpha^2} \Vert \vu_1(t,\cdot)\Vert^{1/2}_{L^2}\, \Vert \vu_1(t,\cdot)\Vert^{1/2}_{\dot{H}^1} \\
\leq & \frac{c}{\alpha^{5/2}} \Vert \vu_1(t,\cdot)\Vert^{1/2}_{L^2}\,\alpha^{1/2} \Vert \vu_1(t,\cdot)\Vert^{1/2}_{\dot{H}^1}\leq  \frac{c}{\alpha^{5/2}} \Vert \vu_1(t,\cdot)\Vert_{\mathcal{H}^{1}_{\alpha}}.
\end{split}    
\end{equation*}
Now, by the point $1)$ in the Proposition \ref{Prop:time-control}  we obtain 
\begin{equation*}
  \frac{c}{\alpha^{5/2}} \Vert \vu_1(t,\cdot)\Vert_{\mathcal{H}^{1}_{\alpha}} \leq \frac{c}{\alpha^{5/2}} \left( \Vert \vu_{0,1} \Vert_{\mathcal{H}^{1}_{\alpha}} e^{- \frac{\beta}{2}\, t} + \frac{2}{\beta} \Vert \fe_1 \Vert_{\mathcal{H}^{1}_{\alpha}} \right).   
\end{equation*}
On the other hand, we recall that  for any external force $\fe \in \mathcal{H}^{1}_{\alpha}(\Rt)$, by Lemma \ref{LemaTech2-Absorbing set} we have that the set  $\ds{\mathcal{B}=\left\{ \vu \in \mathcal{H}^{1}_{\alpha}(\Rt): \Vert \vu \Vert^{2}_{\mathcal{H}^{1}_{\alpha}} \leq \frac{8}{\beta^2} \Vert \fe \Vert^{2}_{\mathcal{H}^{1}_{\alpha}}\right\}}$, is a absorbing set in the sense of Definition \ref{def-Abosorbing-set}. Then, for any initial datum $\vu_{0} \in \mathcal{H}^{1}_{\alpha}(\Rt)$, the solution $\vu(t,x)$ arising from $(\vu_{0}, \fe)$ always lies in the set $\mathcal{B}$ from a time large enough. Thus, getting back to the solution $\vu_1(t,x)$,  without lost  of generality we may suppose that the initial  datum $\vu_{0,1}$ belongs to the absorbing set $\mathcal{B}=\left\{ \vu \in \mathcal{H}^{1}_{\alpha}(\Rt): \Vert \vu \Vert^{2}_{\mathcal{H}^{1}_{\alpha}} \leq \frac{8}{\beta^2} \Vert \fe_1 \Vert^{2}_{\mathcal{H}^{1}_{\alpha}}\right\}$, and  we can write  the estimate $\ds{\Vert \vu_{0,1}\Vert_{\mathcal{H}^{1}_{\alpha}} \leq \frac{\sqrt{8}}{\beta} \Vert \fe_1 \Vert_{\mathcal{H}^{1}_{\alpha}}}$.  \\

We thus have 
\begin{equation*}
   \frac{c}{\alpha^{5/2}} \Vert \vu_1(t,\cdot)\Vert_{\mathcal{H}^{1}_{\alpha}} \leq \frac{c}{\alpha^{5/2}} \left( \frac{\sqrt{8}}{\beta} \Vert \fe_1 \Vert_{\mathcal{H}^{1}_{\alpha}}\, e^{-\frac{\beta}{2} t} + \frac{2}{\beta} \Vert \fe_1 \Vert_{\mathcal{H}^{1}_{\alpha}} \right)  \leq   \frac{c}{\alpha^{5/2} \beta^2} \, \Vert \fe_1 \Vert_{\mathcal{H}^{1}_{\alpha}}, 
\end{equation*} and then,  we can write 
\begin{equation*}
\frac{1}{\alpha^2}  \Vert \vu_1(t,\cdot)\Vert_{L^3} \leq \frac{c}{\alpha^{5/2} \beta} \, \Vert \fe_1 \Vert_{\mathcal{H}^{1}_{\alpha}}. 
\end{equation*}
Finally, gathering these inequalities we get the following estimate 
\begin{equation*}
\left\vert \left\langle (\vw \cdot \vec{\nabla}) \vu_1(t,\cdot) , \vw(t,\cdot)\right\rangle_{\dot{H}^{-1}\times \dot{H}^1} \right\vert \leq  \Vert \vw(t,\cdot) \Vert^{2}_{\mathcal{H}^{1}_{\alpha}} \,   \frac{c}{\alpha^{5/2} \beta} \, \Vert \fe_1 \Vert_{\mathcal{H}^{1}_{\alpha}}.  
\end{equation*}

With this estimate at hand, we get back to the energy equality (\ref{estim-energ-stab}) to write 
\begin{equation*}
\begin{split}
  \frac{1}{2}\frac{d}{dt} \, \Vert \vw(t,\cdot)\Vert^{2}_{\mathcal{H}^{1}_{\alpha}} \leq & 	-\nu \Vert \vw(t,\cdot)\Vert^{2}_{\dot{H}^1} - \alpha^2 \Vert \vw(t,\cdot)\Vert^{2}_{\dot{H}^2} - \beta  \Vert \vw(t,\cdot)\Vert^{2}_{\mathcal{H}^{1}_{\alpha}} - \left\langle \fe_1-\fe_2, \vw \right\rangle_{L^2 \times L^2}\\
 & + \alpha^2 \left\langle \vec{\nabla}\otimes (\fe_1-\fe_2), \vec{\nabla}\otimes \vw \right\rangle_{L^2 \times L^2} + \Vert \vw(t,\cdot) \Vert^{2}_{\mathcal{H}^{1}_{\alpha}} \,   \frac{c}{\alpha^{5/2} \beta} \, \Vert \fe_1 \Vert_{\mathcal{H}^{1}_{\alpha}}.   
\end{split}    
\end{equation*}
As the first and second term in the right side are negatives, we have 
\begin{equation*}
\begin{split}
  \frac{1}{2}\frac{d}{dt} \, \Vert \vw(t,\cdot)\Vert^{2}_{\mathcal{H}^{1}_{\alpha}} \leq &  - \beta  \Vert \vw(t,\cdot)\Vert^{2}_{\mathcal{H}^{1}_{\alpha}} - \left\langle \fe_1-\fe_2, \vw \right\rangle_{L^2 \times L^2}  + \alpha^2 \left\langle \vec{\nabla}\otimes (\fe_1-\fe_2), \vec{\nabla}\otimes \vw \right\rangle_{L^2 \times L^2} \\
  & + \Vert \vw(t,\cdot) \Vert^{2}_{\mathcal{H}^{1}_{\alpha}} \,   \frac{c}{\alpha^{5/2} \beta} \, \Vert \fe_1 \Vert_{\mathcal{H}^{1}_{\alpha}}.   
\end{split}    
\end{equation*}
Then, applying first the Cauchy-Schwarz inequalities and thereafter the Young inequalities, we get 
\begin{equation*}
\begin{split}
  \frac{1}{2}\frac{d}{dt} \, \Vert \vw(t,\cdot)\Vert^{2}_{\mathcal{H}^{1}_{\alpha}} \leq &  - \beta  \Vert \vw(t,\cdot)\Vert^{2}_{\mathcal{H}^{1}_{\alpha}} + \frac{1}{2\beta} \Vert \fe_1 - \fe_2 \Vert^{2}_{L^2}+ \frac{\beta}{2}\Vert \vw(t,\cdot)\Vert^{2}_{L^2}+ \frac{ \alpha^2}{2\beta}\Vert \fe_1-\fe_2 \Vert^{2}_{\dot{H}^1}+ \frac{\alpha^2 \beta}{2}\Vert \vw(t,\cdot)\Vert^{2}_{\dot{H}^1}\\
  & + \Vert \vw(t,\cdot) \Vert^{2}_{\mathcal{H}^{1}_{\alpha}} \,   \frac{c}{\alpha^{5/2} \beta} \, \Vert \fe_1 \Vert_{\mathcal{H}^{1}_{\alpha}} \\
 \leq &-\beta \Vert \vw(t,\cdot)\Vert^{2}_{\mathcal{H}^{1}_{\alpha}} +  \frac{\beta}{2} \left( \Vert \vw(t,\cdot)\Vert^{2}_{L^2}+\alpha^2 \Vert \vw(t,\cdot)\Vert^{2}_{\dot{H}^1}\right) + \frac{1}{2\beta}\left( \Vert \fe_1\Vert^{2}_{L^2}+\alpha^2 \Vert \fe_1\Vert^{2}_{\dot{H}^1} \right) \\  
  & + \Vert \vw(t,\cdot) \Vert^{2}_{\mathcal{H}^{1}_{\alpha}} \,   \frac{c}{\alpha^{5/2} \beta} \, \Vert \fe_1 \Vert_{\mathcal{H}^{1}_{\alpha}} \\
 \leq & \left( -\frac{\beta}{2}+ \frac{c}{\alpha^{5/2}\beta}\Vert \fe_1 \Vert_{\mathcal{H}^{1}_{\alpha}}\right)\Vert \vw(t,\cdot)\Vert^{2}_{\mathcal{H}^{1}_{\alpha}}+ \frac{1}{2\beta} \Vert \fe_1 - \fe_2 \Vert^{2}_{\mathcal{H}^{1}_{\alpha}},
\end{split}    
\end{equation*}
hence, we write 
\begin{equation*}
\begin{split}
 \frac{d}{dt} \, \Vert \vw(t,\cdot)\Vert^{2}_{\mathcal{H}^{1}_{\alpha}} \leq  2\left( -\frac{\beta}{2}+ \frac{c}{\alpha^{5/2}\beta}\Vert \fe_1 \Vert_{\mathcal{H}^{1}_{\alpha}}\right) \Vert \vw(t,\cdot)\Vert^{2}_{\mathcal{H}^{1}_{\alpha}} + \frac{1}{\beta} \Vert \fe_1-\fe_2 \Vert^{2}_{\mathcal{H}^{1}_{\alpha}}.  
\end{split}    
\end{equation*}
We set now the quantity $\ds{\eta(\beta)= 2\left( -\frac{\beta}{2}+ \frac{c}{\alpha^{5/2}\beta}\Vert \fe_1 \Vert_{\mathcal{H}^{1}_{\alpha}}\right)}$, and using the Gr\"onwall inequalities we finally obtain the desired estimate (\ref{Estim-Stab}). Proposition \ref{Prop:long-time-behavior} is now proven. \finpv

With Proposition \ref{Prop:long-time-behavior} at hand, we are able to prove each point stated in Theorem \ref{Th:Stability}. For this, we recall the definition of the expression $\eta(\beta)$ given in (\ref{eta}). 

\subsection*{1) The orbital stability when $\eta(\beta)=0$.} 
In the framework of Proposition \ref{Prop:long-time-behavior}, firs we set  $\fe_1=\fe_2=\fe$, and we  get the estimate
\begin{equation*}
\Vert \vu_1(t,\cdot) - \vu_{2}(t,\cdot)\Vert^{2}_{\mathcal{H}^{1}_{\alpha}} \leq \Vert \vu_{0,1}-\vu_{0,2}\Vert^{2}_{\mathcal{H}^{1}_{\alpha}}e^{\eta(\beta)\, t}. 
\end{equation*}
Then, we take $\eta(\beta)=0$ to obtain the following control
\begin{equation*} 
\Vert \vu_1(t,\cdot) - \vu_{2}(t,\cdot)\Vert^{2}_{\mathcal{H}^{1}_{\alpha}} \leq \Vert \vu_{0,1}-\vu_{0,2}\Vert^{2}_{\mathcal{H}^{1}_{\alpha}},
\end{equation*}
hence, the result stated in this point follows directly. 

\subsection*{2) The characterization of the global attractor when $\eta(\beta)<0$.}

In the first step, we will prove that the uniqueness of the stationary solution constructed in  point $1)$ of the Theorem \ref{Th-Stationary-Solutions}. Let $(\U_1, P_1), \, (\U_2, P_2)\in H^2(\Rt)\times H^1(\Rt)$ be two solutions of the stationary problem (\ref{Bardina-Stat}) associated with the same external force $\fe$.  As $\U_1$ and $\U_2$ are time-independing functions we have $\partial_t \U_1=0$ and $\partial_t \U_2=0$, and thus, $(\U_1, P_1), \, (\U_2, P_2)$ are also two solutions of the equation (\ref{Bardina}), arising from the initial data $\U_1$ and $\U_2$ respectively and with external force $\fe$. Moreover, we also have  $(\U_1, P_1), \, (\U_2, P_2)\in L^{\infty}_{t}\, H^{1}_{x} \cap (L^{2}_{loc})_{t}\, \dot{H}^{2}_{x} \times (L^{2}_{loc})_{t}\, H^{1}_{x}$. \\

By Proposition \ref{Prop:long-time-behavior} ( with $\fe_1=\fe_2=\fe$) we can write: $\ds{\Vert \U_1 - \U_2 \Vert_{\mathcal{H}^{1}_{\alpha}} \leq \Vert \U_1 - \U_2 \Vert_{\mathcal{H}^{1}_{\alpha}} e^{\eta(\beta)\, t}}$. Moreover, as we have  $\eta(\beta) <0$,  for a time $t>0$ large enough we can write $\ds{\Vert \U_1 - \U_2 \Vert_{\mathcal{H}^{1}_{\alpha}} \leq \frac{1}{2} \Vert \U_1 - \U_2 \Vert_{\mathcal{H}^{1}_{\alpha}}}$, hence we obtain the identity $\U_1=\U_2$. Finally, as the pressure term is always related to the velocity field by (\ref{Caract-pression}) we also obtain the identity $P_1=P_2$.  \\

In the second step, we prove now   that the singleton  $\{ \U \}$ is also a global attractor for the equation (\ref{Bardina}) given in Definition \ref{Def:atractor}.  The first point in Definition \ref{Def:atractor} is evident,  so we will focus on the second point.  Let $t\geq 0$. We will prove the identity $S(t)\{ \U \} = \{ \U \}$. Let $\vu \in S(t)\{ \U \}$, \emph{i.e.}, $\vu(t,\cdot)$ is the solution of the equation (\ref{Bardina}), at the time $t$ and  arising from the initial datum $\U$. But, as explained above, we also have that $\U$ is a solution of the equation (\ref{Bardina}) which arises from $\U$. Then,  by Proposition \ref{Prop:long-time-behavior} ( with $ \vu_{1,0}=\vu_{0,2}=\U$ and $\fe_1=\fe_2=\fe$) we can write 
$\ds{\Vert \vu(t,\cdot)-\U \Vert_{\mathcal{H}^{1}_{\alpha}} \leq \Vert \U - \U \Vert_{\mathcal{H}^{1}_{\alpha}}e^{\eta(\beta)\, t}}$, hence we get $\vu=\U$. On the other hand, we also have   $\U \in S(t)\{\U\}$, and it  directly  follows from the fact that taking $\U$ as an initial datum in equation (\ref{Bardina}) then, for all time $t\geq 0$, the same function $\U$ is the unique solution for this problem. We verify now the  third point in Definition \ref{Def:atractor}. Let $\vu_0 \in \mathcal{H}^{1}_{\alpha}(\Rt)$ an initial datum and let $\vu \in L^{\infty}_{t}\,\mathcal{H}^{1}_{\alpha} \cap (L^{2}_{loc})_{t}\,\dot{H}^{2}_{x}$ be the unique solution of equation (\ref{Bardina}), arising from $\vu_0$ and associated with the force $\fe$, which is  given by Theorem \ref{Th-WP}. On the other hand, let $\U \in H^2(\Rt)$ be the unique stationary solution of the problem  (\ref{Bardina}) associated with the same force $\fe$. Thus, always by Proposition \ref{Prop:long-time-behavior} (with $\fe_1=\fe_2=\fe$), for all time $t\geq 0$ we have the estimate 
\begin{equation}\label{key-estim}
 \Vert \vu(t,\cdot)-\U \Vert_{\mathcal{H}^{1}_{\alpha}} \leq \Vert \vu_{0}-\U \Vert_{\mathcal{H}^{1}_{\alpha}} \, e^{\eta(\beta), t}, \quad \eta(\beta)<0.    
\end{equation}
 By this estimate, we may observe that the unique stationary solution $\U$ is asymptotically stable. Actually we have a stronger stability property in the sense that, for any initial datum $\vu_0$ the unique solution $\vu(t,\cdot)$ arising from  $\vu_0$ strongly converges to the stationary solution $\U$ when $t \to +\infty$, and thus, the third point in Definition \ref{Def:atractor} holds.\\

In the third and last step, we will prove the identity $\{ \U \}=\mathcal{A}_{\fe}$. Indeed,  on the one hand, we have that $\{ \U \}$ is  a global attractor of the equation (\ref{Bardina}). On the other hand, recall that by the Theorem \ref{Th-Atractor} we also have the global attractor $\mathcal{A}_{\fe}$ given by this theorem. But, by Lemma $2.18$, page 16 in \cite{Raugel}, we have the uniqueness of the global attractor, provided it satisfies all the points in Definition \ref{Def:atractor}, hence we conclude that $\{ \U \}=\mathcal{A}_{\fe}$. Theorem \ref{Th:Stability} is now proven. \finpv

\subsection*{Proof of Proposition \ref{Prop:profile}} 
Let $\U \in \mathcal{A}_{\fe}$ be the unique solution of the equation (\ref{Bardina-Stat}), and let $\vu \in L^{\infty}_{t}\,(\mathcal{H}^{1}_\alpha)_{x} \cap (L^{2}_{t})_{loc}\,\dot{H}^{2}_{x}$ be a solution of the equation (\ref{Bardina}) arising from an initial datum $\vu_0 \in \mathcal{H}^{1}_{\alpha}(\Rt)$. We define the term 
$\ds{\mathcal{R}_{\vu}=\vu-\U}$ which, to simplify the notation, we shall write as $\mathcal{R}$. \\ 

In order to prove (\ref{time-decay1}), we observe first that as the stationary solution verifies $\partial_t \U=0$, then it is also a solution of the equation (\ref{Bardina}) with initial datum $\U$. Thus, the term $\mathcal{R}$ solves the following equation: 
\begin{equation}
\partial_t \mathcal{R} - \nu \Delta  \mathcal{R} + \P \Big( ( (\vu \cdot \vec{\nabla})\vu - (\U \cdot \vec{\nabla})\U )_\alpha \Big)+ \beta \mathcal{R}=0,\quad \mathcal{R}(0,\cdot)=\vu_0-\U,  \end{equation} and consequently $\mathcal{R}(t,x)$ can be written as the integral form:
\begin{equation}
\mathcal{R}(t,x)=e^{\nu t \Delta}(\vu_0-\U) -\beta \int_{0}^{t}e^{\nu (t-s)} \mathcal{R}(s,x) ds - \int_{0}^{t}e^{\nu (t-s)} \P \Big( ( (\vu \cdot \vec{\nabla})\vu - (\U \cdot \vec{\nabla})\U )_\alpha \Big)(s,x) ds.     
\end{equation}
For all $t>0$, we write 
\begin{equation}\label{estim-R}
\begin{split}
\Vert \mathcal{R}(t,\cdot)\Vert_{L^{\infty}} &  \leq  \Vert  e^{\nu t \Delta}(\vu_0-\U) \Vert_{L^{\infty}}+ \beta \int_{0}^{t} \left\Vert  e^{\nu (t-s)} \mathcal{R}(s,\cdot)\right\Vert_{L^{\infty}}ds \\ 
&+ \int_{0}^{t} \left\Vert e^{\nu (t-s)} \P \Big( ( (\vu \cdot \vec{\nabla})\vu - (\U \cdot \vec{\nabla})\U )_\alpha \Big)(s,\cdot)\right\Vert_{L^{\infty}}ds=I_1(t)+I_2(t)+I_3(t),
\end{split}     
\end{equation}
where,  we will prove now that each term above verify $\ds{\Vert I_{i}(t)\Vert_{L^{\infty}}\leq C \, t^{-3/4}}$, with $t\gg 1$.  We mention here $C>0$ is a generically constant which may change from one estimate to other, but it does not depend on the time variable.\\  

For the first term $I_1(t)$, by well-known properties of the heat kernel $h_{\nu t}(x)$, and moreover, by the Young inequalities (with $1+1/\infty=1/2+1/2$), we directly have 
\begin{equation}\label{I1}
I_1(t) \leq  \Vert  h_{\nu t} \Vert_{L^2} \Vert \vu_0-\U \Vert_{L^2}\leq \frac{c_\nu}{t^{3/4}}\Vert \vu_0 - \U \Vert_{\mathcal{H}^{1}_{\alpha}} \leq C\, t^{-3/4}.     
\end{equation}
Thereafter, to study the second term $I_2(t)$ we write 
\begin{equation*}
I_2(t)\leq \beta \int_{0}^{t/2} \left\Vert  e^{\nu (t-s)} \mathcal{R}(s,\cdot)\right\Vert_{L^{\infty}}ds +   \beta \int_{t/2}^{t} \left\Vert  e^{\nu (t-s)} \mathcal{R}(s,\cdot)\right\Vert_{L^{\infty}}ds =I_{2,1}(t)+I_{2,2}(t). 
\end{equation*}
To estimate the term $I_{2,1}(t)$, using again  the Young inequalities, using the identity $\mathcal{R}=\vu-\U$, and moreover,  by the estimate (\ref{key-estim}) (where we have $\eta(\beta)<0)$, we can write 
\begin{equation*}
\begin{split}
I_{2,1}(t) & \leq  \beta \int_{0}^{t/2} \Vert h_{\nu(t-s)}\Vert_{L^2}\Vert \vu(s,\cdot)-\U\Vert_{L^2} ds  \leq c_{\beta,\nu} \int_{0}^{t/2} \frac{1}{(t-s)^{3/4}} \Vert \vu(s,\cdot)-\U\Vert_{L^2} ds  \\
& \leq \frac{c_{\beta,\nu}}{t^{3/4}} \int_{0}^{t/2} \Vert \vu(s,\cdot)-\U\Vert_{L^2} ds \leq \frac{c_{\beta,\nu}}{t^{3/4}} \int_{0}^{t/2} \Vert \vu(s,\cdot)-\U\Vert_{\mathcal{H}^{1}_{\alpha}} ds \leq  \frac{c_{\beta,\nu}}{t^{3/4}} \Vert \vu_0-\U\Vert_{\mathcal{H}^{1}_{\alpha}} \int_{0}^{t/2} e^{\eta(\beta)\, s} ds \leq  C\, t^{-3/4}. 
\end{split}   
\end{equation*}
These same facts also allow us to treat the term $I_{2,2}(t)$, where, for $t\gg 1$ we have 
\begin{equation*}
\begin{split}
I_{2,2}(t) & \leq  \beta \int_{t/2}^{t} \Vert h_{\nu(t-s)}\Vert_{L^2}\Vert \vu(s,\cdot)-\U \Vert_{L^2}ds \leq  c_{\beta,\nu} \int_{t/2}^{t} (t-s)^{-3/4} e^{\eta(\beta)s} ds \leq   c_{\beta,\nu}\,  e^{\eta(\beta) t /2}  \int_{t/2}^{t}(t-s)^{-3/4} ds \\
&\leq C e^{\eta(\beta) t/2} t^{1/4} \leq  C\, t^{-3/4}. 
\end{split}    
\end{equation*}
We thus have 
\begin{equation}\label{I2}
I_2(t)\leq C\, t^{-3/4}, \quad t\gg 1.     
\end{equation}

Finally, to estimate the first term $I_3(t)$, we recall first that the filtering operator $(\cdot)_\alpha$ given in (\ref{filtering}) can be also defined by a convolution product with  a kernel $K_\alpha$ \cite{Grafakos}. This kernel has good decaying properties in the spatial variable and moreover we have $\Vert K_{\alpha} \Vert_{L^p}<+\infty$, for $1\leq p <+\infty$. Thus, by the Young inequalities (with $1+1/\infty= 5/6+1/6$), and by the boundness of the Leray's projector in the Lebesgue spaces, we have
\begin{equation*}
\begin{split}
I_{3}&\leq \int_{0}^{t} \left\Vert \Big( e^{\nu (t-s)} \P \Big( (\vu \cdot \vec{\nabla})\vu - (\U \cdot \vec{\nabla})\U  \Big)(s,\cdot) \Big)_\alpha \right\Vert_{L^{\infty}} ds\\
&=  \int_{0}^{t} \left\Vert K_\alpha \ast \Big( e^{\nu (t-s)} \P \Big( (\vu \cdot \vec{\nabla})\vu - (\U \cdot \vec{\nabla})\U  \Big)(s,\cdot) \Big) \right\Vert_{L^{\infty}} ds  \\
&\leq  \Vert K_\alpha \Vert_{L^6} \int_{0}^{t}  \left\Vert  e^{\nu (t-s)} \P \Big(  (\vu \cdot \vec{\nabla})\vu - (\U \cdot \vec{\nabla})\U  \Big)(s,\cdot) \right\Vert_{L^{6/5}} ds  \\
&\leq c_\alpha \int_{0}^{t}  \left\Vert  e^{\nu (t-s)}  \Big( (\vu \cdot \vec{\nabla})\vu - (\U \cdot \vec{\nabla})\U  \Big)(s,\cdot) \right\Vert_{L^{6/5}} ds=(a).
\end{split}    
\end{equation*}
Here, as we have $div(\vu)=0$ and $div(\U)=0$, and moreover, applying again the Young inequalities (with $1+5/6=1+5/6$), we can write 
\begin{equation*}
\begin{split}
 (a)&=  c_\alpha \int_{0}^{t}  \left\Vert  e^{\nu (t-s)}  div \Big(\vu \otimes \vu - \U \otimes \U \Big)(s,\cdot) \right\Vert_{L^{6/5}} ds  \leq c_\alpha \int_{0}^{t} \Vert \vec{\nabla}h_{\nu(t-s)} \Vert_{L^{6/5}}\left\Vert \Big(\vu \otimes \vu - \U \otimes \U \Big)(s,\cdot) \right\Vert_{L^1} ds=(b).
\end{split}    
\end{equation*}
Recalling that $\mathcal{R}=\vu-\U$ then we have $\ds{\vu \otimes \vu - \U \otimes \U = \mathcal{R}\otimes \vu+ \U \otimes \mathcal{R}}$,  and moreover, by the H\"older inequalities (with $1=1/2+1/2$), we write 
\begin{equation*}
\begin{split}
(b)&\leq c_{\alpha,\nu}  \int_{0}^{t} (t-s)^{-3/4}\left\Vert \Big(\mathcal{R}\otimes \vu+ \U \otimes \mathcal{R}\Big)(s,\cdot) \right\Vert_{L^1} ds  \\
&\leq c_{\alpha,\nu}  \int_{0}^{t}(t-s)^{-3/4} \Vert \mathcal{R}(s,\cdot)\Vert_{L^2}\Big( \Vert \vu(s,\cdot)\Vert_{L^2}+\Vert \U \Vert_{L^2}\Big)ds\\
&\leq  c_{\alpha,\nu} \int_{0}^{t} (t-s)^{-3/4} \Vert (\vu-\U)(s,\cdot)\Vert_{L^2}\Big( \Vert \vu(s,\cdot)\Vert_{L^2}+\Vert \U \Vert_{L^2}\Big)ds \\
&\leq  c_{\alpha,\nu} \int_{0}^{t} (t-s)^{-3/4} \Vert (\vu-\U)(s,\cdot)\Vert_{\mathcal{H}^{1}_{\alpha}}\Big( \Vert \vu(s,\cdot)\Vert_{\mathcal{H}^{1}_{\alpha}}+\Vert \U \Vert_{\mathcal{H}^{1}_{\alpha}}\Big)ds=(c). 
\end{split}    
\end{equation*}
In this last expression, by (\ref{key-estim}) we have $\ds{\Vert (\vu-\U)(s,\cdot)\Vert_{\mathcal{H}^{1}_{\alpha}}\leq \Vert \vu_0 -\U \Vert_{\mathcal{H}^{1}_{\alpha}} e^{\eta(\beta)s}}$ (with $\eta(\beta)<0$). Moreover, by point $1)$ of Proposition \ref{Prop:time-control} we have $\ds{\Vert \vu(s,\cdot)\Vert^{2}_{\mathcal{H}^{1}_{\alpha}} \leq \Vert \vu_0  \Vert^{2}_{\mathcal{H}^{1}_{\alpha}}e^{-\beta s}+\frac{1}{\beta^2}\Vert \fe \Vert^{2}_{\mathcal{H}^{1}_{\alpha}}}$, hence we can write $\ds{\Vert \vu(s,\cdot)\Vert_{\mathcal{H}^{1}_{\alpha}} \leq \Vert \vu_0 \Vert_{\mathcal{H}^{1}_{\alpha}}+\Vert \fe \Vert_{\mathcal{H}^{1}_{\alpha}}}$. Thus, gathering these estimates, for $t\gg 1$ we obtain 
\begin{equation*}
\begin{split}
(c)&\leq C \int_{0}^{t}(t-s)^{-3/4}e^{\eta(\beta)s} ds= C\int_{0}^{t/2} (t-s)^{-3/4}e^{\eta(\beta)s} ds + C \int_{t/2}^{t}(t-s)^{-3/4}e^{\eta(\beta)s} ds  \leq C\, t^{-3/4}.  
\end{split}    
\end{equation*}
\begin{equation}\label{I3}
I_3(t)\leq C\, t^{-3/4}, \quad t\gg 1.     
\end{equation}

Once we have the estimates (\ref{I1}), (\ref{I2}) and (\ref{I3}), we get back to (\ref{estim-R}), hence we finally get the desired estimate (\ref{time-decay1}). Proposition \ref{Prop:profile} is proven. \finpv 

\section{The damped Navier-Stokes-Bardina's model without external force}\label{Sec:ForceZero}

\subsection*{Proof of Proposition \ref{Prop:force-zero}}

When $\fe=0$, we observe first that $\U=0$ is a solution of the stationary problem (\ref{Bardina-Stat}), and moreover, by point $2)$ of Theorem \ref{Th:Stability} is the unique one since here we have $\eta(\beta)=-\beta <0$. We thus have $\mathcal{A}_{0}=\{ 0 \}$.  \\

On the other hand, for any initial datum $\vu_0 \in \mathcal{H}^{1}_{\alpha}(\Rt)$, let $\vu(t,\cdot)$ be the solution of the equation (\ref{Bardina}) given by Theorem \ref{Th-WP}. Then,   by (\ref{key-estim}) we have  $\ds{\Vert \vu(t,\cdot)\Vert_{\mathcal{H}^{1}_{\alpha}} \leq \Vert \vu_0 \Vert_{\mathcal{H}^{1}_{\alpha}} e^{-\beta\, t}}$, hence we can write  $\ds{\Vert \vu(t,\cdot)\Vert_{L^2} \leq \Vert \vu_0 \Vert_{\mathcal{H}^{1}_{\alpha}} e^{-\beta\, t}}$.  Moreover, by Proposition \ref{Prop:profile} (with $\U=0$) we also have $\Vert \vu(t,\cdot)\Vert_{L^{\infty}} \leq C\, t^{-3/4} \leq C$, with $t\gg 1$. Thus, for $2\leq p <+\infty$, by the interpolation inequalities we write 
\begin{equation*}
\Vert \vu(t,\cdot)\Vert_{L^{p}} \leq c \Vert \vu(t,\cdot)\Vert^{2/p}_{L^2}\Vert \vu(t,\cdot)\Vert^{1-2/p}_{L^{\infty}} \leq C(p,\vu_0)\, e^{- \frac{2\beta}{p}\, t}.      
\end{equation*}
Proposition \ref{Prop:force-zero} is proven. \finpv

\quad\\[5mm]
\begin{flushright}
\begin{minipage}[r]{130mm}
	\textbf{{\large Manuel Fernando Cortez}} (manuel.cortez@epn.edu.ec) \\ 
	\\
	Departamento de Matem\'aticas, Escuela Politécnica Nacional,  Ladron de Guevera E11-253\\
	Quito - Ecuador.\\ 
	\\
	\\
	\\
	\textbf{\large {Oscar Jarr\'in}} (oscar.jarrin@udla.edu.ec)\\
	\\
	Dirección General  de Vinculaci\'on e Investigación  (DGVI)\\
	Universidad de las Américas \\ 
	Calle José Queri s/n y Av. Granados. Bloque 7, Tercer Piso\\
Quito – Ecuador. 
\end{minipage}
\end{flushright}	 


\begin{thebibliography}{40}
\bibitem{Adams} N. A. Adams \& S. Stolz. \emph{ Deconvolution methods for subgrid-scale approximation in large eddy simulation}. Modern Simulation Strategies for Turbulent Flow, R.T. Edwards (2001).
\bibitem{Babin} A. Babin \& M. Vishik. \emph{Attractors of evolution partial differential equations
and estimates of their dimension}. Uspekhi Mat. Nauk 38:3, 133–187 (1983).
English transl. Russian Math. Surveys. 38 (1983)

\bibitem{Babin2}  A. Babin \& M. Vishik. \emph{Atracttors for evolution equations}. Studies in Mathematics and its Applications, 25. North-Holland Publishing Co., Amsterdam (1992).   


\bibitem{Bardina1} J. Bardina, J. Coakley  \&  P. Huang. \emph{Turbulence Modeling Validation, Testing, and Development}.  NASA Technical Memorandum (1997). 
\bibitem{Bardina} J. Bardina, J. H. Ferziger, \& W. C. Reynolds. \emph{Improved subgrid scale models for large eddy simulation}. AIAA paper, 80:80–1357, (1980).
\bibitem{Berselli} L. C. Berselli \& R. R. Lewandowski. \emph{On the Bardina’s Model in the Whole Space}. Journal of Mathematical Fluid Mechanics volume 20, pages1335–1351(2018).
\bibitem{Cao}  Y. Cao, E. M. Lunasin \& Edriss S. Titi. \emph{Global well-posedness of the three-dimensional viscous and inviscid simplified Bardina turbulence models}.  Commun. Math. Sci. 4(4): 823-848 (2006). 
\bibitem{Brieve} S. De Briève, F. Genoud \& S. Rota Nodari. \emph{Orbital stability: analysis meets geometry}. Nonlinear Optical and Atomic Systems, 2146, pp 147-273, Lectures Notes in Mathematics (2015). 

\bibitem{Chamorro} D. Chamorro, O. Jarr\'in \& P.G. Lemarié-Rieusset. \emph{On the Kolmogorov Dissipation Law in a Damped Navier–Stokes Equation}. Journal of Dynamics and Differential Equations (2020).
\bibitem{Chow} F.K. Chow,  S.F. De Wekker \& B. Synder. \emph{Mountain Weather Research and Forecasting: Recent Progress and Current Challenges}. Springer, Berlin (2013).
\bibitem{Chepyzhov} V. Chepyzhov \& A. Ilyin. \emph{On the fractal dimension of invariant sets: applications to Navier-Stokes equations}. Discrete Contin. Dyn. Syst. 10 no. 1-2, 117–135 (2004).

\bibitem{Constantin1} P. Constantin \& C. Foias. \emph{Global Lyapunov exponents, Kaplan–Yorke formulas and the dimension of the attractors for the 2D Navier–Stokes equations}.
Comm. Pure Appl. Math. 38: 1–27 (1985).
\bibitem{Constantin2} P. Constantin \& F. Ramos. \emph{Inviscid limit for damped and driven incompressible Navier-Stokes equations in $\mathbb{R}^2$}. Comm. Math. Phys. 275: 529–551 (2007).
\bibitem{Grafakos} L. Grafakos. \emph{Modern Fourier Analysis}. Second Edition. Springer Series 250 (2008). \bibitem{Ilyn1} A. Ilyin, E.M. Lunasin \& E.S. Titi. \emph{A modified-Leray-$\alpha$ subgrid scale model of turbulence}. Nonlinearity, 19(4):879–897 (2006).
\bibitem{Ilyn} A. Ilyn, K. Patni \& S. Zelik. \emph{Upper bounds for the attractor dimension of damped Navier-Stokes equations in $\R^2$}. Discrete and Continuous Dynamical Systems 26(4):2085-2102 (2016).
\bibitem{Jarrin} O. Jarr\'in. \emph{Deterministic descriptions of turbulence in the Navier-Stokes equations}. Ph.D. thesis at \emph{Université Paris-Saclay}, hal-01821762  (2018). 
\bibitem{Layton2} W. Layton \& R. Lewandowski.  \emph{A simple and stable scale-similarity model for large eddy simulation: energy balance and existence of weak solutions}. Appl. Math. Lett., 16(8):1205–1209, (2003).
\bibitem{Layton1} W. Layton \& R. Lewandowski. \emph{ On a well-posed turbulence model}.  Dicrete and Continuous Dyn. Sys. B, 6, (2006).
111-128.
\bibitem{Liua} H. Liua \& H. Gao. \emph{Decay of solutions for the 3D Navier–Stokes equations with damping}. Applied Mathematics Letters 68: 48–54 (2017).
\bibitem{PGLR1} P.G. Lemarié-Rieusset. \emph{The Navier--Stokes Problem in the 21st Century}. CRC Press, Chapman \& Hall Book (2016). 
\bibitem{Lieb} E.H. Lieb. \emph{Lieb-Thirring Inequalities}.  	arXiv:math-ph/0003039 (2020);
\bibitem{Pedlosky} J. Pedlosky, Geophysical Fluid Dynamics, Springer, New York, (1979).  
\bibitem{Raugel} G. Raugel. \emph{Global Attractors in Partial Differential Equations}. Lectures notes, CNRS et Université Paris-Sud, Analyse Numérique et EDP, UMR 8628 (2006). 
\bibitem{Temam} R. Temam. \emph{Infinite Dimensional Dynamical Systems in Mechanics and Physics, 2nd ed}. Springer-Verlag, New York, (1997).
\bibitem{Temam1} R. Temam. \emph{Navier-Stokes equations.  Theory and numerical analysis. } Reprint of the 1984 edition. AMS Chelsea Publishing, Providence, RI, (2001).
\end{thebibliography}
\end{document}